\documentclass [10pt,a4paper]{article}
\usepackage{amssymb}
\usepackage{enumerate,verbatim}
\usepackage{amsmath}
\usepackage{amsfonts,mathrsfs}
\usepackage{amsthm,wasysym}
\usepackage{epsfig,lpic,tikz}
\usepackage[applemac]{inputenc}
\usepackage[english]{babel}

\usepackage[width=0.9\textwidth]{caption}

\usepackage{pdfsync}
\usepackage{xcolor,hyperref}
\usepackage[T1]{fontenc}
\usepackage{lpic}
\usepackage{mathtools}

\hypersetup{
    colorlinks=true,
    linkcolor=blue,
    filecolor=magenta,      
    urlcolor=blue,
    pdftitle={Overleaf Example},
    pdfpagemode=FullScreen,
    }

\definecolor{green}{rgb}{0,0.5,0}

\usepackage{marvosym}

\newcounter{teoremaganso}
\newcounter{appendix}

\newcounter{coryganso}

\renewenvironment{abstract}{\small\quotation\noindent
 {\bfseries \abstractname .}}{\endquotation \par}

\newenvironment{trivlistparm}[2]{\trivlist
\item[{#1 #2}]}{\endtrivlist}
\newenvironment{proclama}[1]{\trivlistparm{\bfseries}{#1}\itshape}{\endtrivlistparm}
\newenvironment{prooftext}[1]{\trivlistparm{\bfseries}{#1}}{\Qed\endtrivlistparm}
\newenvironment{prova}{\trivlistparm{\bfseries}{Proof.}}{\Qed\endtrivlistparm}

\catcode`\@=11

\def\resetthefootnote{\renewcommand{\thefootnote}{\@arabic\c@footnote} }
\def\@principiremex#1{\trivlist
 \item[\hskip \labelsep{\bfseries #1\ \thetheo.}]\ignorespaces}
\def\opar@principiremex#1[#2]{\trivlist
 \item[\hskip \labelsep{\bfseries #1\ \thetheo\ (#2).}]\ignorespaces}

\newcommand{\newTHEOremrom}[2]{\newenvironment{#1}{\refstepcounter{theo}\@ifnextchar[{\opar@principiremex{#2}}
{\@principiremex{#2}}}{\qedB\endtrivlist}} \catcode`\@=12
\DeclareMathSymbol{\square}{\mathord}{AMSa}{"03}
\newcommand{\qedB}{\nopagebreak\hspace*{\fill}$\square$\par}
\newcommand{\Qed}{\nopagebreak\hspace*{\fill}{\vrule width6pt height6pt depth0pt}\par}

\newtheorem {theo} {Theorem} [section]
\newtheorem {prop} [theo] {Proposition}
\newtheorem {cory} [theo] {Corollary}
\newtheorem {lem} [theo] {Lemma}
\newtheorem {bigtheo} [teoremaganso] {Theorem}

\newTHEOremrom {defi} {Definition}
\newTHEOremrom {obs} {Remark}
\newTHEOremrom {ex} {Example}

\newcommand{\refc}[1]{\mbox{$(\ref{#1})$}}
\newcommand{\secc}[1]{Section~\ref{#1}}
\newcommand{\apc}[1]{Appendix~\ref{#1}}
\newcommand{\teoc}[1]{Theorem~\ref{#1}}
\newcommand{\propc}[1]{Proposition~\ref{#1}}
\newcommand{\coryc}[1]{Corollary~\ref{#1}}
\newcommand{\lemc}[1]{Lemma~\ref{#1}}
\newcommand{\defic}[1]{Definition~\ref{#1}}
\newcommand{\obsc}[1]{Remark~\ref{#1}}

\newcommand{\figc}[1]{Figure~\ref{#1}}

\newcommand{\N}{\ensuremath{\mathbb{N}}}
\newcommand{\Z}{\ensuremath{\mathbb{Z}}}
\newcommand{\R}{\ensuremath{\mathbb{R}}}

\newcommand{\T}{\boldsymbol{T}}

\newcommand{\C}{\ensuremath{\mathbb{C}}}

\newcommand{\F}{\ensuremath{\mathcal{F}}}

\newcommand{\RP}{\ensuremath{\mathbb{RP}}}
\newcommand{\no}{{\ensuremath{{\hat\mu}}}}
\newcommand{\np}{{\ensuremath{{\nu}}}}

\newcommand{\Sc}{\ensuremath{\mathbb{S}}}
\newcommand{\Dc}{\mathbb D}
\newcommand{\cc}{\ensuremath{\mathscr{C}}}

\newcommand{\al}{\ensuremath{a}}
\newcommand{\be}{\ensuremath{b}}
\newcommand{\ga}{\ensuremath{\gamma}}

\def\map#1#2#3{\mbox{${#1}\!:{#2}\longrightarrow{#3}$}}

\newcommand{\sist}[2]{
  \left\{\!
   \begin{array}{l}
    \dot x=#1 \\[2pt] \dot y=#2
   \end{array}
  \right.
}

\newcommand{\op}{\ensuremath{\mbox{\rm o}}}
\newcommand{\PA}{\mathscr{P}}
\newcommand{\out}{\Pi}
\newcommand{\mf}[1]{\mathfrak{#1}}

\newcommand{\mr}[1]{\mathrm{#1}}

\newcommand{\blue}[1]{{\color{blue}#1}}

\newcommand{\dsp}{\displaystyle}
\newcommand{\gorro}{\hat}

\title{\textbf{The criticality of reversible quadratic centers \\ at the outer boundary of its period annulus}
\footnotetext{This work has been partially funded by the Ministry of Science, Innovation and Universities of Spain through the grants PGC2018-095998-B-I00 and MTM2017-86795-C3-2-P and by the Agency for Management of University and Research Grants of Catalonia through the grants 2017SGR1725 and 2017SGR1617.
}
}

\author{D. Mar\'{\i}n and J. Villadelprat
\\*[.1truecm]
{\small \textsl{Departament de Matem{\`a}tiques, Edifici Cc,
Universitat Aut{\`o}noma de Barcelona,}}\\*[-.05truecm]
{\small\textsl{08193 Cerdanyola del Vall\`es (Barcelona), Spain}}
\\*[-.05truecm]
{\small \textsl{Centre de Recerca Matem\`atica, Edifici Cc, Campus de Bellaterra,}}\\*[-.05truecm]
{\small \textsl{08193 Cerdanyola del Vall\`es (Barcelona), Spain}}
\\*[.1truecm]
{\small \textsl{Departament d'Enginyeria Inform{\`a}tica i Matem{\`a}tiques, ETSE,}}
\\*[-.05truecm]
{\small \textsl{Universitat Rovira i Virgili, 43007 Tarragona, Spain}}}
                  
\date{\vspace{-5ex}}

\begin{document}
\maketitle

\begin{abstract}
This paper deals with the period function of the reversible quadratic centers 
\begin{equation*}
X_{\np}=-y(1-x)\partial_x+(x+Dx^2+Fy^2)\partial_y,
\end{equation*}  
where $\np=(D,F)\in\R^2.$ Compactifying the vector field to $\Sc^2$, the boundary of the period annulus has two connected components, the center itself and a polycycle. We call them the inner and outer boundary of the period annulus, respectively. We are interested in the bifurcation of critical periodic orbits from the polycycle $\out_\np$ at the outer boundary. A critical period is an isolated critical point of the period function. The criticality of the period function at the outer boundary is the maximal number of critical periodic orbits of $X_\np$ that tend to $\out_{\np_0}$ in the Hausdorff sense as $\np\to\np_0.$ This notion is akin to the cyclicity in Hilbert's 16th Problem. Our main result (Theorem A) shows that the criticality at the outer boundary is at most 2 for all $\np=(D,F)\in\R^2$ outside the segments $\{-1\}\times [0,1]$ and $\{0\}\times [0,2]$. With regard to the bifurcation from the inner boundary, Chicone and Jacobs proved in their seminal paper on the issue that the upper bound is 2 for all $\np\in\R^2.$ In this paper the techniques are different because, while the period function extends analytically to the center, it has no smooth extension to the polycycle. We show that the period function has an asymptotic expansion near the polycycle with the remainder being uniformly flat with respect to~$\np$ and where the principal part is given in a monomial scale containing a deformation of the logarithm, the so-called \'Ecalle-Roussarie compensator. More precisely, Theorem~A follows by obtaining the asymptotic expansion to fourth order and computing its coefficients, which are not polynomial in~$\np$ but transcendental. Theorem~A covers two of the four quadratic isochrones, which are the most delicate parameters to study because its period function is constant. The criticality at the inner boundary in the isochronous case is bounded by the number of generators of the ideal of all the period constants but there is no such approach for the criticality at the outer boundary. A crucial point to study it in the isochronous case is that the flatness of the remainder in the asymptotic expansion is preserved after the derivation with respect to parameters. We think that this constitutes a novelty that is of particular interest also in the study of similar problems for limit cycles in the context of Hilbert's 16th Problem. 
Theorem~A also reinforces the validity of a long standing conjecture by Chicone claiming that the quadratic centers have at most two critical periodic orbits. A less ambitious goal is to prove the existence of a uniform upper bound for the number of critical periodic orbits in the family of quadratic centers. By a compactness argument this would follow if one can prove that the criticality of the period function at the outer boundary of any quadratic center is finite. Theorem~A leaves us very close to this existential result.
\end{abstract}

\tableofcontents

\section{Introduction and main results}

A singular point~$p\in\R^2$ of a planar differential system
\[
 \sist{f(x,y),}{g(x,y),}
\] 
is a \emph{center} if it has a punctured neighbourhood that consists entirely of periodic orbits surrounding~$p$. 
The \emph{period annulus} of the center is the  largest punctured neighbourhood with this property and we denote it by~$\PA$. The period annulus is an open subset of~$\R^2$ that may be unbounded. For this reason we embed~$\PA$ in~$\RP^2$ and, abusing notation, we denote the boundary of the resulting set by $\partial\PA$. Clearly the center $p$ belongs to~$\partial\PA$ and in what follows we call it the \emph{inner boundary} of 
the period annulus. We also define the \emph{outer boundary} of the period annulus 
to be $\Pi\!:=\partial\PA\setminus\{p\}$, which is a nonempty compact subset of $\RP^2.$  The subject of our study is the \emph{period function} of the center, that assigns to each periodic 
orbit in~$\PA$ its period. Since the period function is defined on the set of periodic orbits in~$\PA,$ in order to study its qualitative properties we need to parametrize
this set. This can be done by taking a transverse section
to the vector field $X=f(x,y)\partial_x+g(x,y)\partial_y$ on~$\PA$, for instance an orbit of the orthogonal vector
field~$X^\bot$. To fix ideas let us suppose that $\{\gamma_s\}_{s\in (0,1)}$ is such a
parametrization where $s\approx 0$ corresponds to the periodic orbits near~$p$ and $s\approx 1$ to the ones near~$\out$. Then the map \map{P}{(0,1)}{(0,+\infty)} defined by $P(s)\!:=\!\{\mbox{period of
$\gamma_s$}\}$ provides the qualitative
properties of the period function that we are concerned with and one can readily show by using the Implicit Function Theorem that it is as smooth as $X$. It is also well-known that if $X$ is analytic and the center~$p$ is non-degenerate then $P$ extends analytically to $s=0.$ Let us advance that, on the contrary, $P$ does not extend smoothly to $s=1.$ The \emph{critical periods} are the isolated critical points of~$P$, i.e. $\hat s\in (0,1)$
such that $P'(\hat s)=0$ and $P'(s)\neq 0$ if $0<|s-\hat s|<\varepsilon.$
In this case, more geometrically, we shall say
that~$\gamma_{\hat s}$ is a \emph{critical periodic orbit} of $X$. One can easily see that the property of being a critical periodic orbit does not 
depend on the particular parametrization of the set of periodic orbits used, see \obsc{local_ind}. The study of the critical periodic orbits is another issue arising from the famous Hilbert's 16th Problem and it has strong parallelisms with the research on limit cycles, from both the conceptual and technical point of views. In this regard we can mention for instance that the isochronicity problem (i.e., to decide whether a center has a constant period function) is the counterpart of the center-focus problem. The renowned conjecture claiming that a quadratic differential system can have at most four limit cycles has also an analogue in the context of the period function and it was posed by C. Chicone~\cite{Chi_mathscinet}. More specifically this conjecture asserts that if a quadratic center has some critical periodic orbit then by an affine transformation and a constant rescaling of time it can be brought to Loud normal form
\begin{equation}\label{BDF}
 \sist{-y+Bxy,}{x+Dx^2+Fy^2,}
\end{equation}
and that this center has at most two critical periodic orbits for any $(B,D,F)\in\R^3$. In fact there is much analytic evidence that this conjecture is true (see \cite{CopGav,Jordi2,Zhao} for instance).  

The problems that we are interested in take place when the vector field $X$ depends on parameters. To fix notation, let $U$ be an open subset of $\R^N$ and consider a family of planar vector fields $\{X_\mu,\mu\in U\}$ such that each $X_\mu$ has a center $p_\mu$ with period annulus $\PA_\mu$. Let us denote the period function of the center~$p_\mu$ by $P(\,\cdot\,;\mu)$ and observe that, given some $\mu_0\in U$, the number of critical periodic orbits of~$X_\mu$ can vary as we perturb $\mu\approx\mu_0.$ Under some regularity assumptions on the dependence of $\PA_\mu$ with respect to $\mu$ it can be proved (see \lemc{diag}) that the emergence/disappearance of critical periodic orbits can only occur from three different places: 
\begin{enumerate}[$(a)$]
\item Bifurcations at the inner boundary of the period annulus (i.e., the center $p_\mu$).
\item Bifurcations at the outer boundary of the period annulus (i.e., the polycycle $\out_\mu$).
\item Bifurcations at the interior of the period annulus $\PA_\mu.$
\end{enumerate}
Chicone and Jacobs give in their seminal paper \cite{Chicone} a complete description of the bifurcations from the inner boundary for the whole family of quadratic centers. In this case the parameter $\mu$ are the coefficients of the vector field and since the center is non-degenerate $P(s;\mu)$ extends analytically to $s=0$, so that one can consider its Taylor series $P(s;\mu)=\sum_{i=0}^\infty a_i(\mu)s^i$ at $s=0$, whose coefficients $a_i$ belong to the polynomial ring~$\R[\mu]$. On account of this the result about the bifurcations from the isochronous centers (see \cite[Theorem 2.2]{Chicone}), which are the most difficult ones to study, follows by analyzing the ideal $(a_1,a_2,\ldots)$ of all Taylor coefficients exactly as N. Bautin does in \cite{Bautin} to study the bifurcations of limit cycles from the quadratic centers. In the present paper we resume our study of the bifurcations from the outer boundary that we initiated in \cite{MMV03,MMV2}. Let us recall that the differential system \refc{BDF} has no critical periodic orbits if $B=0,$ see \cite[Theorem 1]{GGV}. By means of a rescaling the case $B\neq 0$ can be brought to $B=1$, i.e., 
\begin{equation}\label{sist_loud}
X_{\np}\!:=-y(1-x)\partial_x+(x+Dx^2+Fy^2)\partial_y\text{ with $\np\!:=(D,F).$}
\end{equation}   
Here we already adopt the parameter notation that we shall use throughout the paper, which is devoted to the bifurcation of critical periodic orbits from the outer boundary in the family $\{X_\np,\np\in\R^2\}.$ Since each vector field~$X_\np$ is polynomial we can consider its Poincaré compactification $p(X_\np)$, see \cite[\S 5]{ADL}, which is an analytic vector field on the sphere~$\Sc^2$ topologically equivalent to~$X_\np.$ The outer boundary~$\out_\np$ becomes then a polycycle of $p(X_\np)$ that can be studied using local charts of~$\Sc^2$, but even so the period function $P(s;\np)$ cannot be smoothly extended to $s=1.$ For the family under consideration we show that $P(s;\np)$ has an asymptotic expansion at $s=1$ with the remainder being uniformly flat with respect to~$\np$ and where the principal part is given in a monomial scale containing a deformation of the logarithm, the so-called \'Ecalle-Roussarie compensator. Our main theorem follows by obtaining the asymptotic expansion to fourth order and computing its coefficients, which are not polynomial in~$\np$ but transcendental (more concretely, they are hypergeometric functions). To this end we strongly rely on the tools that we develop in our recent papers \cite{MV20a,MV20b,MV21}. 
The results that we obtain in the present paper can be viewed conceptually as the analogue for the outer boundary of the work carried out by Chicone and Jacobs in~\cite{Chicone} on the bifurcation of critical periodic orbits from the inner boundary of the quadratic centers. That being said, the proofs of the results on the outer boundary are technically tougher than the ones on the inner boundary because~$\out_\np$ is a polycycle 
and the period function $P(s;\np)$ cannot be analytically extended there. By way of example, to determine the parameters that vanish simultaneously two coefficients in the asymptotic expansion at $s=1$ takes 5 pages of computations dealing with a hypergeometric function (see \apc{apB}), whereas the same problem for the Taylor series at $s=0$ can be solved readily by taking resultants because the coefficients are polynomials. 

\begin{figure}[t]
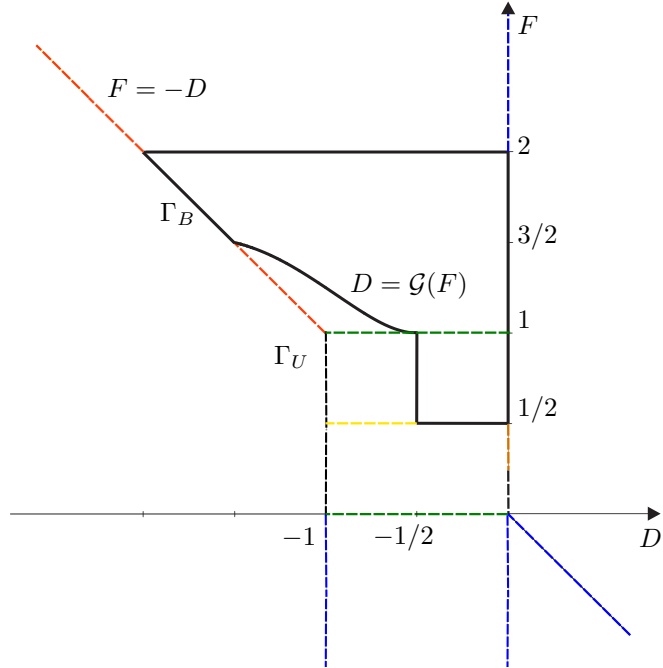

 \centering
 \begin{lpic}[l(0mm),r(0mm),t(0mm),b(5mm)]{dib2b(1)}
   \lbl[l]{13,76.5;$F=-D$}
   \lbl[l]{67,85;$F$}
   \lbl[l]{83,17;$D$}
   \lbl[l]{67,34;$1/2$}
   \lbl[l]{67,46;$1$}  
   \lbl[l]{67,57;$3/2$}  
   \lbl[l]{67,69;$2$}  
   \lbl[l]{48,17;$-1/2$}  
   \lbl[l]{36,17;$-1$}
   \lbl[l]{20,60;$\Gamma_B$}  
   \lbl[l]{35,41;$\Gamma_U$}         
   \lbl[l]{45.0,50.5;$D=\mathcal G(F)$}
 \end{lpic}     
\caption{The thick (closed) curve $\Gamma_B$ consists of local bifurcation values of the period function 
              at the outer boundary according to \cite{MMV2}, where 
              the curve that joins $\left(-\frac{3}{2},\frac{3}{2}\right)$ and 
              $\left(-\frac{1}{2},1\right)$ is the graphic of an analytic function $D=\mathcal G(F)$, see \obsc{fun}.
              The dotted lines $\Gamma_U$ correspond to  
              parameters that remained unspecified in that paper and we colour the subsequent improvements  
              obtained in \cite{TopaV,MMV20,MMSV,MSV,MV,Mariana,Jordi1}. 
              The parameters outside $\Gamma_B\cup\Gamma_U$ are 
              local regular values by the result in \cite{MMV2}.
              }\label{diagrama}
 \end{figure}

In this paper we use the notion of \emph{criticality} of the period function at the outer boundary which, roughly speaking, is the number of critical periodic orbits that can emerge or disappear from $\out_\np$ as we perturb $\np$ slightly. It is defined in exactly the same way as the notion of \emph{cyclicity} of a limit periodic set, which is used to study the bifurcation of limit cycles in the context of Hilbert's 16th Problem, see~\cite{Roussarie} for instance. Before giving its precise definition, and the statement of our main contribution, we enumerate the previous results about the bifurcation of critical periodic orbits from the outer boundary~$\out_\np$ for the family $\{X_\np,\np\in\R^2\}.$ In this regard we stress that these results are given according to the dichotomy between local regular value and local bifurcation value (of the period function at the outer boundary) that we introduce in our early paper~\cite{MMV2}. This notion (see \defic{def2}) enables to obtain a structure theorem for the bifurcation diagram of the period function in its full domain (see \lemc{diag}), but it has the inconvenience of not being so quantitative and geometric as the criticality. In order to simplify the exposition for the moment we can think that $\np_0\in\R^2$ is a local regular value 
if and only if the criticality of the period function at~$\out_{\np_0}$ is zero (i.e., no critical periodic orbit bifurcates from $\out_{\np_0}$ as we perturb $\np\approx\np_0).$ 
That said, let~$\Gamma_U$ be the union of dotted straight
lines in \figc{diagrama}, whatever its colour is. Consider also the thick curve~$\Gamma_B.$ (Here the
subscripts~$B$ and~$U$ stand for bifurcation and unspecified respectively.)
 Then according to \cite[Theorem A]{MMV2} the open set $\R^2\setminus\left(\Gamma_B\cup\Gamma_U\right)$ corresponds to local regular values and~$\Gamma_B$ consists of local bifurcation values (of the period function at the outer boundary). In that paper we also conjecture that any parameter in $\Gamma_U$ is regular, except for the segment $\{0\}\!\times\!\left[0,\frac{1}{2}\right]$ in the vertical axis, that should consist of bifurcation values. Since the formulation of this conjecture there has been some progress in the study of the parameters in~$\Gamma_U$: 

\begin{itemize}

\item From the results in \cite{TopaV,Jordi1} it follows that the parameters in blue are regular. In these papers the 
         authors determine a region $M$ in the parameter plane for which the corresponding center has a globally
         monotonous period function (i.e., it has no critical periodic orbits). 
         The parameters that we draw in blue are inside the interior of~$M$, which prevents 
         the bifurcation of critical periodic orbits.
         
\item Along the straight line $F=-D$ there is a breaking of a heteroclinic connection between two
         hyperbolic saddles at the outer boundary. From the results in \cite{MMV20} it follows that the
         parameters in red are regular.
         
\item Along the two segments in green it occurs a saddle-node bifurcation at the outer boundary of the period annulus. 
         An asymptotic expansion of the Dulac time of this type of unfolding is obtained
         in~\cite{MMSV} and as an application 
         it is proved that the parameters in the segment $(-1,0)\times\{1\}$, with the exception of $(-\frac{1}{2},1)$, are 
         local regular values. A subsequent refinement of this approach shows in~\cite{MSV} that the segment  
         $(-1,0)\times\{0\}$ also consists of local regular values.
         
\item By \cite[Theorem B]{MV} the parameters in brown, more precisely the segment 
          $\{0\}\!\times\!\left[\frac{1}{4},\frac{1}{2}\right]$, are local bifurcation values of the period function at the outer
          boundary. 

\item Along the segment $(-1,0)\times\{\frac{1}{2}\}$ there is a resonant saddle at $\out_\np$ and
        the parameters in yellow are local 
        regular values at the outer boundary of the period annulus according to \cite[Corollary~B]{Mariana}.

\end{itemize} 
As we already explained, these results are addressed to solve the dichotomy between local regular value and local bifurcation value (of the period function at the outer boundary). Beyond this dichotomy a challenging problem is the computation of the exact number of critical periodic orbits that can bifurcate from the outer boundary, which constitutes the counterpart of the result by Chicone and Jacobs~\cite{Chicone} about the bifurcation from the inner boundary. The following is the precise definition of the number that we aim to compute for the quadratic centers, where 
$d_{H}$ stands for the Hausdorff distance between compact sets of~$\RP^2.$

\begin{defi}\label{criticality}
Consider a $\mathscr C^{\infty}$ family $\{X_{\mu},\mu\in U\}$ of planar vector fields with a center and fix some $\mu_0\in U$. Suppose that the outer boundary of the period annulus varies continuously at $\mu_0\in U$, meaning that
$d_{H}(\out_{\mu},\out_{\mu_0})$ tends to zero as $\mu\to\mu_0.$ Then, setting 
\[
 N(\delta,\varepsilon)=\sup\left\{\text{\# critical periodic orbits $\gamma$ of $X_{\mu}$ in $\PA_{\mu}$ with $d_{H}(\gamma,\out_{\mu_0})\leqslant\varepsilon$ and $\|\mu-\mu_0\|\leqslant\delta$}\right\},
\]
the \emph{criticality} of $(\out_{\mu_0},X_{\mu_0})$ w.r.t. the deformation $X_{\mu}$ is $\mathrm{Crit}\bigl((\out_{\mu_0},X_{\mu_0}),X_{\mu}\bigr)\!:=\inf_{\delta,\varepsilon}N(\delta,\varepsilon).$
\end{defi}
We stress that in this definition the vector field $X_\mu$ is not required to be polynomial but $\mathscr C^{\infty}$. This is so because in order to define the outer boundary $\Pi_\mu$ of the period annulus $\PA_\mu$ of $X_\mu$ we do not compactify the vector field but only the set~$\PA_\mu$ and to this end there is no need that $X_\mu$ is polynomial. Certainly $\mathrm{Crit}\bigl((\out_{\mu_0},X_{\mu_0}),X_{\mu}\bigr)$ may be infinite but, if it is not, then it gives the maximal number of critical periodic orbits of $X_{\mu}$ that tend to $\out_{\mu_0}$ in the Hausdorff sense as $\mu\to\mu_0.$ Related with this issue we point out that the contour of the period annulus~$\PA_{\mu_0}$ may change for $\mu\approx\mu_0.$ The assumption that the period annulus varies continuously ensures that this change does not occur abruptly. In this regard note that $X_{\mu}=-y\partial _x+(x+\mu x^3+x^5)\partial_y$, with $\mu\in\R,$ is a polynomial family of vector fields with a center at the origin for which the outer boundary does not vary continuously at $\mu=2.$ This is so because the period annulus~$\PA_{\mu}$ is the whole plane for $\mu<2$, whereas it is bounded for $\mu=2$ (see~\cite{ManVil06} for details). In this example 
$\mathrm{Crit}\bigl((\out_{\mu_0},X_{\mu_0}),X_{\mu}\bigr)$, as introduced in \defic{criticality}, does not give the number of critical periodic orbits bifurcating from $\Pi_\mu$ as $\mu\to\mu_0.$
Let us mention that this assumption is also required in \cite{MRV16,MRV18}, where the authors obtain several results addressed to bound the criticality at the outer boundary of families of vector fields of potential type, i.e., $-y\partial_x+V'(x)\partial_y$.

Let us remark at this point that if Chicone's conjecture about the number of critical periodic orbits of the quadratic centers is true then $\mathrm{Crit}\bigl((\out_{\np_0},X_{\np_0}),X_{\np}\bigr)\leqslant 2$ for all $\np\in\R^2,$ see \refc{sist_loud}. In this paper, by applying our recent results from \cite{MV20a,MV20b,MV21}, we prove the following:

\begin{bigtheo}\label{Loud}
Let $\{X_{\np},\np\in\R^2\}$ be the family of quadratic vector fields given in~\refc{sist_loud} and consider the period function of the center at the origin. Then the following assertions hold:
\begin{enumerate}[$(a)$]
\item $\mathrm{Crit}\bigl((\out_{\np_0},X_{\np_0}),X_{\np}\bigr)=0$ 
        for $\np_0\notin\Gamma_B\cup\big\{D=-1,F\in[0,1]\big\}\cup\big\{D=0, F\in[0,\frac{1}{2}]\big\}$.

\item $\mathrm{Crit}\bigl((\out_{\np_0},X_{\np_0}),X_{\np}\bigr)=1$ 
        for $\np_0\in\Gamma_B\setminus\left(\left\{D=0\right\}\cup\left\{(-2,2),
       (\mathcal G(\frac{4}{3}),\frac{4}{3})\right\}\right)$. 
       
\item $\mathrm{Crit}\bigl((\out_{\np_0},X_{\np_0}),X_{\np}\bigr)\geqslant 1$ 
         for $\np_0 \in\big\{D=0,F\in[\frac{1}{4},2]\big\}$.
         
\item $\mathrm{Crit}\bigl((\out_{\np_0},X_{\np_0}),X_{\np}\bigr)=2$ 
         for $\np_0\in\big\{(-2,2), (\mathcal G(\frac{4}{3}),\frac{4}{3})\big\}$.   
         
\item There is a $\cc^1$ curve arriving 
         at $\nu= (\mathcal G(\frac{4}{3}),\frac{4}{3})$ tangent to $\Gamma_B$ and
         there is a $\cc^0$ curve with a exponential flat contact with $\{F=2\}$ at $\nu=(-2,2),$
         consisting both of local bifurcation values of the period function at the interior.       
\end{enumerate}
\end{bigtheo}

There are some papers containing results related with assertion $(b)$ in \teoc{Loud} to be referred. Thus, by \cite[Theorem A]{RV}, $\mathrm{Crit}\bigl((\out_{\np_0},X_{\np_0}),X_{\np}\bigr)=1$ for any $\np_0=(\mathcal G(F_0),F_0)$ with  $F_0\in (\frac{4}{3},\frac{3}{2})$. This is a piece of the curve that joins $\left(-\frac{3}{2},\frac{3}{2}\right)$ and 
$\left(-\frac{1}{2},1\right)$, see \figc{diagrama}, and in this regard observe that the criticality is 2 for $\np_0=(\mathcal G(\frac{4}{3}),\frac{4}{3}).$ Furthermore, it is proved in \cite[Theorem B]{Rojas} that if $\np_0=(D_0,2)$ with $D_0\in (-2,0)\setminus\{-\frac{1}{2}\}$ then $\mathrm{Crit}\bigl((\out_{\np_0},X_{\np_0}),X_{\np}\bigr)=1$. The same conclusion is true for any $\np_0=(-F_0,F_0)$ with $F_0\in [\frac{3}{2},2)$ thanks to \cite[Theorem C]{MMV20}. In that paper it is also partially proved the claim about the parameter $\np_0=(-2,2)$ in assertion~$(e)$ of \teoc{Loud}. Apart from these references to previous results we also want to point out the following issues with regard to the statement and proof of \teoc{Loud}:

\begin{itemize}

\item As expected, the study of the bifurcation of critical periodic orbits, either from the inner or the outer boundary, is much more delicate when we perturb an isochronous center. By the result of W.S. Loud, see \cite{Loud}, we know that there are four nonlinear quadratic isochrones,
\begin{equation}\label{isochronous}
 \np_1=(0,1),\; \np_2=(-1/2,2),\; \np_3=(-1/2,1/2)\text{ and }\np_4=(0,1/4),
\end{equation}
which are located in $\Gamma_B.$ Chicone and Jacobs prove, see \cite[Theorem 3.1]{Chicone}, that the criticality of each isochrone $\np_i$ at the inner boundary of its period annulus (i.e., the center itself) is one. The proof of this follows by finding a finite set of generators for the ideal of \emph{all} the coefficients of the Taylor series of $P(s;\np)$ at $s=0$. In the present paper we are able to show that $\np_2$ and $\np_3$ have criticality one also at the outer boundary (i.e., the polycycle), see Propositions~\ref{I2} and~\ref{I4} respectively. A crucial point to see this is that, as we prove in \cite{MV20b}, the flatness of the remainder in the asymptotic expansion at $s=1$ is preserved after the derivation with respect to parameters. This constitutes the cornerstone to obtain \lemc{Fdiv}, which enables us to perform a convenient division in the space of coefficients and proceed then as in the proof of Bautin \cite[\S 3]{Bautin} for the analogous result about the bifurcation of limit cycles from the center. The isochrones $\np_1$ and $\np_4$ cannot be analyzed following this approach because the polycycle at the outer boundary is not hyperbolic. 

\item It is well-known, see \cite[Theorem 3.2]{Chicone}, that the criticality at the inner boundary of any quadratic center is at most two and that this maximum criticality is achieved at three parameter values, the so-called Loud points, which we give in \refc{Li}. For consistency with Chicone's conjecture, each one of these three parameters should have a ``twin'' where the maximum criticality at the outer boundary is attained. In this paper we identify two of these twin parameters, see assertion $(d)$ in \teoc{Loud}. We conjecture that each pair of twins is connected by a curve that consists of local bifurcation values at the interior, see Remarks~\ref{figdobe_com} and~\ref{figdoblebis_com}.

\item The local bifurcation values of the period function can only occur at the inner boundary (i.e., the center), at the outer boundary (i.e., the polycycle) or at the interior of the period annulus, see \lemc{diag}. (With regard to the latter, its counterpart in the context of Hilbert's 16th Problem is the bifurcation from a semi-stable limit cycle, which is characterized by the sudden emergence of a double limit cycle that gives rises to two hyperbolic limit cycles with different stability, see~\cite[\S 13.3]{Hale} for instance). As occurs with limit cycles, the identification of this third type of local bifurcation value is out of reach for the moment and only partial results have been obtained. Thus, in a joint paper with P.~Marde\v si\'c we prove (see \cite[Theorem 4.3]{MMV2}) that at each Loud point there exists a germ of analytic curve that consists of local bifurcation values at the interior. Since $P(s;\np)$ extends analytically to $s=0$, this follows readily by applying the Weierstrass Preparation Theorem. In the present paper, see assertion $(e)$ in \teoc{Loud}, we show the existence of two germs of curve which also consists of local bifurcation values at the interior and that are the mirror image at the outer boundary (i.e., at $s=1$) of those previously obtained in~\cite{MMV2}, see Figures~\ref{figdoble} and~\ref{doblebis}.

\end{itemize}

In another vein it is well-known (see \cite[\S 2.2]{Roussarie} for details) that the problem of proving the \emph{existence} of a uniform bound for the number of limit cycles in a given family, for instance Hilbert's 16th Problem, can be replaced by a local problem that consists in showing that the cyclicity of each limit periodic set within the family is finite. The proof of this is by a compactness argument and it does not provide an algorithm to compute an explicit upper bound even if we had an explicit bound for the cyclicity of every limit periodic set. In any case this gives a program for solving the existential Hilbert's 16th Problem that has been posed and implemented for the quadratic vector fields by R. Roussarie and his collaborators (see~\cite{DRR1,Roussarie1}). One can of course transfer this problem to the period function and ask for the existence of a uniform bound for the number of critical periodic orbits in the family of quadratic centers. Similarly as it occurs in the context of limit cycles, an affirmative answer would follow if one can prove that the criticality of the period function at the outer boundary of any quadratic center is finite, cf. \lemc{finitud}. On account of \teoc{Loud} we are not very far from proving the existence of this uniform bound for the family of reversible quadratic centers. It will follow in particular if one can prove the validity of the following more specific conjecture:

\begin{proclama}{Conjecture.}
Let $\{X_{\np},\np\in\R^2\}$ be the family of quadratic vector fields given in~\refc{sist_loud} and consider the period function of the center at the origin. Then the following assertions are true:
\begin{enumerate}[$(a)$]
\item $\mathrm{Crit}\bigl((\out_{\np_0},X_{\np_0}),X_{\np}\bigr)=0$ for $\np_0\in\{D=-1,F\in[0,1]\}$.
\item $\mathrm{Crit}\bigl((\out_{\np_0},X_{\np_0}),X_{\np}\bigr)=1$ for $\np_0\in\{D=0,F\in (0,2]\}$.
\item $\mathrm{Crit}\bigl((\out_{\np_0},X_{\np_0}),X_{\np}\bigr)=2$ for $\np_0=(0,0)$. 
\item There is a curve of local bifurcation values of the period function at the interior arriving to $\np=(0,0)$ 
         tangent to $D=0$.
\end{enumerate}
\end{proclama}

As a matter of fact to show the existence of a uniform bound for the number of critical periodic orbits of the reversible quadratic centers it suffices to verify that $\mathrm{Crit}\bigl((\out_{\np_0},X_{\np_0}),X_{\np}\bigr)$ is finite for all 
$\np_0=(D_0,F_0)$ inside the segments $\{-1\}\times [0,1]$ and $\{0\}\times [0,2]$.
To put this into context let us recall that the differential system~\refc{BDF} has no critical periodic orbits if $B=0$ by \cite[Theorem 1]{GGV}. On the other hand, apart from the reversible one, there are essentially three other families of quadratic centers: the Hamiltonian, the codimension four $Q_4$ and the generalized Lotka-Volterra systems $Q_3^{LV}.$ According to Chicone's conjecture the number of critical periodic orbits should be zero for the centers in these three families. This is known to be true for the Hamiltonian and $Q_4$ families thanks to the results of Coppel and Gavrilov \cite{CopGav} and Zhao \cite{Zhao}, respectively. With regard to the family $Q_3^{LV}$ it is proved in \cite{Jordi2} that, except for a subset of codimension one in the parameter plane, the criticality at the outer boundary is zero. It is clear then that any contribution to the proof of the above conjecture will constitute a very significant step forward to the existence of a uniform bound for the number of critical periodic orbits in the \emph{whole} family of quadratic centers. Let us mention in this respect that along $D=-1$ and $D=0$  the singularity at the outer boundary of the period annulus is nilpotent. In this situation the results of \cite{MV20b,MV21} do not apply and new techniques must be developed. 

The paper is organized in the following way. In \secc{crvsbif} we recall the definition of local bifurcation value at the outer boundary, that we introduce in our early paper~\cite{MMV2} to study the bifurcation diagram of the period function of the family $\{X_\np,\np\in\R^2\}$, and we prove several results that relate it with the criticality. We also show how to study the criticality by means of a suitable parametrization of the set of periodic orbits near the outer boundary.  \secc{sec2} is devoted to the asymptotic expansion of the period function near the outer boundary, which is the cornerstone in the proof of \teoc{Loud}. To this end we prove three results that are addressed to three different parameter subsets according to the contour of the period annulus. As one might expect the proofs of these results are rather long and technical. Furthermore they are based on previous tools developed in \cite{MV20a,MV20b,MV21} that need to be introduced appropriately. For these reasons, to ease the paper's readability we defer some proofs to \apc{appA}. In \secc{distinguished} we study three distinguished parameters. On one hand the two isochrones for which we succeed in proving that the criticality is one (see Propositions~\ref{I2} and~\ref{I4}) and, on the other hand, the parameter $\np=(\mathcal G(\frac{4}{3}),\frac{4}{3})$, which is also rather special because it has criticality two (see \propc{doble}). Due to the novel approach of its proof we think that each one of these results is of particular interest in the context of Hilbert's 16th Problem. \secc{provaA} is entirely devoted to the proof of \teoc{Loud}. Next, in \apc{appA} we prove the results stated in \secc{sec2} that we mentioned before and in \apc{ApBeta} we are concerned with the integral representation of the Beta and hypergeometric functions, which usually appear as coefficients in the asymptotic expansions that we obtain. Finally \apc{apB} is addressed to prove a technical result that is used to study the vanishing set of two coefficients.



\section{Criticality vs bifurcation}\label{crvsbif}

In this section we recap the notion of local bifurcation value of the period function at the outer boundary as we introduced in our early paper \cite{MMV2}. We relate it with the criticality, which is a more quantitative and geometric definition, and prove a general result connecting both notions, see \lemc{BZ}. More specifically our aim is to take advantage in the present paper of the results that we obtained in~\cite{MMV2} with regard to the period function of Loud's centers \refc{sist_loud} and that are not stated using the notion of criticality. Related with this issue, our goal in this section is also to clarify the usage of a parametrization of the period function near the outer boundary to compute its criticality, see \lemc{CZ}. Finally we give a sufficient condition in order that a parameter is a local bifurcation value of the period function at the interior, see \lemc{lema-interior}.

Several results in this section are equally valid in the finitely smooth class $\cc^k$, $k\in\N,$ the infinitely smooth class $\cc^\infty$ and the analytic class $\cc^\omega$. For simplicity in the exposition we write $\cc^\varpi$ with the wild card $\varpi\in\N\cup\{\infty,\omega\}.$ Our first result is addressed to the regularity properties of the 
map $(p,\mu)\mapsto\hat P(p;\mu)$ that assigns to each $\mu\in U$ and $p\in\PA_\mu$ the period $\hat P$ of the periodic orbit of $X_\mu$ passing through the point $p.$ The result is given under a technical assumption concerning the existence of a continuous parametrization $\sigma(s;\mu)$ of the period annulus $\PA_\mu$ near its outer boundary $\out_\mu$. We point out that from now on, in contrast with the notation used in the introduction, for the sake of convenience $s=0$ corresponds to $\out_\mu$ and $s=1$ to the center.

\begin{lem}\label{obert}
Let us fix $\varpi\in\N\cup\{\infty,\omega\}$ and consider a $\cc^\varpi$ family of planar vector fields $\{X_\mu\}_{\mu\in U}$ such that, for each $\mu\in U$, $X_\mu$ has a center $p_\mu\in\R^2$ with period annulus $\PA_\mu$. Suppose that there exists a continuous map $\sigma:(0,\delta)\times U\to\R^2$ verifying, for each fixed $\mu\in U$, that

\begin{enumerate}[$(a)$]

\item the map $\sigma(\,\cdot\,;\mu):(0,\delta)\to\R^2$ is $\mathscr C^1$,

\item the vectors $\partial_s\sigma(s;\mu)$ and $X_\mu(\sigma(s;\mu))$ are linearly independent 
         for all $s\in (0,\delta),$ and

\item for each compact set $K\subset\PA_\mu\cup\{p_\mu\}$ there exists $s_K>0$ such that 
        $\sigma(s;\mu)\in\PA_\mu\setminus K$ for all $s\in (0,s_K).$
        
\end{enumerate}
Then the following assertions hold:
\begin{enumerate}[$1.$]
\item $\mathscr U=\bigcup\limits_{\mu\in U}\PA_\mu\times\{\mu\}$ is an open subset of $\R^2\times U$, and
\item the map $(p,\mu)\mapsto \hat P(p;\mu)=\{\text{period of the periodic orbit of $X_\mu$ passing through $p$}\}$ 
         is $\cc^\varpi$ on $\mathscr U.$ 
\end{enumerate}
\end{lem}

\begin{prova}
We consider the family $\{X_\mu\}_{\mu\in U}$ as a single $\mathscr C^\varpi$ vector field $Y$ on $\R^2\times U$ whose trajectories are contained in the submanifolds $\mu=\text{constant}$. Denote the flow of $Y$ by $\phi(t;p,\mu)=(\varphi(t;p,\mu),\mu)$. In order to prove the first assertion, for a given $(p_0,\mu_0)\in\mathscr U$ we must show that there is an open subset~$V$ of $\R^2\times U$ such that $(p_0,\mu_0)\in V\subset\mathscr U.$ We claim that this is true in the particular case that there exists $s_0\in (0,\delta)$ such that $\sigma(s_0;\mu_0)=p_0.$ Indeed, due to the assumption in $(b)$, note that $Y$ is transverse to
\[
 \Sigma_\varepsilon\!:=\{(\sigma(s;\mu),\mu);\text{ $|s-s_0|<\varepsilon$ and $\|\mu-\mu_0\|<\varepsilon$}\}
\] 
for all $\varepsilon>0$ small enough and that $(p_0,\mu_0)\in\Sigma_\varepsilon.$
Then, since $\sigma:(0,\delta)\times U\to\R^2$ is continuous, by the flow box theorem (and shrinking $\varepsilon>0$ if necessary) it follows that 
\[
 V\!:=\bigcup_{t\in (-\varepsilon,\varepsilon)} \phi(t;\Sigma_\varepsilon)
\]
is an open subset of $\R^2\times U$. Furthermore, since $\mathscr U$ is invariant by $\phi$ and $\Sigma_\varepsilon\subset\mathscr U$ by construction, we have that $(p_0,\mu_0)\in V\subset\mathscr U$ and this proves the claim. Let us consider now an arbitrary $p_0\in\PA_{\mu_0}$. Denote the periodic orbits of $X_{\mu_0}$ passing through $q\!:=\sigma(\delta/2;\mu_0)$ and $p_0$ by $\gamma_q$ and $\gamma_{p_0}$, respectively. For each $\mu\in U$ we take the orthogonal vector field to $X_\mu$, say $X_\mu^\bot$, pointing inward the periodic orbits in $\PA_\mu$. We consider the family $\{X_\mu^\bot\}_{\mu\in U}$ as a single $\mathscr C^\varpi$ vector field $\hat Y$ on $\R^2\times U$ and denote its flow by $\hat\phi(t;p,\mu)=(\hat\varphi(t;p,\mu),\mu)$. 
Note that $p_\mu$ is also a singular point for $X_{\mu}^\bot$ that, by applying the Poincaré-Bendixson Theorem (see for instance \cite{ADL}), it is easy to show to be asymptotically stable. Observe moreover that $\hat\phi(t;\mathscr U)\subset\mathscr U$ for all $t\geqslant 0.$ We define $\Gamma\!:=\{\hat\varphi(t;q,\mu_0);\,t\geqslant 0\}\subset\PA_{\mu_0},$ which is clearly a transverse section for $X_{\mu_0}$, and distinguish two cases:
\begin{itemize}
\item Case 1: $\Gamma\cap\gamma_{p_0}\neq\emptyset$. In this case there exist $t_1,t_2\geqslant0$ such that 
         $\hat\varphi(t_2;q,\mu_0)=\varphi(t_1;p_0,\mu_0)\in\Gamma.$ Since $q=\sigma(\delta/2;\mu_0),$ 
         on account of the claim we can take 
         an open neighbourhood $V_1$ of 
         $(q,\mu_0)$ in $\R^2\times U$ with $V_1\subset\mathscr U.$ 
         Then, by the continuity of solutions 
         with respect to initial conditions, there exists an open neighbourhood  $V_2$ of $(p_0,\mu_0)$ such that 
         $\hat\phi\big(-t_2;\phi(t_1;V_2)\big)\subset V_1.$ Thus 
         $V_2\subset \phi\big(-t_1;\hat\phi(t_2;V_1)\big)\subset\mathscr U,$ where the second inclusion follows due to 
         the $\hat\phi(t;\mathscr U)\subset\mathscr U$ for all $t\geqslant 0$ and 
         $\phi(t;\mathscr U)=\mathscr U$ for all $t\in\R,$ together with the fact that $V_1\subset\mathscr U.$  

\item Case 2: $\Gamma\cap\gamma_{p_0}=\emptyset$. Note that in this case
        $\text{Int}(\gamma_q)\subset \text{Int}(\gamma_{p_0})$. (Here, given a Jordan curve $\gamma\subset\R^2$,  $\text{Int}(\gamma)$ denotes the bounded connected component of $\R^2\setminus\{\gamma\}$.) Moreover, by the assumption in~$(c)$ and taking 
        $K=\overline{\text{Int}(\gamma_{p_0})}$, 
        there exists $s_1\in (0,\delta/2)$ satisfying that $\sigma(s_1;\mu_0)\in\PA_{\mu_0}\setminus K.$ Therefore, since 
        $q=\sigma(\delta/2;\mu_0)\in\text{Int}(\gamma_{p_0})$, by continuity there exists $s_2\in (s_1,\delta/2)$ such that 
        $\sigma(s_2;\mu_0)\in\gamma_{p_0}.$ Consequently 
        $\sigma(s_2;\mu_0)=\varphi(t_3;p_0,\mu_0)$ for some $t_3\in\R$ and on the other hand, 
        again on account of the claim,
        there exists an open neighbourhood $V_3$ of 
        $(\sigma(s_2;\mu_0),\mu_0)$ in $\R^2\times U$ with $V_3\subset\mathscr U$. Thus, exactly as before, 
        by continuity 
        of solutions with respect to initial conditions, there is an open neighbourhood~$V_4$ of $(p_0,\mu_0)$ 
        such that $V_4\subset \phi(-t_3;V_3)\subset\mathscr U.$

\end{itemize}
This proves the validity of the first assertion. 

Let us prove now that the function \map{\hat P}{\mathscr U}{(0,+\infty)} defined by 
\[
(p,\mu)\mapsto \hat P(p;\mu)=\{\text{period of the periodic orbit of $X_\mu$ passing through $p$}\} 
\]
is $\mathscr C^\varpi$. In what follows we shall use the notation $p=(x,y)$ for the components of a point of $\R^2.$
We fix $(\hat p,\hat\mu)\in\mathscr U$ and suppose that the period of the periodic orbit of $X_{\hat\mu}$ passing through $\hat p=(\hat x,\hat y)\in\PA_{\hat\mu}$ is $\hat\tau>0.$ Then, due to $X(\hat p;\hat\mu)\!:=X_{\hat\mu}(\hat p)\neq(0,0)$,
 there is $i\in\{1,2\}$ such that 
\[
\partial_t\varphi_i(\hat\tau;\hat p,\hat\mu)=\partial_t\varphi_i(0;\hat p,\hat\mu)=X_i(\hat p;\hat\mu)\neq 0.
\]
For simplicity in the exposition let us suppose that $X_1(\hat p;\hat\mu)>0.$ In this case we can apply the Implicit Function Theorem to the equation $\varphi_1(t;p,\mu)=x$ at $(t,p,\mu)=(\hat\tau,\hat p,\hat\mu)$ in order to obtain a $\cc^\varpi$ positive function ${S}(p;\mu)$ in a open neighbourhood $W\subset\mathscr U$ of $(\hat p,\hat\mu)$ verifying ${S}(\hat p;\hat\mu)=\hat\tau$ and 
\begin{equation}\label{A4eq1}
 \left.\varphi_1(t;p,\mu)\right|_{t={S}(p;\mu)}=x\text{ for all $(p,\mu)\in W.$}
\end{equation}
Clearly we can assume that $W$ is a cube $Q({\varepsilon_1})$ with center $(\hat p,\hat\mu)$ and edge length $\varepsilon_1>0$. We diminish~$\varepsilon_1$ if necessary so that $X_1(p;\mu)>0$ for all $(p,\mu)\in Q({\varepsilon_1})$. Furthermore, thanks to ${S}(\hat p;\hat\mu)=\hat\tau$ together with the continuity of ${S}$ and $\phi$, we can take $\varepsilon_2\in (0,\varepsilon_1)$ such that 
\[
 \left.\phi(t;p,\mu)\right|_{t={S}(p;\mu)}\in Q({\varepsilon_1})\text{ for all $(p,\mu)\in Q({\varepsilon_2}).$}
\]
We claim that $\hat P=S$ on $Q({\varepsilon_2})$. Clearly the claim will follow once we show that 
\[
 \left.\varphi_2(t;p,\mu)\right|_{t={S}(p;\mu)}=y\text{ for all $(p,\mu)\in Q(\varepsilon_2)$.}
\]
By contradiction, suppose that there exists $(\bar p,\bar\mu)\in Q({\varepsilon_2})$ such that $\varphi_2\big(t;\bar p,\bar\mu\big)|_{t={S}(\bar p;\bar\mu)}\neq \bar y$. Due to $Q(\varepsilon_2)\subset\mathscr U,$ the trajectory of $X_{\bar\mu}$ passing through $\bar p$ is a periodic orbit which, for simplicity in the exposition, we assume to travel clockwise around the center~$p_{\bar\mu}$ (the other case follows verbatim). That being said we consider the piece of trajectory 
\begin{align*}
&\ell\!:=\left\{\varphi(t;\bar p,\bar\mu);\,t\in [0,{S}(\bar p,\bar\mu)]\right\}
\intertext{and the vertical segment, recall \refc{A4eq1},}
&\Gamma\!:=\left\{(1-s)\bar p+s\varphi\big({S}(\bar p,\bar\mu);\bar p,\bar\mu\big);\,s\in (0,1)\right\}\subset\{x=\bar x\}.
\end{align*}
Arguing on the phase portrait of $X_{\bar\mu}$, due to $X_1(p;\bar\mu)>0$ for all $p\in\Gamma$, if $\varphi_2\big({S}(\bar p,\bar\mu);\bar p,\bar\mu\big)< \bar y$ then interior of the Jordan curve $\ell\cup\Gamma$ is a positively but not negatively invariant subset of $\PA_{\bar\mu}.$ Similarly, if $\varphi_2\big({S}(\bar p,\bar\mu);\bar p,\bar\mu\big)>\bar y$ then we obtain a negatively invariant subset of $\PA_{\bar\mu}$ which is not positively invariant. In both cases we get a contradiction with the fact that $\PA_{\bar\mu}$ is foliated by periodic orbits of~$X_{\bar\mu}$ and $Q(\varepsilon_2)\subset\mathscr U.$ Consequently $\varphi_2\big({S}( p, \mu); p, \mu\big)=  y$ for all $(p,\mu)\in Q(\varepsilon_2)$ and so the validity of the claim follows. Since $Q(\varepsilon_2)$ is an open neighbourhood of an arbitrary point of $\mathscr U$ and $S$ is $\cc^\varpi$ in $Q(\varepsilon_2)$, the claim implies the second assertion in the statement.
\end{prova}

The previous result is addressed to a family $\{X_{\mu}\}_{\mu\in U}$ of vector fields and this is the reason why we require the existence of a local transverse section near the outer boundary of the period annulus~$\out_\mu$ that behaves well with respect to parameters. That being said, \lemc{obert} can be applied to a single vector field $X$ without this requirement because a trajectory of the orthogonal vector field $X^\bot$ already provides a transverse section in the whole period annulus. Thus in order to assert that $p\mapsto \hat P(p;\mu)$ is $\cc^\varpi$ on $\PA_\mu$ for each fixed $\mu\in U$, it is not necessary to verify the existence of a continuous map $\sigma:(0,\delta)\times U\to\R^2$ satisfying $(a)$, $(b)$ and $(c).$ 

\begin{obs}\label{local_ind}
If $X$ is a $\cc^\varpi$ vector field, $\varpi\in\N\cup\{\infty,\omega\},$ with a center then the period function $\hat P$ is a first integral for the flow of $X$ on the period annulus $\mathscr P$ that, by \lemc{obert}, is $\cc^\varpi$. Consequently the scalar product $\nabla\hat P(p)\cdot X(p)$ is zero for all $p\in\mathscr P$.
This implies that if $\gamma$ is a critical periodic orbit of $X$ then the gradient $\nabla\hat P$ vanishes on $\gamma.$ Indeed, if \map{\sigma}{(0,1)}{\mathscr P} is a $\cc^\varpi$ transverse section to $X$ on $\mathscr P$ and 
$P(s)\!:=\hat P(\sigma(s))$ then $P'(s)=\nabla\hat P(\sigma(s))\cdot\sigma'(s).$ 
Thus, since $\nabla\hat P(\sigma(s))\cdot X(\sigma(s))=0$, the transversality of $\sigma$ implies that $P'(s)=0$ if, and only if, $\nabla\hat P(\sigma(s))=(0,0).$ This shows in particular that the condition for $\gamma$ to be a critical periodic orbit is local and independent of the particular transverse section used to parametrize the set of critical periodic orbits near $\gamma.$ 
\end{obs}

We define next the notion that enable us to effectively study the criticality at the outer boundary.

\begin{defi}\label{Z}
Let $U$ be an open set of $\R^N$ and consider a family of functions $\{h(\,\cdot\,;\mu)\}_{\mu\in U}$ on $(0,\varepsilon).$ Given any $\mu_\star\in U$ we define $\mathcal Z_0(h(\,\cdot\,;\mu),\mu_\star)$ to be the smallest integer $n$ having the property that there exist $\delta>0$ and a neighbourhood~$V$ of $\mu_\star$ such that for every $\mu\in V$ the function $h(s;\mu)$ has no more than $n$ isolated zeros on $(0,\delta)$ counted with multiplicities. 
\end{defi}

The hypothesis with regard to the local transverse section in our next result are slightly stronger than in the previous one because we require the continuity at $s=0$ and that $\sigma(0;\mu)$ belongs to the outer boundary $\out_\mu$ for all $\mu\in U,$ cf. assumption $(c)$ in \lemc{obert}. We also remark that in the statement $\hat P(p;\mu)$ stands for the period of the periodic orbit of $X_\mu$ passing through $p\in\PA_\mu.$

\begin{lem}\label{CZ}
Let us consider a $\cc^\omega$ family $\{X_{\mu}\}_{\mu\in U}$ of planar polynomial vector fields such that, for each $\mu\in U$, $X_\mu$ has a center $p_\mu\in\R^2$ with period annulus $\PA_\mu$. 
Let $\out_\mu\subset\RP^2$ be the outer boundary of~$\PA_\mu.$ Suppose there exists a continuous map $\sigma:[0,\delta)\times U\to\RP^2$ verifying that, for each $\mu\in U,$
\begin{enumerate}[$(a)$]

\item the map $\sigma(\,\cdot\,;\mu):(0,\delta)\to\R^2$ 
         is $\mathscr C^1$, 
\item the vectors $\partial_s\sigma(s;\mu)$ and $X_\mu(\sigma(s;\mu))$ are linearly independent 
         for all $s\in (0,\delta)$,
         
\item $\sigma(s;\mu)\in\PA_\mu$ for all 
         $s\in (0,\delta)$ and $\sigma(0;\mu)\in\out_\mu$.
                  
\end{enumerate}
Then, for each fixed $\mu_\star\in U$, the following assertions hold:
\begin{enumerate}[1.]
\item The Hausdorff distance between the outer boundaries $\out_\mu$ and $\out_{\mu_\star}$ tends 
         to zero as $\mu\to\mu_\star.$ 
\item If $P(s;\mu)\!:=\hat P(\sigma(s;\mu);\mu)$ for all $(s,\mu)\in (0,\delta)\times U$, then         
\begin{enumerate}[$(2a)$]         
         \item 
         $\mathrm{Crit}\bigl((\out_{\mu_\star},X_{\mu_\star}),X_{\mu}\bigr)\leqslant\mathcal Z_0(P'(\,\cdot\,;\mu),\mu_\star).$
         \item $\mathrm{Crit}\bigl((\out_{\mu_\star},X_{\mu_\star}),X_{\mu}\bigr)\geqslant n$ if for each 
                  open neighbourhood $V$ of $\mu_\star$ and $\delta>0$ there exist $n$ different numbers 
                 $s_1,s_2\ldots,s_n\in (0,\delta)$ and $\hat\mu\in V$ such that $P'(s_i;\hat\mu)=0$ for $i=1,2,\ldots,n.$
         \item $\mathrm{Crit}\bigl((\out_{\mu_\star},X_{\mu_\star}),X_{\mu}\bigr)=0$ if, and only if, 
                   $\mathcal Z_0(P'(\,\cdot\,;\mu),\mu_\star)=0.$                 
\end{enumerate}
\end{enumerate}
\end{lem}

\begin{prova}
To show the first assertion note that, since $X_{\mu}$ is polynomial, we can consider its Poincaré compactification $p(X_\mu)$, see \cite[\S 5]{ADL} for details, which is an analytic vector field on the sphere~$\Sc^2$ topologically equivalent to $X_\mu.$ The outer boundary~$\out_\mu$ becomes then a polycycle of $p(X_\mu)$ that can be studied using local charts of $\Sc^2.$ On account of this, the fact that $d_H(\out_\mu,\out_{\mu_{\star}})\to 0$ as $\mu\to\mu_{\star}$ follows by the continuity of $\mu\mapsto\sigma(0;\mu)\in\out_\mu$ together with the continuity with respect to initial conditions and parameters of the trajectories of $p(X_\mu).$ The interested reader is referred to \cite[Lemma 22, p. 110]{Roussarie} for a related result for limit periodic sets.

With regard to the upper bound in $(2a)$ it is clear that if $\mathcal Z_0(P'(\,\cdot\,;\mu),\mu_\star)=+\infty$ then there is nothing to be proved. So let us assume that $\mathcal Z_0(P'(\,\cdot\,;\mu),\mu_\star)=\ell\in\Z_{\geq 0}$ and argue by contradiction. If $\mathrm{Crit}\bigl((\out_{\mu_\star},X_{\mu_\star}),X_{\mu}\bigr)\geqslant\ell+1$ then there exist $\ell+1$ sequences $\{\gamma^k_{\mu_{i}}\}_{i\in\N}$, $k=1,2,\ldots,\ell+1$, where
$\gamma^1_{\mu_{i}},\gamma^2_{\mu_{i}},\ldots,\gamma^{\ell+1}_{\mu_{i}}$ are different critical periodic orbits of~$X_{\mu_{i}}$ for each $i\in\N$, such that ${\mu_{i}}\to\mu_\star$ and $d_H(\gamma^k_{\mu_{i}},\Pi_{\mu_\star})\to 0$ as $i\to +\infty$. Then, due to $d_H(\out_\mu,\out_{\mu_{\star}})\to 0$ as $\mu\to\mu_{\star}$ and
 \[
 d_H(\gamma^k_{\mu_{i}},\Pi_{\mu_{i}})\leqslant d_H(\gamma^k_{\mu_{i}},\Pi_{\mu_\star})+d_H(\Pi_{\mu_{i}},\Pi_{\mu_\star}),
 \]
we have $d_H(\gamma^k_{\mu_{i}},\Pi_{\mu_{i}})\to 0$ as $i\to+\infty$ for each $k=1,2,\ldots,\ell+1.$ Since $\sigma(0;{\mu_{i}})\in\out_{\mu_{i}}$ and there is a one-to-one correspondence between zeros of $P'(s;\mu_i)$
arbitrarily near $s=0$ and critical periodic orbits of~$X_{\mu_i}$ arbitrarily close to~$\out_{\mu_i}$ (cf. \cite[Lemma 22]{Roussarie}), this implies that there exist $\ell+1$ sequences of positive numbers $\{s_i^k\}_{i\in\N}$, $k=1,2,\ldots,\ell+1,$ such that $P'(s_i^k;\mu_i)=0$ and $\#\{s_i^1,s_i^2,\ldots,s_i^{\ell+1}\}=\ell+1$ for each $i\in\N$, and $\lim_{i\to+\infty}s_i^k=0$ for each $k=1,2,\ldots,\ell+1.$ This clearly contradicts that $\mathcal Z_0(P'(\,\cdot\,;\mu),\mu_\star)=\ell$, see \defic{Z}. The assertion in $(2b)$ follows similarly. Indeed, on account of the assumption and the above mentioned one-to-one correspondence between zeros of $P'(s;\mu)$ near $s=0$ and critical periodic orbits of~$X_{\mu}$ close to~$\out_{\mu_i},$ we can construct ${n}$ sequences $\{\gamma^k_{\mu_{i}}\}_{i\in\N}$, $k=1,2,\ldots,{n}$, where
$\gamma^1_{\mu_{i}},\gamma^2_{\mu_{i}},\ldots,\gamma^{{n}}_{\mu_{i}}$ are different critical periodic orbits of~$X_{\mu_{i}}$ for each $i\in\N$, such that ${\mu_{i}}\to\mu_\star$ and $d_H(\gamma^k_{\mu_{i}},\Pi_{\mu_i})\to 0$ as $i\to +\infty$. Then, using that $d_H(\out_\mu,\out_{\mu_{\star}})\to 0$ as $\mu\to\mu_{\star}$, we can assert that $\lim_{i\to+\infty}d_H(\gamma^k_{\mu_{i}},\Pi_{\mu_{0}})=0$ for each $k=1,2,\ldots,n$, which implies $\mathrm{Crit}\bigl((\out_{\mu_\star},X_{\mu_\star}),X_{\mu}\bigr)\geqslant n$, as desired. Finally the assertion in~$(2c)$ follows easily from the ones in $(2a)$ and $(2b).$ This completes the proof of the result.
\end{prova}

Next we introduce the notion of global transverse section for a family of period annuli. Roughly speaking it is a transverse section, joining the center with some point at the outer boundary of the period annulus, that behaves well with the parameters. 

\begin{defi}\label{gts}
Let us fix $\varpi\in\N\cup\{\infty,\omega\}$ and consider a $\mathscr C^\varpi$ 
family $\{X_\mu\}_{\mu\in U}$ of planar vector fields such that, for each $\mu\in U,$ $X_\mu$ has a center $p_\mu\in\R^2$ with period annulus $\PA_\mu$. Let $\out_\mu\subset\RP^2$ be the outer boundary of $\PA_\mu.$ A \emph{global transverse section} for the family of period annuli $\{\PA_\mu\}_{\mu\in U}$ is a continuous map $\sigma:[0,1]\times U\to\RP^2$ verifying that
\begin{enumerate}[$(a)$]

\item the map $\sigma(\,\cdot\,;\mu):[0,1]\to\mathbb{RP}^2$ is $\mathscr C^\varpi$ for each $\mu\in U$,

\item the vectors $\partial_s\sigma(s;\mu)$ and $X_\mu(\sigma(s;\mu))$ are linearly independent for all $(s,\mu)\in (0,1)\times U$ and the map $\partial_s\sigma:(0,1)\times U\to\R^2$ is continuous, 

\item $\sigma(s;\mu)\in\PA_\mu$ for all $s\in (0,1)$, $\sigma(0;\mu)\in\out_\mu$ and $\sigma(1;\mu)= p_\mu$.

\end{enumerate}
When such a global transverse section exists we say that the family of period annuli $\{\PA_\mu\}_{\mu\in U}$ varies continuously.
\end{defi}

\begin{obs}\label{+CZ}
The period annulus of the family of Loud's quadratic centers given in \refc{sist_loud} varies continuously in the sense of \defic{gts}. Indeed, it follows from the proof of \cite[Lemma 3.2]{MMV2} that 
\[
 \mu=(D,F)\mapsto \xi_\mu\!:=\sup\{t>0;\,(s,0)\in\PA_\mu\text{ for all $s\in (0,t)$}\} 
\] 
is a well-defined continuous function on $\R^2.$ Moreover the point $(\xi_\mu,0)$ belongs to $\out_\mu$ and 
$0<\xi_\mu\leqslant 1$ for all $\mu\in\R^2.$ Then $\sigma(s;\mu)=((1-s)\xi_\mu,0)$ for $(s,\mu)\in [0,1]\times\R^2$ is clearly a global transverse section. In particular, since the Loud's system is polynomial, the outer boundary of the period annulus varies continuously in the Hausdorff sense by the first assertion in \lemc{CZ}. 
\end{obs}

Note, see $(b)$ in \defic{gts}, that we also require $(s,\mu)\mapsto\partial_s\sigma(s;\mu)$ to be continuous. The reason for this is because if we define $P(s;\mu)=\hat P(\sigma(s;\mu);\mu)$ then $(s,\mu)\mapsto\partial_s P(s;\mu)$ is a continuous function by \lemc{obert}. This continuity is a key point in the forthcoming results.  Before that we summarize in the next statement the properties that we get for $P(s;\mu)$ as a consequence of \lemc{obert} and \defic{gts}.

\begin{cory}\label{pc}
Let us fix $\varpi\in\N\cup\{\infty,\omega\}$ and consider a $\cc^\varpi$ family of planar vector fields $\{X_\mu\}_{\mu\in U}$ such that, for each $\mu\in U$, $X_\mu$ has a center $p_\mu\in\R^2$ with period annulus $\PA_\mu$. Assume that the family of period annuli varies continuously and let $\sigma:[0,1]\times U\to\RP^2$ be a global transverse section for $\{\PA_\mu\}_{\mu\in U}$. If $P(s;\mu)\!:=\hat P(\sigma(s;\mu);\mu)$ for all $(s,\mu)\in (0,1)\times U$ then the following holds:
\begin{enumerate}[$(a)$]

\item $P(\,\cdot\,;\mu)\in\cc^\varpi((0,1))$ for each $\mu\in U,$ and 

\item $P$ and $\partial_s P$ are continuous functions on $(0,1)\times U$.

\end{enumerate}
\end{cory}

\begin{defi}\label{period_param}
Under the assumptions of \coryc{pc}, we say that $P(s;\mu)=\hat P(\sigma(s;\mu);\mu)$, which is defined for $(s,\mu)\in (0,1)\times U,$ is a \emph{global parametrization of the period function}. In contrast,  
\[
(p,\mu)\mapsto \hat P(p;\mu)=\{\text{period of the periodic orbit of $X_\mu$ passing through $p$}\}
\]
is defined on $\bigcup_{\mu\in U}\PA_\mu\times\{\mu\},$ which is not so easy to handle. 
\end{defi}

One of the main goals in the present section is to relate the concept of local bifurcation value of the period function, as introduced in~\cite{MMV2}, with the notion of criticality, see \defic{criticality}. As we will see the first one concerns with the qualitative properties of the period function, whereas the second is more geometric and quantitative. In doing so we will be able to take advantage of the results about the bifurcation diagram of the period function of the Loud's centers that we obtained in~\cite{MMV2}. For reader's convenience we next recall the definition of local bifurcation value of the period function. 

\begin{defi}\label{def0}
Let $\{I_\mu\}_{\mu\in U}$ be a continuous family of intervals in $\R$, i.e., such that $I_\mu=(\ell(\mu),r(\mu))$ with $\ell,r\in\cc^0(U)$, and consider a continuous family of functions $\{\map{F_{\mu}}{I_{\mu}}{\R}\}_{\mu\in U}.$ We say that $\mu_0\in U$ is a
\emph{regular value of the family} $\{\map{F_{\mu}}{I_{\mu}}{\R}\}_{\mu\in U}$ if there exist a
neighbourhood~$V$ of~$\mu_0$ and an isotopy $\{\map{h_{\mu}}{I_{\mu}}{I_{\mu_0}}\}_{\mu\in V},$ with $h_{\mu_0}=id,$ such that
 \begin{equation}\label{sgn}
  \text{sgn}\biggl({F}_{\mu}(s)\biggr)=\text{sgn}\biggl({F}_{\mu_0}\bigl(h_{\mu}(s)\bigr)\biggr)
  \text{ for all $s\in I_{\mu}$ and $\mu\in V,$}
\end{equation}
where $\mathrm{sgn}:\R\to\{-1,0,1\}$ is the extended sign function. 
 A parameter $\mu_0$ which is not regular is called a \emph{bifurcation value.}
\end{defi}

The endpoints of $I_\mu$, the domain of definition of $F_\mu$, depend continuously on $\mu,$ so that $\cup_{\mu\in{U}}I_{\mu}\!\times\!\{\mu\}$ is an open subset of $\R\times{U}.$ Thus, by a continuous family of functions $\{\map{F_{\mu}}{I_{\mu}}{\R}\}_{\mu\in U}$, we mean that the map $(s,\mu)\mapsto F_\mu(s)$ is continuous on $\cup_{\mu\in{U}}I_{\mu}\!\times\!\{\mu\}$. Next we particularize the previous definition to study the period function. To this aim note that, by \coryc{pc}, if $\{X_\mu\}_{\mu\in U}$ is a $\cc^1$ family of vector fields with a center such that the corresponding family of period annuli varies continuously, and we set $P(s;\mu)=\hat P(\sigma(s;\mu);\mu)$, then $\{\partial_s P(\,\cdot\,;\mu)\}_{\mu\in U}$ is a continuous family of functions on $(0,1).$ 

\begin{defi}\label{def2}
Consider a $\cc^1$ family of planar vector fields $\{X_\mu\}_{\mu\in U}$ such that, for each $\mu\in U$, $X_\mu$ has a center $p_\mu\in\R^2$ with period annulus $\PA_\mu$, that we suppose to vary continuously.
  \begin{enumerate}[$(a)$]
   \item We say that $\mu_0\in{U}$ is a \emph{regular} (respectively, \emph{bifurcation})
         \emph{value of the period function} if for some global parametrization of the period
         function \map{P}{(0,1)\times U}{(0,+\infty)} we have that $\mu_0$ is a
         regular (respectively, bifurcation) value of the family $\{P'(\,\cdot\,;\mu)\!:(0,1)\to \R\}_{\mu\in U}.$
   \item We say that $\mu_0\in{U}$ is a \emph{local regular
         value of the period function at the interior} if there is some global parametrization of the
         period function \map{P}{(0,1)\times U}{(0,+\infty)} satisfying that for each $c\in (0,1)$ there 
         exists a continuously varying
         neighbourhood~$I_{\mu}(c)$ of~$c$ in $(0,1)$ such that $\mu_0$ is a regular value
         of the family $\{P'(\,\cdot\,;\mu)\!:I_\mu(c)\to \R\}_{\mu\in U}.$ A
         parameter which is not a local regular value at the interior is called a
         \emph{local bifurcation value at the interior.}
   \item We say that $\mu_0\in{U}$ is a \emph{local regular
         value of the period function at the outer} (respectively, \emph{inner}) \emph{boundary}
         if for some global parametrization of the period function \map{P}{(0,1)\times U}{(0,+\infty)}
         there exists a continuously varying
         neighbourhood~$I_{\mu}(c)$ of~$c=0$ (respectively, $c=1$) such that $\mu_0$ is a regular value
         of the family $\{\map{P'(\,\cdot\,;\mu)}{I_{\mu}(c)\cap(0,1)}{\R}\}_{\mu\in{U}}.$ A
         parameter which is not a local regular value at the outer (respectively,
         inner) boundary is called a \emph{local bifurcation value at the outer}
         (respectively, \emph{inner}) \emph{boundary.}
 \end{enumerate}
\end{defi}

\begin{obs}
Let us make the following easy observations with regard to the previous definitions:
\begin{enumerate}[$(a)$]

\item One can replace ``some global parametrization'' by ``any global
parametrization''. Indeed, suppose that $\mu_0\in U$ is a regular value for $\{P'(\,\cdot\,;\mu)\}_{\mu\in U}$ where $P(s;\mu)=\hat P(\sigma(s;\mu);\mu)$ and consider another global parametrization $\bar P(s;\mu)=\hat P(\bar\sigma(s;\mu);\mu)$ of the period function, see \defic{period_param}. If we denote by $\tau_\mu(s)$ the Poincaré map from the transverse section $\Sigma$ given by $s\mapsto\sigma(s;\mu)$ to the transverse section~$\bar\Sigma$ given by $s\mapsto\bar\sigma(s;\mu)$ then~$\tau_\mu$ is an increasing diffeomorphism and $P(s;\mu)=\bar P(\tau_\mu(s);\mu),$ so that  $P'(s;\mu)=\bar P'(\tau_\mu(s);\mu)\tau'_\mu(s).$ On account of this and following the notation in \defic{def0}, $\bar h_\mu\!:=\tau_{\mu_0}\circ h_\mu\circ\tau^{-1}_\mu$ is a suitable isotopy in order to show that $\mu_0$ is a regular value for the family $\{\bar P'(\,\cdot\,;\mu\}_{\mu\in U}$ because 
\[
 \text{sgn}\big(\bar P'(s;\mu)\big)
 =
 \text{sgn}\big(P'(\tau_\mu^{-1}(s);\mu)\big)
 =
 \text{sgn}\big(P'\big((h_\mu\circ\tau_\mu^{-1})(s);\mu_0\big)\big)
 =
 \text{sgn}\big(\bar P'\big(\bar h_\mu(s);\mu_0\big)\big),
\]
where we use that $\tau'(s)>0.$ 

\item In order to study if a parameter is a local regular value at the outer boundary it is not necessary to consider a global transverse section \map{\sigma}{[0,1]\times U}{\RP^2} for
the family of period annuli. Indeed, see point~$(c)$ in \defic{def2}, it suffices to take a local parametrization \map{\sigma}{[0,\delta)\times U}{\RP^2}. Similarly, to study the local regular values at the inner boundary it suffices to take a local parametrization \map{\sigma}{(1-\delta,1]\times U}{\R^2}.
\end{enumerate}
\end{obs}

As expected, $\mu_0$ is a bifurcation value of the period function if, and only if, $\mu_0$ is either a local bifurcation value at the inner boundary, at the outer boundary or at the interior. This is stated in the following result 
and the interested reader is referred to \cite[Lemma 2.7]{MMV2} for the proof. 

\begin{lem}\label{diag}
Let us consider a $\cc^1$ family of analytic planar vector fields $\{X_\mu\}_{\mu\in U}$ such that, for each $\mu\in U$, $X_\mu$ has a center $p_\mu\in\R^2$ with period annulus $\PA_\mu$, that we suppose to vary continuously. Then the bifurcation diagram of the period function is the union of the local bifurcation
diagrams at the inner and outer boundary and in the interior.
\end{lem}

Under the assumptions and notation in \coryc{pc}, a sufficient condition for 
$\mu_\star\in U$ to be a local regular value of the period function at the interior is that $P'(s;\mu_\star)\neq 0$ for all $s\in (0,1).$ This follows easily by the continuity of $(s,\mu)\mapsto P'(s;\mu)$ on $(0,1)\times U$ and a compactness argument. In case that this function is $\cc^1$ then another sufficient condition is that 
$P'(\,\cdot\,;\mu_\star)$ has only simple zeros because the application of the Implicit Function Theorem provides the appropriate isotopies. Hence, in this context, the set of local bifurcation values of the period function at the interior is contained in 
\[
 \Delta\!:=\{\mu\in U;\text{ there exists $s\in(0,1)$ such that } P'(s;\mu)=P''(s;\mu)=0\}.
\] 
If $P'(s;\mu)$ was polynomial in $s$ (which is certainly not true) then $\Delta$ would consist of those parameters~$\mu_\star$ such that the  \emph{discriminant} of $P'(s;\mu_\star)$ is equal to zero. (Recall that the discriminant of $q\in\R[x]$ is the resultant between $q(x)$ and $q'(x)$, see for instance \cite{Cox}.) One may expect on the other hand that the parameters in~$\Delta$ are always local bifurcation values of the period function at the interior. However this is not always the case and the following toy models show that some additional assumptions are needed to this end.

\begin{ex} Setting $N=1,$ we take $P'$ to be $F(s;\mu)=(s-\mu)^2$ and $U=(0,1).$ Then it is clear that any $\mu\in U$ is a local regular value of $F$ at the interior (i.e., there are no local bifurcation values) but we have that $\Delta=U.$ Note that in this case the interior of $\Delta$ is non-empty. 
\end{ex}

\begin{ex} Setting $N=2,$ we take $P'$ to be $F(s;\mu)=(s-\mu_1)^3-\mu_2$ and $U=(0,1)^2.$ Then again it turns out that any $\mu=(\mu_1,\mu_2)\in U$ is a local regular value of $F$ at the interior, whereas $\Delta=\{\mu\in U:\mu_2=0\}.$ Observe that in this case the interior of $\Delta$ is empty but $F(\,\cdot\,;\mu)$ has zeroes of multiplicity 3.
\end{ex}
The following result provides us with an analytical tool to study the local bifurcation values of the period function at the interior. We emphasize that it has the natural hypothesis in view of the previous discussion. 

\begin{lem}
\label{lema-interior}
Let $\{X_\mu\}_{\mu\in U}$ be an analytic family of planar vector fields such that, for each $\mu\in U$, $X_\mu$ has a center $p_\mu\in\R^2$ with period annulus $\PA_\mu$. Assume that the family of period annuli varies continuously and let $\sigma:[0,1]\times U\to\RP^2$ be a global transverse section for $\{\PA_\mu\}_{\mu\in U}$. Setting $P(s;\mu)=\hat P(\sigma(s;\mu);\mu)$ for all $(s,\mu)\in (0,1)\times U$, suppose additionally that 
\begin{enumerate}[$(a)$]
\item the interior of $\Delta$ $($as a subset of $U\subset \R^N)$ is empty, and
\item for each $\mu\in U$, the zeros of $P'(\,\cdot\,;\mu)$ have at most multiplicity 2.
\end{enumerate}
Then each $\mu\in\Delta$ is a local bifurcation value of the period function at the interior.
\end{lem}

\begin{prova}
Note first that, by \coryc{pc}, the function $P(\,\cdot\,;\mu)$ is analytic on $(0,1)$ for each $\mu\in U.$ 
Let us take any $\mu_0\in\Delta.$ Then there exists $s_0\in (0,1)$ such that $P'(s_0;\mu_0)=P''(s_0;\mu_0)=0$ and, by the hypothesis in~$(b)$, $P'''(s_0;\mu_0)\neq 0$. Consequently $P'(\,\cdot\,;\mu_0)$ has a local extremum at $s=s_0$ and so there exists $\varepsilon>0$ small enough such that $P'(\,\cdot\,;\mu_0)$ has the same sign $+1$ or $-1$ on $(s_0-\varepsilon,s_0+\varepsilon)\setminus\{s_0\}.$ 
Assume by contradiction that $\mu_0$ is a local regular value of the period function at the interior. Then, taking $c=s_0$ in $(b)$ of \defic{def2}, we can consider a neighbourhood $V$ of $\mu_0$, a continuously varying neighbourhood $I_\mu$ of $s_0$ in $(0,1)$ and an isotopy $h_\mu:I_\mu\to I_{\mu_0}$ for $\mu\in V$, with $h_{\mu_0}=\mathrm{id}$, verifying the equality in~\refc{sgn}. Since $\Delta$ has empty interior we can take $\hat\mu\in V\setminus\Delta$ and define $\hat s\!:=h_{\hat\mu}^{-1}(s_0)\in I_{\hat\mu}$. On account of this, particularizing~\refc{sgn} with $\mu=\hat\mu$ and $s=\hat s$ we deduce that $P'(\hat s;\hat\mu)=0$. Accordingly, due to $\hat\mu\notin\Delta,$ it follows that $P''(\hat s;\hat\mu)\neq 0$. Therefore the function $s\mapsto P'(s;\hat\mu)$ changes sign at $s=\hat s.$ 
This contradicts~\refc{sgn} taking $\mu=\hat\mu$ and $s\approx\hat s$ because $P'(\,\cdot\,;\mu_0)$ has the same sign on $(s_0-\varepsilon,s_0+\varepsilon)\setminus\{s_0\}.$
\end{prova}

In the statement of our next result $p(X)$ stands for the Poincar\'e compactification in $\Sc^2$ of a planar polynomial vector field $X$, see \cite[\S 5]{ADL} for details. Recall also that any polycycle of an analytic vector field can be desingularized giving a polycycle with only hyperbolic or semi-hyperbolic vertices. By a hyperbolic polycycle we mean that its desingularization does not have semi-hyperbolic vertices (i.e., saddle-nodes).

\begin{lem}\label{BZ}
Consider a $\cc^\omega$ family of planar polynomial vector fields $\{X_\mu\}_{\mu\in U}$ such that, for each $\mu\in U$, $X_\mu$ has a center $p_\mu\in\R^2$ with period annulus $\PA_\mu$, that we suppose to vary continuously. 
Then the following assertions hold for any given $\mu_\star\in U$:
\begin{enumerate}[$(a)$]
\item If $\mathrm{Crit}\bigl((\out_{\mu_\star},X_{\mu_\star}),X_{\mu}\bigr)=0$ then $\mu_\star$ is a local regular 
         value of the period function at the outer boundary. 
\item Assuming additionally that the 
        outer boundary $\out_{\mu_\star}$ is a hyperbolic polycycle of $p(X_{\mu_\star})$, 
        if $\mu_\star$ is a local regular
        value of the period function at the outer boundary then
        $\mathrm{Crit}\bigl((\out_{\mu_\star},X_{\mu_\star}),X_{\mu}\bigr)=0$. 
\end{enumerate}
\end{lem}

\begin{prova}
Since the family of period annuli varies continuously,  see \defic{gts}, we can take a global transverse section $\sigma:[0,1]\times U\to\RP^2$ and consider the global parametrization of the period function given by $P(s;\mu)\!:=\hat P(\sigma(s;\mu);\mu)$ for $(s,\mu)\in (0,1)\times U,$ see \coryc{pc}.

In order to show $(a)$ note that if $\mathrm{Crit}\bigl((\out_{\mu_\star},X_{\mu_\star}),X_{\mu}\bigr)=0$ then $\mathcal Z_0(P'(\,\cdot\,,\mu),\mu_\star)=0$ by assertion $(2c)$ in \lemc{CZ}. This implies, see \defic{Z}, the existence of $\delta>0$ and a neighbourhood $V$ of $\mu_\star$ such that $P'(s;\mu)$ does not vanish on $(0,\delta)\times V$. Hence, since $(s,\mu)\mapsto P'(s;\mu)$ is continuous thanks to $(b)$ in \coryc{pc}, the function $P'(s;\mu)$ has constant sign on $(0,\delta)\times V$. Thus, see Definitions \ref{def0} and \ref{def2}, taking $I_\mu=(0,\delta)$ and $h_\mu=id$ it follows that $\mu_\star$ is a regular value of the family $\{P'(\,\cdot\,;\mu):I_\mu\to\R\}_{\mu\in U}$ as desired. This shows the validity of the assertion in $(a)$.

Let us turn next to the assertion in $(b)$. If $\mu_\star$ is a local regular
value of the period function at the outer boundary then there exist a neighbourhood $V$ of~$\mu_\star$, a continuous strictly positive function $\mu\mapsto\varepsilon_\mu$ on $V$ and an isotopy $\{\map{h_\mu}{(0,\varepsilon_\mu)}{(0,\varepsilon_{\mu_\star})}\}_{\mu\in V}$ such that $\mathrm{sgn}(P'(s;\mu))=\mathrm{sgn}(P'(h_\mu(s));\mu_\star)$ for all $s\in (0,\varepsilon_\mu)$ and $\mu\in V$. From this point we distinguish two cases:

\begin{enumerate}[1.]

\item If the center of $X_{\mu_\star}$ is not isochronous then, by applying \cite[Theorem 1.1]{MS07}, the zeros of $P'(s;\mu_\star)$ do not accumulate to $s=0.$ Let us remark that to apply this result we take into account that the transverse section $\sigma(\,\cdot\,;\mu_\star)$ is analytic at $s=0$, see \defic{gts}, and the hypothesis that~$\out_{\mu_\star}$ is a hyperbolic polycycle of $p(X_{\mu_\star})$. Hence there exists $\rho>0$ such that $P'(s;\mu_\star)\neq 0$ for all $s\in  (0,\rho).$ Thus, since we can suppose without loss of generality that $\varepsilon_{\mu_\star}\in (0,\rho)$ and $\delta\!:=\inf\{\varepsilon_\mu:\mu\in V\}>0$, it follows that $P'(s;\mu)\neq 0$ on $(0,\delta)\times V,$ which implies (see \defic{Z}) that $\mathcal Z_0(P'(\,\cdot\,,\mu),\mu_\star)=0$. Therefore, by assertion $(2c)$ in \lemc{CZ}, $\mathrm{Crit}\bigl((\out_{\mu_\star},X_{\mu_\star}),X_{\mu}\bigr)=0$.

\item If the center of $X_{\mu_\star}$ is isochronous then $P'(\,\cdot\,;\mu_\star)\equiv 0.$ Hence $\mathrm{sgn}(P'(s;\mu))=\mathrm{sgn}(P'(h_\mu(s));\mu_\star)=0$ for all $s\in (0,\varepsilon_\mu)$ and $\mu\in V$. Thus $P'(\,\cdot\,;\mu)$ has not isolated zeros for all $\mu\in V$ and consequently, see \defic{Z}, $\mathcal Z_0(P'(\,\cdot\,,\mu),\mu_\star)=0$. Then $\mathrm{Crit}\bigl((\out_{\mu_\star},X_{\mu_\star}),X_{\mu}\bigr)=0$ by $(2c)$ in \lemc{CZ}.

\end{enumerate}
This shows $(b)$ and completes the proof of the result.
\end{prova}

We conclude this section by showing that, as we explain in the introduction, \teoc{Loud} leaves us very close to the proof of the existence of an upper bound for the number of critical periodic orbits in the family $\{X_\np,\np\in\R^2\}$. In this respect we note that there are parameter values $\np\in\R^2$ for which $X_\np$ has another center $p_\nu$ apart from the one at the origin (see for instance \cite[Figure 4]{MMV2}). The bound also holds for the critical periodic orbits of this second center because one can always find an invertible affine transformation \map{g}{\R^2}{\R^2}
with $g(p_\np)=(0,0)$ such that 
the push-forward of $X_\np$ by $g$ verifies $g_*(X_\np)=\beta X_{\hat\np}$ for some $\hat\np\in\R^2$ and $\beta\neq 0.$

\begin{lem}\label{finitud}
Consider the family of vector fields $\{X_\np,\np\in\R^2\}$ given in \refc{sist_loud}. If $\mathrm{Crit}\bigl((\out_{\np_0},X_{\np_0}),X_{\np}\bigr)$ is finite for every $\np_0\in\R^2$ then there exists $N\in\N$ such that the center at the origin of $X_\np$ has at most $N$ critical periodic orbits for all $\np\in\R^2.$
\end{lem}

\begin{prova}
By \lemc{obert}, $\mathscr U\!:=\bigcup_{\np\in\R^2}\PA_\np\times\{\np\}$ is an open subset of $\R^2\times\R^2$ and the map $$(p,\np)\mapsto \hat P(p;\np)=\{\text{period of the periodic orbit of $X_\np$ passing through $p$}\}$$ is analytic on $\mathscr U.$ We define $P(s;\np)\!:=\hat P\big(\big((1-s)\xi_\np,0\big);\np\big)$ for each $(s,\np)\in (0,1)\times\R^2$, see \obsc{+CZ}, which provides us with a suitable global parametrization of the period function. Let us note in particular that $\partial_s^kP(s;\np)$ is a continuous function on $(0,1)\times\R^2$ for each $k\in\N$.
Moreover, by \cite[Theorem A]{TopaV}, we know that if $\np=(D,F)\notin K\!:=[-7,2]\times [0,4]$ then then the center at the origin of $X_\np$ has no critical periodic orbits. Consequently, if for each fixed $\np\in\R^2$ we define $N_\np$ to be the number of isolated zeros of $s\mapsto P'(s;\np)$ on the interval $(0,1)$ counted without multiplicities,  the result will follow once we prove that
\[
 \sup_{\np\in K}(N_\np)<+\infty.
\]
Let us advance that this will be a consequence of the compactness of $[0,1]\times K.$ With this end in view we fix any $(s_\star,\np_\star)\in [0,1]\times K$ and observe that three different situations may occur:
\begin{enumerate}[$(a)$]

\item {\bf Case $\mathbf{s_\star=1}$.} As a consequence of the result of Chicone and Jacobs, see 
         \cite[Theorem 3.1]{Chicone}, there exist $\varepsilon,\delta>0,$ depending on $\np_\star,$ such 
         if $\np\in B_{\varepsilon}(\np_\star)\!:=\{\np\in\R^2:\|\np-\np_\star\|<\varepsilon\}$ then the number of 
         isolated roots of $P'(s;\np)=0$ with $s\in (1-\delta,1)$ is at most 2 (counted with multiplicities).
         
\item {\bf Case $\mathbf{s_\star=0}$.} Since $\ell\!:=\mathrm{Crit}\bigl((\out_{\np_\star},X_{\np_\star}),X_{\np}\bigr)<+\infty$ 
         by assumption, $(2b)$ in \lemc{CZ} implies that there exist
         $\varepsilon,\delta>0$ (depending on $\np_\star$ again) such 
         if $\np\in B_{\varepsilon}(\np_\star)$ then the number of 
         isolated roots of $P'(s;\np)=0$ with $s\in (0,\delta)$ is at most $\ell$ (counted without multiplicities). 
         
\item {\bf Case $\mathbf{s_\star\in (0,1)}$.} 

         \begin{enumerate}[$(c1)$]
         
         \item If the center of $X_{\np_\star}$ is not isochronous then there exists 
                 $k\in\N$, depending on $(s_\star,\np_\star)$, such 
                 that $\partial_s^{k}P(s_\star;\np_\star)\neq 0.$ By continuity there is 
                 a neighbourhood $V$ of $(s_\star,\np_\star)$ such that  
                 $\partial_s^{k}P(s;\np)\neq 0$ for all $(s,\np)\in V.$ Hence the application of Rolle's 
                 Theorem shows that there exist $\varepsilon,\delta>0$ such 
                 if $\np\in B_{\varepsilon}(\np_\star)$ then the number of 
                 roots of $P'(s;\np)=0$ with $s\in (s_\star-\delta,s_\star+\delta)$ is at most $k$ (counted with multiplicities).
               
        \item Let us suppose finally that the center of $X_{\np_\star}$ is isochronous. Since 
                $\big(((1-s)\xi_{\np_\star},0),{\np_\star}\big)\in\mathscr U$, and by taking for instance the
                flow of the orthogonal vector field $X_\np^\bot$, there exists a transverse section 
                $\bar s\mapsto \sigma(\bar s;\np)$ given by 
                an analytic map
                \[
                \map{\sigma}{(-\delta_1,\delta_1)\times B_{\varepsilon_1}(\np_\star)}{\mathscr U}
                \] 
                and such that $\sigma(0;\np_\star)=\big(((1-s_\star)\xi_{\np_\star},0),{\np_\star}\big).$ We then define 
                $\bar P(\bar s;\np)\!:=\hat P(\sigma(\bar s;\np))$, which is clearly analytic on 
                $(-\delta_1,\delta_1)\times B_{\varepsilon_1}(\np_\star)$. We can thus compute its Taylor's series 
                at $\bar s=0$,
                \[
                 \bar P(\bar s;\np)=\sum_{i=0}^\infty a_i(\np)\bar s^i,
                \]
                where each $a_i$ is an analytic function on $B_{\varepsilon_1}(\np_\star)$ with $a_i(\np_\star)=0$.
                Working in the local ring $\R\{\np\}_{\np_\star}$ of convergent power series at $\np_\star$, we consider the
                ideal $\mf B\!:=\big(a_i,i\in\N\big)$. The ring is Noetherian and so there exists $\ell\in\N$ 
                such that $\mf B=(a_1,a_2,\ldots,a_\ell).$ Verbatim the proof of Chicone and Jacobs for 
                \cite[Theorem 2.2]{Chicone} (see also the result of Roussarie in \cite[\S 4.3.1]{Roussarie} for a similar
                result for the displacement map), there exist analytic functions $h_i(\bar s;\np)$ in a neighbourhood of 
                $(0,\np_\star)$ with $h_i(0;\np)\equiv 1$ for $i=1,2,\ldots,\ell$ such that we can write
                \[
                 \bar P'(\bar s;\np)=\sum_{i=1}^\ell a_i(\np)\bar s^{i-1}h_i(\bar s;\np).
                \]
                Now, setting $\psi_i(\bar s;\np)\!:=\bar s^{i-1}h_i(\bar s;\np)$ and proceeding just like the 
                proof of  \cite[Theorem 2.2]{Chicone}, one can apply the well-known 
                derivation-division algorithm and use recursively Rolle's Theorem to show that 
                there exist $\delta_2,\varepsilon_2>0$ small enough such that
                if $\np\in B_{\varepsilon_2}(\np_\star)$ then the ordered set $(\psi_1,\psi_2,\ldots,\psi_\ell)$ is an extended
                complete Cheybshev system for $\bar s\in (-\delta_2,\delta_2)$, see~\cite{KS} for a definition. 
                Accordingly
                if $\np\in B_{\varepsilon_2}(\np_\star)$ then either $\bar P'(\,\cdot\,;\np)\equiv 0$ or $P'(\bar s;\np)=0$ has
                at most $\ell-1$ roots with $\bar s\in (-\delta_2,\delta_2)$ counted with multiplicities.                               
                Using the original
                parametrization of the period function, this shows the existence of $\delta_3,\varepsilon_3>0$ small
                enough such that 
                if $\np\in B_{\varepsilon_3}(\np_\star)$ then
                the number of isolated roots of 
                $P'(s;\np)=0$ with $s\in (s_\star-\delta_3,s_\star+\delta_3)$ is at 
                most $\ell-1$ taking multiplicities into account. 
         \end{enumerate}
\end{enumerate}
Since in each one of the possible cases there is a neighbourhood of $(s_\star,\np_\star)$ where the number of critical periods is finite, the result follows by taking a finite subcover of $[0,1]\times K.$
\end{prova}


\section{Asymptotic expansion of the period function}\label{sec2}

From now on we focus on the quadratic family $\{X_\np,\np\in\R^2\}$ given in \refc{sist_loud} and study the period function of the center at the origin. In this section we give its asymptotic expansion near the outer boundary $\out_\np$ for parameters $\nu$ inside three specific sets (see \figc{diagrama}):
\begin{align*}
 \Gamma_1=&\textstyle\big\{D=-\frac{1}{2}, F\in(\frac{1}{2},1)\big\}\cup\big\{F=\frac{1}{2}, D\in(-1,0)\big\},
 \\[3pt]
 \Gamma_2=&\textstyle\big\{F=2,D\in (-2,0)\big\}\cup\big\{D=\mathcal G(F):F\in(1,\frac{3}{2})\big\}
 \intertext{and}
 \Gamma_3=&\big\{F=1, D\in(-1,0)\big\}.
\end{align*}
In all the cases the period annulus $\PA_\np$ is unbounded. Since the vector field $X_\nu$ is polynomial, in order to study the behaviour of the trajectories near infinity one can use its Poincaré compactification $p(X_\nu),$ which is an analytic vector field on the sphere~$\Sc^2$ topologically equivalent to $X_\nu$, see \cite[\S 5]{ADL} for details. The outer boundary~$\out_\nu$ is a polycycle of $p(X_\nu)$ that can be studied using local charts of $\Sc^2.$ In doing so one obtains (see~\cite[Figure 4]{MMV2}) the different phase portraits in the dehomogeneized Loud's family $\{X_\nu,\nu\in\R^2\}$. For the parameter values studied in this section it occurs that the polycycle~$\out_{\nu}$ of~$p(X_{\nu})$ is hyperbolic if $\nu\in \Gamma_1\cup\Gamma_2$ and has a saddle-node singularity if $\nu\in\Gamma_3.$ With regard to the phase portrait, it happens that the affine part of~$\out_\nu$ is a straight line for $\nu\in\Gamma_1,$ whereas it is a branch of a hyperbola for $\nu\in\Gamma_2.$ These are the reasons why we split the parameters under consideration in these three subsets, which are studied in the forthcoming subsections. Concerning the behaviour of the period function near $\out_\nu$, the dichotomy between local regular value and local bifurcation value (see \defic{def2}) is solved for any $\np\in \Gamma_1\cup\Gamma_2\cup\Gamma_3$ thanks to the results in \cite{MMV03,MMV2,MMSV,Mariana}. In these papers it is computed the asymptotic expansion of the period function to second order, which usually suffices to tackle the regular/bifurcation dichotomy. However in order to study the criticality we need here to go further and compute the third, and even the fourth, order expansion. Let us advance that the asymptotic expansions for $\nu\in\Gamma_1$ are given in \propc{coeficients_Loud1} and the ones for $\nu\in\Gamma_2$ in \propc{prop52}. Being the proof of both results rather long and technical, for the sake of paper's readability we postpone them to \apc{appA}, where we also summarize the fundamental results and definitions from~\cite{MV20a,MV20b,MV21} that we shall use here. Among them we point out the notion of $L$-flatness $\F_L^\infty(\np_0)$, see \defic{defi2}, used in the remainder, and the \'Ecalle-Roussarie compensator $\omega(s;\alpha),$ that is a deformation of the logarithm used in the monomial scale in which the asymptotic expansion is given.

\begin{defi}\label{defi_comp}
The function defined for $s>0$ and $\alpha\in\R$ by means of
 \[
  \omega(s;\alpha)\!:=
  \left\{
   \begin{array}{ll}
    \frac{s^{-\alpha}-1}{\alpha} & \text{if $\alpha\neq 0,$}\\[2pt]
    -\log s & \text{if $\alpha=0,$}
   \end{array}
  \right.
 \]
is called the \emph{\'Ecalle-Roussarie compensator}. In the sequel we shall also use the  notation $\omega_\alpha(s)=\omega(s;\alpha)$.
\end{defi}
The asymptotic expansion for $\nu\in\Gamma_3$ is given in \propc{prop53} and its proof is of a different nature due to the occurrence of a saddle-node bifurcation at the polycycle.  

\subsection{Study of $\{D=-1/2, F\in(1/2,1)\}$ and $\{F=1/2, D\in(-1,0)\}$}\label{subsec_2.1}

\begin{figure}[t]
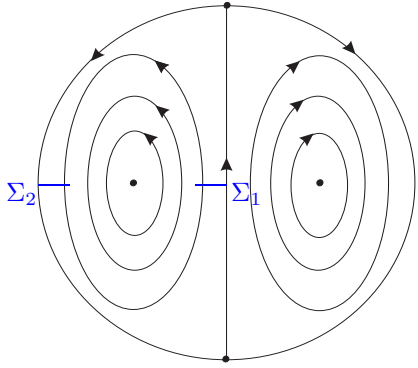

 \centering
 \begin{lpic}[l(0mm),r(0mm),t(0mm),b(0mm),]{fig1(0.8)}
  \lbl[l]{-5,28;$\blue{\Sigma_2}$}  
  \lbl[l]{32,28;$\blue{\Sigma_1}$}    
 \end{lpic}
 \caption{Phase portrait in the Poincaré disc of $X_\np$ for $\np\in V$, where for convenience the center at $(0,0)$ is shifted to the left and the vertical invariant line is $\{x=1\}$.}\label{fig1}
 \end{figure}

\figc{fig1} shows the phase portrait in the Poincaré disc of the vector field $X_\np$ in \refc{sist_loud} for $\np$ varying inside 
 \[
  V\!:=\left\{(D,F)\in\R^2: D\in (-1,0),\,F\in (0,1)\right\}.
 \]
We take transverse sections $\Sigma_1$ and $\Sigma_2$ parametrized by $s\mapsto (1-s,0)$ and $s\mapsto(-1/s,0)$ with $s>0$, respectively, and 
define $T(s;\np)$ to be the time that spends the solution of $X_\np$ with initial condition at $(1-s,0)\in\Sigma_1$ to arrive at $\Sigma_2.$ Thanks to the symmetry of $X_\np$ with respect to $\{y=0\}$, it turns out that 
the period of the periodic orbit passing through $(1-s,0)\in\Sigma_1$ is precisely $2T(s;\np)$. Consequently the emergence/disappearance of critical periodic orbits from $\out_\np$ corresponds to zeros of $T'(s;\np)$ bifurcating from $s= 0,$ more concretely to the number  $\mathcal Z_0(T'(\,\cdot\,;\np),\np_\star)$ as introduced in \defic{Z}. A key point to study these bifurcations is that $T(s;\np)$ is the Dulac time associated to the passage through a hyperbolic saddle, which is at infinity (see \figc{fig1} again). 
Therefore we can apply \cite[Theorem A]{MV20b} to obtain the asymptotic expansion of $T(s;\np)$ at $s=0$ and use then \cite[Theorem A]{MV21} to compute its first coefficients $T_{ij}(\np)$. Next result gathers all this information, where 
$\Gamma(\,\cdot\,)$ denotes the gamma function.

\begin{prop}\label{coeficients_Loud1}
Let $T(s;\np)$ be the Dulac time of the passage from $\Sigma_1$ to $\Sigma_2$ of the saddle at infinity of the vector field $X_\np$ in \refc{sist_loud} for $\np\in V$. Then the coefficients $T_{00},$ $T_{01},$ $T_{10}$ and $T_{20}$ in its asymptotic expansion at $s=0$ are meromorphic functions on $V$ that can be written as 
\[
\begin{array}{lcl}
T_{00}(\np)=
\frac{\pi}{2\sqrt{F(D+1)}}, & &
T_{01}(\np)
=\rho_1(\np)\frac{\Gamma\big(-\frac{\lambda}{2}\big)}{\Gamma\big(\frac{1-\lambda}{2}\big)},
\\[10pt]
T_{10}(\np)\textstyle=
\rho_2(\np)(2D+1)\frac{\Gamma\big(1-\frac{1}{2\lambda}\big)}{\Gamma\big(\frac{3}{2}-\frac{1}{2\lambda}\big)},
& &
T_{20}(\np)=\frac{\sqrt{\pi}}{\sqrt{2F}}\frac{\Gamma\big(\frac{1}{2}-\frac{1}{\lambda}\big)}{\Gamma\big(1-\frac{1}{\lambda}\big)}+\rho_3(\np)(2D+1),
\end{array}
\]
where $\lambda(\np)=\frac{F}{1-F}$ is the hyperbolicity ratio of the saddle,  
\[\textstyle
\rho_1(\np)=\frac{\sqrt{\pi}}{2(1-F)}\left(\frac{F}{D+1}\right)^{\frac{\lambda+1}{2}}\left(\frac{D}{F-1}\right)^{\frac{\lambda}{2}} \text{ and }\rho_2(\np)=\frac{\sqrt{\pi}}{2\sqrt{F(D+1)^3}},
\]
and $\rho_3$ is an analytic function on $V\cap \left\{\frac{2}{3}<F<1\right\}.$ In addition the following holds:

\begin{enumerate}[$(a)$]

\item If $\np_0\in V\cap\left\{\frac{2}{3}<F<1\right\}$ then, for all $\upsilon>0$ small enough,
\[
 T(s;\np)=T_{00}(\np)+T_{10}(\np)s+T_{20}(\np)s^2+\F_{L_0-\upsilon}^\infty(\np_0)
\]  
with $L_0=\min\big(3,\lambda(\np_0)\big)$. Moreover $T_{10}(-\frac{1}{2},F)=0$ and $T_{20}(-\frac{1}{2},F)> 0$ for all $F\in(\frac{2}{3},1)$.

\item If $\np_0\in V\cap\left\{\frac{1}{2}<F<\frac{2}{3}\right\}$ then, for all $\upsilon>0$ small enough, 
\[ 
 T(s;\np)=T_{00}(\np)+T_{10}(\np)s+T_{01}(\np)s^\lambda+\F^\infty_{2-\upsilon}(\np_0).
\]
Furthermore
$T_{10}(-\frac{1}{2},F)=0$ and $T_{01}(-\frac{1}{2},F)>0$ for all $F\in(\frac{1}{2},\frac{2}{3})$.

\item If $\np_0\in V\cap\left\{F=\frac{2}{3}\right\}$ then, for all $\upsilon>0$ small enough, 
\[
T(s;\np)=T_{00}(\np)+T_{10}(\np)s+T_{201}^2(\np)s^2\omega_{2-\lambda}(s)+T_{200}^2(\np)s^2+\F_{3-\upsilon}^\infty(\np_0),
\]
where $T_{200}^2$ and $T_{201}^2$ are analytic functions in a neighbourhood of $V\cap\{F=\frac{2}{3}\}$. Moreover 
$T_{10}(-\frac{1}{2},\frac{2}{3})=0$ and $T_{201}^2(-\frac{1}{2},\frac{2}{3})\neq 0$.

\item If $\np_0\in V\cap\left\{F=\frac{1}{2}\right\}$ then, for all $\upsilon>0$ small enough,
\[
 T(s;\np)=T_{00}(\np)+T_{101}^1(\np)s\omega_{1-\lambda}(s)+T_{100}^1(\np)s+\F_{2-\upsilon}^\infty(\np_0),
 \]
where
\[
T_{101}^1(\np)=-\rho_4(\np)(F-1/2)^2\text{ and }T_{100}^1(\np)=\rho_5(\np)(D+1/2)+\rho_6(\np)(F-1/2)
\]
for some analytic positive functions $\rho_i$ in a neighbourhood of $V\cap\{F=\frac{1}{2}\}$
with
$\rho_5(-\frac{1}{2},\frac{1}{2})=\rho_6(-\frac{1}{2},\frac{1}{2}).$
\end{enumerate}
\end{prop}

As we already explained, the proof of this result is postponed to \apc{appA}. The monomial order in each one of these asymptotic expansions is with respect to the strict partial order $\prec_{\np_0}$ given in \cite{MV20a}. Let us recall in its regard that we write $f\prec_{\np_0} g$ in case that 
 \[
  \lim_{(s,\np)\to (0,\np_0)}\frac{g(s;\lambda)}{f(s;\lambda)}=0.
 \]
For the monomials under consideration this order is preserved after derivation with respect to $s$, and so it is the good flatness properties of the remainder.  
Thus, as it occurs with the Taylor's series of an analytic function, an upper bound for the number of zeros of $T'(s;\np)$ that can bifurcate from $s=0$ follows by identifying the first non-vanishing coefficient in the asymptotic expansion. For the proof and a precise statement of this result, which essentially follows by using the well-known derivation-division algorithm, the reader is referred to \cite[Theorem C]{MV20a}.


\subsection{Study of  $\{F=2,D\in (-2,0)\}$ and  $\{D=\mathcal G(F):F\in(1,3/2)\}$}\label{subsec_2.2}

\begin{figure}[t]
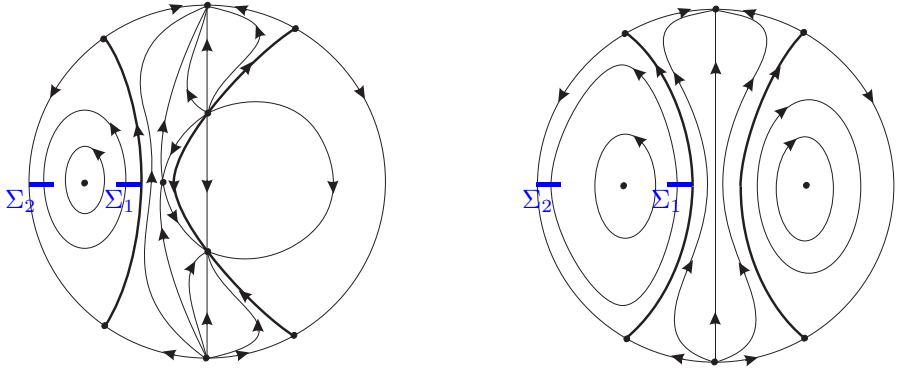

 \centering
 \begin{lpic}[l(0mm),r(0mm),t(0mm),b(0mm)]{dib3(1)}
  \lbl[l]{-3,22;$\blue{\Sigma_2}$}  
  \lbl[l]{10,22;$\blue{\Sigma_1}$}    
  \lbl[l]{65,22;$\blue{\Sigma_2}$}  
  \lbl[l]{82,22;$\blue{\Sigma_1}$}   
 \end{lpic}
 \caption{Phase portrait in the Poincaré disc of $X_\np$ for $\np\in W$ with $D<-1$ (left) and $D>-1$ (right), where the center at $(0,0)$ is shifted to the left,  the vertical invariant line is $\{x=1\}$ and the hyperbola $\{\frac{1}{2}y^2-q(x)=0\}$ appears in boldface type.}\label{fig4}
 \end{figure}

\figc{fig4} shows the phase portrait in the Poincaré disc of the vector field $X_\np$ in \refc{sist_loud} for $\np$ varying inside
\[
  W\!:=\{(D,F)\in\R^2:F+D>0,\,D<0\text{ and }F>1\}.
\]
In this case the outer boundary of the period annulus of the center at $(0,0)$ is 
 contained in the union of the line at infinity and an invariant hyperbola $\mathcal C\!:=\{\frac{1}{2}y^2-q(x)=0\}$, where $q(x)=ax^2+bx+c$ with
\begin{equation*}
a\!:=\frac{D}{2(1-F)},\mbox{ }
b\!:=\frac{D-F+1}{(1-F)(1-2F)}\mbox{ and }
c\!:=\frac{F-D-1}{2F(1-F)(1-2F)}.
\end{equation*}
One can verify that if $\np\in W$ then $q$ has two distinct real zeros, that we shall denote by $p_1$ and $p_2$ taking $p_1<p_2.$ 
That being said, we place two transverse sections $\Sigma_1$ and $\Sigma_2$ parametrized by $s\mapsto (p_1-s,0)$ and $s\mapsto (-1/s,0)$ with $s>0$, respectively, and define $T(s;\np)$ to be the time that takes to the solution of~$X_\np$ with initial condition at $(p_1-s,0)\in\Sigma_1$ to arrive at $\Sigma_2.$ Then $T(s;\np)$ is the Dulac time associated to the passage through a hyperbolic saddle at infinity, so that we can apply the results in \cite{MV20b,MV21} to obtain its asymptotic expansion at $s=0.$ This is important for the proof of \teoc{Loud} because, exactly as in the previous case, the symmetry of $X_\np$ with respect to $\{y=0\}$ implies that 
the period of the periodic orbit passing through $(1-s,0)\in\Sigma_1$ is $2T(s;\np)$. With regard to our next result we remark that $\frac{1-p_2}{1-p_1}<1$ for all $\np\in W$, which is relevant since the hypergeometric function  
${}_{2}F_1(a,b;c;\,\cdot\,)$ is holomorphic on $\C\setminus [1,+\infty)$, see Appendix~\ref{apB}. 
Let us also mention that $\mathrm B(\cdot,\cdot)$ is the beta function.
\begin{prop}\label{prop52}
Let $T(s;\np)$ be the Dulac time of the passage from $\Sigma_1$ to $\Sigma_2$ of the saddle at infinity of the vector field $X_\np$ in \refc{sist_loud} for $\np\in W$. Then the coefficients $T_{00},$ $T_{01},$ $T_{10}$ and $T_{20}$ in its asymptotic expansion at $s=0$ are meromorphic functions on $W$ that can be written as 
\begin{align*}
T_{00}(\np)
&\textstyle=\frac{\sqrt{2}}{\sqrt{a}(1-p_1)}\,{}_2F_1\big(1,\frac{1}{2};\frac{3}{2};\frac{1-p_2}{1-p_1}\big),\\[4pt]
T_{01}(\np)&\textstyle=\rho_1(\np)\mathrm B\big(-\lambda,\frac{1}{2}\big),\\[4pt]
T_{10}(\np)&\textstyle=\rho_2(\np)\mathrm B\big(1-\frac{1}{\lambda},-\frac{1}{2}\big){}_2F_1\big(-1-\frac{1}{\lambda},-\frac{1}{2};\frac{1}{2}-\frac{1}{\lambda};\frac{1-p_2}{1-p_1}\big)
\intertext{and}
T_{20}(\np)&\textstyle=\rho_3(\np)\mathrm B\big(1-\frac{2}{\lambda},-\frac{3}{2}\big){}_2F_1\big(-\frac{2}{\lambda}-3,-\frac{3}{2};-\frac{1}{2}-\frac{2}{\lambda};\frac{1-p_2}{1-p_1}\big)+\rho_4(\np)T_{10}(\np),
\end{align*}
where $\lambda(\np)=\frac{1}{2(F-1)}$ is the hyperbolicity ratio of the saddle and, for $i=1,2,3,4$, $\rho_i$ is an analytic positive function on $W$. In addition the following holds:
\begin{enumerate}[$(a)$]

\item If $\np_0\in W\cap\left\{1<F<\frac{5}{4}\right\}$ then, for all $\upsilon>0$ small enough, 
          \[
          T(s;\np)=T_{00}(\np)+T_{10}(\np)s+T_{20}(\np)s^2+\F_{L_0-\upsilon}^\infty(\np_0)
          \]
          with $L_0=\min(3,\lambda(\np_0))$. Moreover $T_{20}(\np)\neq 0$ for all 
          $\np\in W\cap\left\{1<F<\frac{5}{4}\right\}$ such that $T_{10}(\np)=0.$
         
\item If $\np_0\in W\cap\left\{\frac{5}{4}<F<\frac{3}{2}\right\}$ then, for all $\upsilon>0$ small enough,
         \[
          T(s;\np)=T_{00}(\np)+T_{10}(\np)s+T_{01}(\np)s^\lambda+T_{20}(\np)s^2+\F_{L_0-\upsilon}^\infty(\np_0)
          \]
          with $L_0=\lambda(\np_0)+1,$ and there exists a unique 
          $\np_\star\in W\cap\left\{\frac{5}{4}<F<\frac{3}{2}\right\}$ such that $T_{10}(\np_\star)=0$ and 
          $T_{01}(\np_\star)=0$. 
         Furthermore
          $T_{20}(\np_\star)<0$, the gradients of $T_{01}$ and $T_{10}$ at $\np_\star$ are linearly independent, and
          $\np_\star=(D_\star,\frac{4}{3})$ with $D_\star=\mathcal G(\frac{4}{3})\approx -1.128$.
           
\item If $\np_0\in W\cap\left\{F=\frac{5}{4}\right\}$ then, for all $\upsilon>0$ small enough,
         \[
          T(s;\np)=T_{00}(\np)+T_{10}(\np)s+T_{201}^2(\np)s^2\omega_{2-\lambda}(s)
             +T^2_{200}(\np)s^2+\F_{3-\upsilon}^\infty(\np_0),
         \]
        where $T_{200}^2$ and $T_{201}^2$ are analytic functions in a neighbourhood of  
       $W\cap\{F=\frac{5}{4}\}$. Moreover        
        $T_{10}(D,\frac{5}{4})=0$ if and only if $D=-1$, and $T_{201}^2(-1,\frac{5}{4})\neq 0$.
        
\item If $\np_0\in W\cap\left\{F=2\right\}$ then, for all $\upsilon>0$ small enough,
        \[
        T(s;\np)=T_{00}(\np)+T_{01}(\np)s^\lambda+T^{\frac{1}{2}}_{101}(\np)s\omega_{1-2\lambda}(s)
        +T_{100}^{\frac{1}{2}}(\np)s+\F_{3/2-\upsilon}^\infty(\np_0),
        \]
       where $T_{100}^{\frac{1}{2}}$ and $T_{101}^{\frac{1}{2}}$ are analytic functions in a neighbourhood of $W\cap\{F=2\}$.
       Moreover $T_{01}(D,2)=0$ for all $D\in (-2,0)$, $T_{101}^{\frac{1}{2}}(D,2)=0$ if and only if $D=-\frac{1}{2}$, and 
       the gradients of $T_{01}$ and $T_{101}^{\frac{1}{2}}$ are linearly independent at $(-\frac{1}{2},2)$.
       
\end{enumerate}
\end{prop}

The proof of this result is postponed to \apc{appA}. 

\begin{obs}\label{fun}
The asymptotic expansions in \propc{prop52} were already given in \cite[Theorem 3.6]{MMV2} but only to second order. In that result it is given, among others, the expression of the coefficient $T_{10}(\np)$ in terms of a definite improper integral. Furthermore, see \cite[Proposition 3.11]{MMV2}, it is proved by applying the Implicit Function Theorem that the set of those $\np\in W_1\!:=W\cap\{F<3/2\}$ such that $T_{10}(\np)=0$ is the graphic of an analytic function $D=\mathcal G(F).$
This is the function that appears in assertion $(b)$ of \propc{prop52}. Thanks to the results in \apc{ApBeta} we can now identify the improper integral as a hypergeometric function, so that we can write
\[\textstyle
  \big\{\np\in W_1:D=\mathcal G(F)\big\}
  =\big\{\np\in W_1:
  {}_2F_1\big(-1-\frac{1}{\lambda},-\frac{1}{2};\frac{1}{2}-\frac{1}{\lambda};\frac{1-p_2}{1-p_1}\big)=0
  \big\},
\]
where $p_1$ and $p_2$ with $p_1<p_2$ are the real roots of $q(x)=0$ and $\lambda(\np)=\frac{1}{2(F-1)}$.
\end{obs}

\begin{obs}
In the statement of \propc{prop52} we refer to some positive functions $\rho_i\in\cc^\omega(W).$ Let us mention that in the proof we show that  
\[
\begin{array}{lcl}
 \rho_1(\np)=\frac{1}{2\sqrt{2a}}\frac{(p_2-p_1)^{\frac{1}{2(F-1)}}}{(F-1)(1-p_1)^{\frac{F}{F-1}}}
& &
 \rho_2(\np)=\frac{1}{2\sqrt{2a}}\frac{1}{(p_2-p_1)(1-p_1)} \\[15pt]
 \rho_3(\np)=\frac{3}{8\sqrt{2a}}\frac{1}{(p_2-p_1)^2(1-p_1)} & &
 \rho_4(\np)=\frac{p_1-1+2F(p_2-p_1)}{(p_2-p_1)(p_1-1)}
\end{array}
\]
We do not use the explicit expressions in this paper but they may be relevant for future applications.
\end{obs}


\subsection{Study of $\{F=1, D\in(-1,0)\}$}\label{subsec_2.3}

Our aim in this section is to study the period function of the center at the origin of $X_\np$ for $\np=(D,F)$ with $F\approx 1$ and $D\in(-1,0)$. 
\begin{figure}[t]
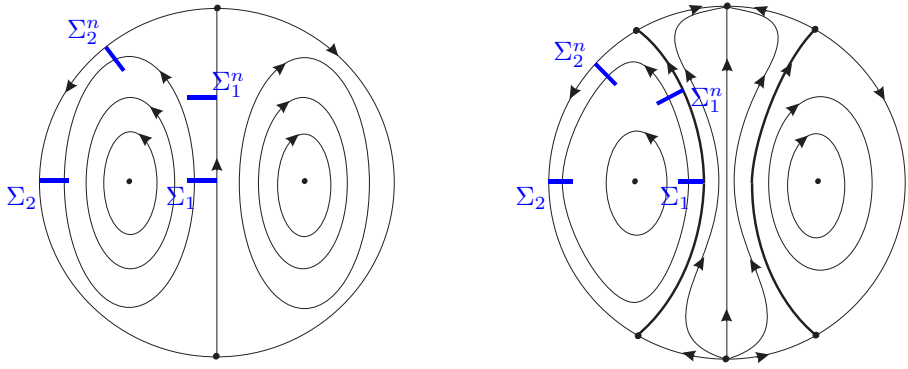

 \centering
 \begin{lpic}[l(0mm),r(0mm),t(0mm),b(0mm)]{dib4(1)}
  \lbl[l]{-4,22;$\blue{\Sigma_2}$}
  \lbl[l]{4,44;$\blue{\Sigma_2^n}$}
  \lbl[l]{23,37;$\blue{\Sigma_1^n}$}
  \lbl[l]{68,41;$\blue{\Sigma_2^n}$}
  \lbl[l]{86,35;$\blue{\Sigma_1^n}$}           
  \lbl[l]{17,22;$\blue{\Sigma_1}$}    
  \lbl[l]{63,22;$\blue{\Sigma_2}$}  
  \lbl[l]{82,22;$\blue{\Sigma_1}$}   
 \end{lpic}
 \caption{Phase portrait of $X_\np$ for $\np=(D,F)\in (-1,0)\times (0,+\infty)$ in the Poincaré disc with $F\leqslant 1$ (left) and $F>1$ (right). In this case, contrary to the previous ones, the singularity at infinity for $F=1$ is not a hyperbolic saddle but a saddle-node.}\label{fig6}
 \end{figure}
To this end we introduce transverse sections $\Sigma_1$ and $\Sigma_2$ parametrized, respectively, by means of
\[
 \sigma_1(s;\np)\!:=\left\{
  \begin{array}{cl}
   (1-s,0) & \text{ if $F\leqslant 1,$} \\[2pt]
   (p_1-s,0) & \text{ if $F> 1,$}
  \end{array}
 \right.
 \text{ and }
 \sigma_2(s;\np)\!:=(-1/s,0),
\]
where recall that $q(x)=a(x-p_1)(x-p_2)$ with $p_1<p_2$ for $F>1$. One can also check that $\lim_{F\to 1^+}p_1=1.$ For each $\np=(D,F)\in (-1,0)\times (0,+\infty)$ we define $T(s;\np)$ as the time that spends the solution of $X_\np$ starting at $\sigma_1(s;\np)\in\Sigma_1$ to arrive at $\Sigma_2.$ 
A key feature of this Dulac time is that the singularity for $F=1$ is not a hyperbolic saddle but a saddle-node. Our next result gives the asymptotic expansion of $T(s;\np)$ at $s=0$ for $F\approx 1$. We point out that this is relevant for the proof of \teoc{Loud} because the period of the periodic orbit of $X_\np$ passing through $\sigma_1(s;\np)\in\Sigma_1$ is precisely $2T(s;\np)$ due to the symmetry of the vector field.

\begin{prop}\label{prop53}
Let $T(s;\np)$ be the Dulac time of the passage from $\Sigma_1$ to $\Sigma_2$ of the saddle-node unfolding at infinity $\{X_\np\}.$  
Then there is an open neighbourhood~$\mathscr U$ of~$(-1,0)\times\{1\}$ such that
\begin{equation*}
T(s;\np)=T_{0}(\np)+T_{1}(\np) s+T_{2}(\np)s^2+s^2h(s;\np),
\end{equation*}
where $T_0,T_1,T_2\in\cc^0(\mathscr U)$ and, setting $\Theta=s\partial_s$, $\lim_{s\to 0^+}\Theta^kh(s;\np)=0$ uniformly on compact sets of $\mathscr U$ for $k=0,1,2$. Moreover $T_{1}(D,1)=0$ if, and only if, $D=-\frac{1}{2}$. Finally $T_{2}(-1/2,1)\neq 0.$
\end{prop}

\begin{prova}
To study the saddle-node bifurcation that occurs at infinity we work in the projective plane $\RP^2$ and perform the change of coordinates 
\[
 (u,v)=p(x,y)\!:=\left(\frac{1-x}{y},\frac{1}{y}\right).
\] 
The meromorphic extension of $X_\np$ in these coordinates is given by
\[
\tilde X_\np\!:=p_*X_\np=\frac{1}{v}\Big(uP(u,v;\np)\partial_u+vQ(u,v;\np)\partial_v\Big)
\]  
with
\begin{align*}
P(u,v;\np)&=1-F-Du^{2}+(2D+1)uv-(D+1)v^{2}\\
\intertext{and}
Q(u,v;\np)&=-F-Du^{2}+(2D+1)uv-(D+1)v^{2}.
\end{align*}
Our first goal is to show that we can bring locally the saddle-node unfolding to a convenient normal form in order that we can apply the tools developed in \cite{MMSV} to study the asymptotic expansion of its Dulac map and Dulac time. With this aim, 
some long but easy computations show that the local analytic change of coordinates given by
\begin{equation*}
(z,w)=\Psi(u,v)\!:=\left(\frac{u}{\sqrt{g(u,v)}},\frac{v}{\sqrt{g(u,v)}}\right),
\end{equation*}
where $g(u,v)\!:= {\frac { \left( 2D+1\right) }{  \left( 2F-1 \right)D }}uv-{\frac 
{ \left( D+1 \right) }{  2  F D}}{v}^{2}-\frac{1}{2D},$ brings the vector field $\tilde X_\np$ to
\[
\bar X_\np\!:=\frac{1}{w\bar U(z,w;\np)}\bigl(z(z^{2}-2(F-1))\partial_{z}-w(2F-z^{2})\partial_{w}\bigr),
\]
with $\bar U(z,w;\np)\!:=\left({\frac { \left( 2D+1 \right)}{2(2F-1)}} zw-{\frac {
\left( D+1 \right) }{4F}}{w}^{2}-\frac{D}{2}\right)^{-\frac{1}{2}}$.
A technical assumption in order to apply the results from \cite{MMSV} is that for each $\np$ the Taylor's series of $(z,w)\mapsto\bar U(z,w;\np)$ at $(0,0)$ is absolutely convergent for all $(z,w)\in [-1,1]^2$. This is not fulfilled unless we perform a rescaling which is only well defined provided that~$\np$ varies inside a compact subset of $(-1,0)\times (0,+\infty)$ and this forces us to work locally. For this reason, as a first step in the proof, we will show a local version of the statement. More concretely, that for each $\np_\star=(D_\star,1)$ with $D_\star\in (-1,0)$ there exists an open ball $B_{\np_\star}=\{\np\in\R^2:\|\np-\np_\star\|<\delta\}$ 
such that 
\begin{equation*}
T(s;\np)=T_0^{\np_\star}(\np)+T_{1}^{\np_\star}(\np) s+T_{2}^{\np_\star}(\np)s^2+s^2h^{\np_\star}(s;\np),
\end{equation*}
with $T_i^{\np_\star}$ continuous functions on $B_{\np_\star}$ and $\lim_{s\to 0^+}\Theta^kh^{\np_\star}(s;\nu)=0$ uniformly on~$B_{\np_\star}$ for $k=0,1,2$. To begin with we take $\delta>0$ small enough so that the closure of $B_{\np_\star}$ is inside $(-1,0)\times (\frac{1}{2},+\infty)$ and define 
\begin{equation}\label{p3eq7}
 r\!:=\inf\left\{\left|\frac{(D+1)}{DF}w^2-\frac{2(2D+1)}{D(2F-1)}zw\right|^{-1/2}:|z|\leqslant 1, |w|\leqslant 1
 \text{ and } \|\np-\np_\star\|\leqslant\delta\right\},
\end{equation}
which is strictly positive.
The pull-back of $\bar X_\np$ by the rescaling $\rho(z,w)\!:=(rz,rw)$ can now be written as in \cite[Eq. 13]{MMSV} because
one can easily verify that
\[
 \rho^*\bar X_\np=\frac{1}{wU(z,w;\np)}\left(z(z^2-\varepsilon)\partial_z-w(2F/r^2-z^2)\partial_w\right)
 \text{ with $\varepsilon\!:=2(F-1)/r^2$}
\]
and where the Taylor's series of
\[
U(z,w;\np)\!:=\frac{\bar U(rz,rw;\np)}{r}
=\frac{-D}{2r}\left(1+r^2\left(\frac{(D+1)}{2FD}w^2-\frac{(2D+1)}{D(2F-1)}zw\right)\right)^{-\frac{1}{2}}
\]
at $(z,w)=(0,0)$ is absolutely convergent for all $(z,w)\in [-1,1]^2$ and $\np\in B_{\np_\star}$ since, on account of \refc{p3eq7},
\[ 
r^2\left|\frac{(D+1)}{2FD}w^2-\frac{(2D+1)}{D(2F-1)}zw\right|\leqslant\frac{1}{2}\text{ for all $(z,w,\np)\in[-1,1]^2\times B_{\np_\star}.$}
\]
In these new rescaled coordinates, that we still denote by $(z,w)$ for simplicity, the period annulus is inside the quadrant $\{w\geqslant 0,z\geqslant\vartheta_{\varepsilon}\}$ where 
\[
 \vartheta_{\varepsilon}\!:=\left\{
  \begin{array}{cl}
   \sqrt{\varepsilon} & \text{ if $\varepsilon\geqslant 0$},\\
   0 & \text{ if $\varepsilon<0$}.
  \end{array}
 \right.
\]
Setting $\Psi_\star\!:=\rho^{-1}\circ\Psi,$
we take two auxiliary transverse sections $\Sigma_1^n\!:=\Psi_\star^{-1}(\{w=1\})$ and $\Sigma_2^n\!:=\Psi_\star^{-1}
(\{z=1\})$ parameterized by $\sigma_1^n(s;\np)\!:=\Psi_\star^{-1}(s+\vartheta_{\varepsilon},1)$ and $\sigma_2^n(s;\np)\!:=\Psi_\star^{-1}
(1,s)$, respectively (see \figc{fig6}). We define $\mathcal T(s;\np)$ and $\mathscr D(s;\np)$ to be the Dulac time and Dulac map of $\tilde X_\np$ from $\Sigma_1^n$ to $\Sigma_2^n,$ respectively. We remark that, by construction, $\mathcal T(s;\np)$ is the time that the solution of $\rho^*\bar X_\np$ starting at the point $(s+\vartheta_{\varepsilon},1)$ spends to arrive at $\{z=1\}$ and that the intersection point is precisely $(z,w)=(1,\mathscr D(s;\np))$. In this regard, since $w=\mathscr D(z;\np)$ is a trajectory of the vector field $z(z^{2}-\varepsilon)\partial_{z}+w(2F/r^2-z^{2})\partial_{w}$, see \cite[p. 6417]{MMSV} for details, the application of $(b)$ in Corollary~A of \cite{MMSV} with $\{\mu=2,\ell=k=2,\lambda=2F/r^2\}$ shows that 
\begin{equation}\label{p3eq2}
 \mathscr D(s;\np)=s^2\mathcal I(B_{\np_\star}),
\end{equation}
by shrinking $\delta>0$ if necessary. Here, and in what follows, $\mathcal I(B_{\np_\star})$ stands for some function $h(s;\np)$ verifying that $\lim_{s\to 0^+}\Theta^kh(s;\np)=0$ uniformly on $\np\in B_{\np_\star}$ for $k=0,1,2$.
Furthermore, by applying Corollary B in the same paper with $\{\mu=2,\ell=k=2\}$ and shrinking $\delta>0$ again we can assert that
\begin{equation}\label{p3eq3}
 \mathcal T(s;\np)=b_0(\np)+b_1(\np)s+b_2(\np)s^2+s^2\mathcal I(B_{\np_\star})
\end{equation}
with $b_i\in\mathscr C^0(B_{\np_\star})$ for $i=0,1,2.$ 
Working in the original $(x,y)$ coordinates, we consider next the transition times $T_1(\,\cdot\,;\np)$ and $T_2(\,\cdot\,;\np)$ of~$X_\np$ from $\Sigma_1$ to $p^{-1}(\Sigma_1^n)$ and from $p^{-1}(\Sigma_2^n)$ to $\Sigma_2$, respectively. We define moreover $R(\,\cdot\,;\np)$ to be the transition map from~$\Sigma_1$ to~$p^{-1}(\Sigma_1^n)$. Accordingly
\begin{equation}\label{p3eq4}
 T(s)=T_1(s)+\big(\mathcal T\circ R\big)(s)+\big(T_2\circ\mathscr D\circ R\big)(s),
\end{equation}
where we omit the dependence on $\np$ for the sake of shortness. By \cite[Lemma 3.2]{MMV03}, we have that $T_2(s;\np)$ an analytic function at $\{0\}\times B_{\np_\star}$ with $T_2(0;\np)=0.$ Observe at this point that, setting
\[
 \xi_\np\!:=\left\{
  \begin{array}{cl}
   0 & \text{ if $F\leqslant 1,$} \\
   1-p_1 & \text{ if $F>1,$}
  \end{array}
 \right.
\]
we can write the parametrization of $\Sigma_1$ as $\sigma_1(s;\np)=(1-\xi_\np-s,0).$ We claim that there exist two functions $f(\hat s;\np)$ and $g(\hat s;\np)$, analytic at $\{0\}\times B_{\np_\star},$ such that
\begin{equation}\label{p3eq1}
 T_1(s;\np)=f(s+\xi_\np;\np)\text{ and }R(s;\np)=g(s+\xi_\np;\np)-\vartheta_\varepsilon.
\end{equation}
To show this let us consider two additional transverse sections $\hat\Sigma_1$ and $\hat\Sigma_1^n$ parameterized respectively by $\hat\sigma_1(\hat s;\np)\!:=(1-\hat s,0)$ and $\hat\sigma_1^n(\hat s;\np)\!:=(\Psi_\star\circ p)^{-1}(\hat s,1),$ which clearly are analytic at $\{0\}\times B_{\np_\star}$. Moreover it is clear that they are related with $\Sigma_1$ and $\Sigma_1^n$ through $\sigma_1(s;\np)=\hat\sigma_1(s+\xi_\np;\np)$ and $\sigma_1^n(s;\np)=\hat\sigma_1^n(s+\vartheta_\varepsilon;\np).$ That being said, the claim follows by noting that if we choose $f(\hat s;\np)$ and $g(\hat s;\np)$ to be, respectively, the transition time and transition map of $X_\np$ from $\hat\Sigma_1$ to $\hat\Sigma_1^n$, which are clearly analytic at $\{0\}\times B_{\np_\star},$ then the equalities in~\refc{p3eq1} hold. Note moreover that $g(\xi_\np;\np)=\vartheta_\varepsilon$ since $R(0;\np)=0.$ On account of the claim, by considering the second order Taylor's development of $f(\hat s;\np)$ and $g(\hat s;\np)$ at $\hat s=\xi_\np$, respectively, we get 
\[
 T_1(s;\np)=a_0(\np)+a_1(\np)s+a_2(\np)s^2+s^2\mathcal I(B_{\np_\star})\text{ and }R(s;\np)=c_1(\np)s+c_2(\np)s^2+s^2\mathcal I(B_{\np_\star})
\]
with $a_i,c_i\in\mathscr C^0(B_{\np_\star})$ and where we also use that $\np\mapsto \xi_\np$ is a continuous function. The combination of the second expression above with \refc{p3eq2} and \refc{p3eq3} yields
\[
 \big(\mathscr D\circ R\big)(s)=s^2\mathcal I(B_{\np_\star})\text{ and }
 \big(\mathcal T\circ R\big)(s)=\hat b_0(\np)+\hat b_1(\np)s+\hat b_2(\np)s^2+s^2\mathcal I(B_{\np_\star}),
\]
respectively, with $\hat b_i\in\mathscr C^0(B_{\np_\star}).$ Summing up, since $\big(T_2\circ\mathscr D\circ R\big)(s)=s^2\mathcal I(B_{\np_\star})$ due to $T_2(0;\np)=0,$  from~\refc{p3eq4} we can assert that 
\begin{equation}\label{p3eq5}
T(s;\np)=T_0^{\np_\star}(\np)+T_1^{\np_\star}(\np)s+T_2^{\np_\star}(\np)s^2+s^2h^{\np_\star}(s;\np)
\end{equation}
for some functions $T^{\np_\star}_i$ that are continuous on $B_{\np_\star}$ and some 
$h^{\np_\star}\in\mathcal I(B_{\np_\star})$. This concludes the proof of the local version of the statement, in which we remark that the coefficients $T_i^{\np_\star}(\np)$ and the remainder $s^2h^{\np_\star}(s;\np)$ depend by construction on $\np_\star.$
Our next step will be to globalize them and to this end we define 
\[
 \mathscr U\!:=\bigcup_{\np_\star\in (-1,0)\times\{1\}}B_{\np_\star}
\]
which is clearly an open neighbourhood of $(-1,0)\times\{1\}$. Let us consider now any 
$\np_1,\np_2\in (-1,0)\times\{1\}$ such that $B_{\np_1}\cap B_{\np_2}\neq 0.$ Then, from \refc{p3eq5}, we get that
\[
 T_0^{\np_1}(\np)-T_0^{\np_2}(\np)+\big(T_1^{\np_1}(\np)-T_1^{\np_2}(\np)\big) s+
 \big(T_2^{\np_1}(\np)-T_2^{\np_2}(\np)\big) s^2+s^2\big(h^{\np_1}(s;\np)-h^{\np_2}(s;\np)\big)=0
\]
for all $s>0$ small enough and $\np\in B_{\np_1}\cap B_{\np_2}$. Since $h^{\np_1}-h^{\np_2}\in\mathcal I(B_{\np_1}\cap B_{\np_2}),$ taking the limit $s\to 0^+$ on both sides of the above equality we deduce that 
$T_0^{\np_1}=T_0^{\np_2}$ on $B_{\np_1}\cap B_{\np_2}$. Similarly, but taking the first and second derivatives with respect to $s$, respectively, we get that $T_1^{\np_1}=T_1^{\np_2}$ and $T_2^{\np_1}=T_2^{\np_2}$ on $B_{\np_1}\cap B_{\np_2}$. Hence, for $i=0,1,2,$ the local functions $T_i^{\np_\star}\in\mathscr C^0(B_{\np_\star})$ for $\np_\star\in 
(-1,0)\times\{1\}$ glue together into a well defined continuous function $T_i$ on~$\mathscr U$. Exactly the same argument shows that the local functions 
$h^{\np_\star}\in\mathcal I(B_{\np_\star})$ for $\np_\star\in 
(-1,0)\times\{1\}$ glue together into a well defined function $h(s;\np)$ satisfying that $\lim_{s\to 0^+}\Theta^kh(s;\nu)=0$ uniformly on compact sets of $\mathscr U$ for $k=0,1,2$. To show this last assertion it suffices to take a finite subcover $B_{\np_1}\cup\ldots\cup B_{\np_n}$ of the given compact subset of $\mathscr U$ and use that $h|_{B_{\np_i}}\in\mathcal I(B_{\np_i})$ for $i=1,2,\ldots,n.$ 

So far we have proved the first assertion in the statement. Let us turn to the proof of the second one. To this end the key point is that for those $\np_0\in \mathscr U\cap\{F<1\}$ we can also apply $(a)$ in \propc{coeficients_Loud1}
to obtain that
\begin{equation}\label{p3eq6}
 T(s;\np)=T_{00}(\np)+T_{10}(\np)s+T_{20}(\np)s^2+\F_{3-\upsilon}^\infty(\np_0),
\end{equation}
where, setting $\lambda(\np)=\frac{F}{1-F}$, 
\[\textstyle
 T_{10}(\np)=\frac{\sqrt{\pi}(2D+1)}{2\sqrt{F(1+D)^3}}
\frac{\Gamma\left(1-\frac{1}{2\lambda}\right)}{\Gamma\left(\frac{3}{2}-\frac{1}{2\lambda}\right)}
\text{ and }T_{20}(-\frac{1}{2},F)=\frac{\sqrt{\pi}}{\sqrt{2F}}\frac{\Gamma\big(\frac{1}{2}-\frac{1}{\lambda}\big)}{\Gamma\big(1-\frac{1}{\lambda}\big)}.
\]
Hence, since $T_i\in\mathscr C^0(\mathscr U)$, from~\refc{p3eq5} and~\refc{p3eq6} we can assert that
\[
 T_1(D,1)=\lim_{F\to 1^-} T_{10}(D,F)=\frac{2D+1}{(1+D)^{3/2}},
\]
where we also use that $\lim_{F\to 1^-}\frac{\Gamma\left(1-\frac{1}{2\lambda}\right)}{\Gamma\left(\frac{3}{2}-\frac{1}{2\lambda}\right)}=\frac{\Gamma(1)}{\Gamma(\frac{3}{2})}=\frac{2}{\sqrt{\pi}}$. Consequently, as desired, $T_1(D,1)=0$ if and only if $D=-\frac{1}{2}$. The same argument shows that 
\[
 T_2(-1/2,1)=\lim_{F\to 1^-} T_{20}(-1/2,F)=\frac{\sqrt{\pi}}{\sqrt{2}}\frac{\Gamma\big(\frac{1}{2}\big)}{\Gamma\big(1\big)}\neq 0,
\]
and this completes the proof of the result. 
\end{prova}

\section{Distinguished cases}\label{distinguished}

This section is devoted to study three specific parameters. Recall that among the quadratic centers there are four nonlinear isochrones, see \refc{isochronous}. Chicone and Jacobs show in \cite[Theorem 3.1]{Chicone} that the criticality of the period function at the inner boundary (i.e., the center) of  $\PA$ is exactly 1 for each one of the nonlinear isochrones. In this section we prove that for two of them, namely $\nu=(-\frac{1}{2},2)$ and $\nu=(-\frac{1}{2},\frac{1}{2}),$
the criticality at the outer boundary (i.e., the polycycle) is also 1, see Propositions~\ref{I2} and~\ref{I4}, respectively. 
In the same vein it is also well-known that the criticality at the inner boundary of any quadratic center is at most two, see \cite[Theorem 3.2]{Chicone}. This maximum criticality is achieved in three parameter values, the so-called Loud points, which following the notation in \cite{Chicone} are given by $\nu=L_i$ with
 \begin{equation}\label{Li}
 \textstyle
  L_1\!:=\left(-\frac{3}{2},\frac{5}{2}\right),\;
  L_2\!:=\left(\frac{-11+\sqrt{105}}{20},\frac{15-\sqrt{105}}{20}\right)\text{ and }
  L_3\!:=\left(\frac{-11-\sqrt{105}}{20},\frac{15+\sqrt{105}}{20}\right).
 \end{equation}
As we already explained in the introduction, we conjecture that the criticality at the outer boundary of any quadratic center is at most two, and that there are only three parameter values where this maximum criticality is attained. In this paper we identify and prove the validity of the conjecture for  two of these parameters. We investigate one of them in this section, see \propc{doble}. The other one was already studied in \cite{MMV20} and we postpone its treatment until the proof of \teoc{Loud}. 

The following is a sort of division theorem within the class of flat functions that will be used to study the criticality at the outer boundary for the above-mentioned isochrones. In its statement, and in what follows, we use the notation $0_n=(0,0,\ldots,0)\in\R^n$ for the sake of shortness.

\begin{lem}\label{Fdiv}
Let us fix $K\in\N\cup\{\infty\}$, $L\geqslant 0$ and $n\in\N.$ If $f(s;\mu_1,\ldots,\mu_n)\in\F^K_L(0_n)$ verifies that 
\[
 f(s;\mu_1,\ldots,\mu_{k-1},0,\ldots,0)\equiv 0\text{, for some $k\in\{1,2,\ldots,n\}$,}
\]
then there exist $f_{k},\ldots,f_n\in\F_L^{K-1}(0_n)$ such that $f=\sum_{i=k}^n\mu_if_i$.
\end{lem}

\begin{prova}
We proceed by induction on $n\in\N.$ For the base case $n=1$ we take $f(s;\mu_1)\in\F_L^K(0_1)$ with $f(s;0)\equiv 0$ and define $f_1(s;\mu_1)\!:=\int_0^1\partial_2f(s;\mu_1t)dt,$ so that $f=\mu_1 f_1.$ To show that $f_1\in \F_L^{K-1}(0_1)$  we use that, by hypothesis (see \defic{defi2}), for every $\nu=(\nu_0,\nu_1)\in\Z^2_{\ge 0}$ with $|\nu|=\nu_0+\nu_1\leqslant K-1$ there exist a neighbourhood $V\subset\R$ of $0$ and $C,s_0>0$ such that $|\partial^{\nu_0}_s\partial^{\nu_1+1}_{\mu_1}f(s;\mu_1)|\leqslant Cs^{L-\nu_0}$ for every $\mu_1\in V$ and $s\in (0,s_0)$. On account of this and applying the Dominated Convergence Theorem \cite[Theorem~11.30]{Rudin},
\[
 |\partial^\nu f_1(s;\mu_1)|
 \leqslant\int_0^1\left|\partial^\nu(\partial_{2}f(s;\mu_1t))\right|dt
 \leqslant\int_0^1|\partial_s^{\nu_0}\partial_{2}^{\nu_1+1}f(s;\mu_1t)|t^{\nu_1}dt
 \leqslant \frac{C}{\nu_1+1}s^{L-\nu_0}
\]
for every $\mu_1\in V$ and $s\in (0,s_0).$ Hence $f_1\in\F_L^{K-1}(0_1)$. To prove the inductive step we suppose that $n>1$ and consider $f(s;\mu_1,\ldots,\mu_n)\in\F^K_L(0_n)$ verifying that 
$f(s;\mu_1,\ldots,\mu_{k-1},0,\ldots,0)\equiv 0$ for some $k\in\{1,2,\ldots,n\}$. It is clear that we can write
\begin{align}\label{l1eq1}
&f(s;\mu_1,\ldots,\mu_{n-1},\mu_n)-f(s;\mu_1,\ldots,\mu_{n-1},0)=\mu_nf_n(s;\mu_1,\ldots,\mu_n)\\
\intertext{with}\notag
&f_n(s;\mu_1,\ldots,\mu_n)\!:=\int_0^1\partial_{n+1}f(s;\mu_1,\ldots,\mu_{n-1},\mu_nt)dt.
\end{align} 
Similarly as for the base case, taking $f\in\F_L^{K}(0_n)$ into account,  one can easily show that $f_n\in\F_L^{K-1}(0_n)$. Since $f(s;\mu_1,\ldots,\mu_{n-1},0)|_{\mu_k=\ldots=\mu_{n-1}=0}\equiv 0,$ by the inductive hypothesis there exist 
\begin{align*}
& f_k(s;\mu_1,\ldots,\mu_{n-1}),\ldots,f_{n-1}(s;\mu_1,\ldots,\mu_{n-1})\in\F_L^{K-1}(0_{n-1})\\
\intertext{such that}
&f(s;\mu_1,\ldots,\mu_{n-1},0)=\sum_{i=k}^{n-1}\mu_if_i(s;\mu_1,\ldots,\mu_{n-1}).
\end{align*} 
Due to $\F_L^{K-1}(0_{n-1})\subset\F_L^{K-1}(0_n),$ see \defic{defi2}, the combination of this identity with \refc{l1eq1} shows that $f=\sum_{i=k}^n\mu_i f_i$ with $f_{k},\ldots,f_n\in\F_L^{K-1}(0_n)$ as desired. This shows the inductive step and concludes the proof of the result. 
\end{prova}

We state next our first result about the bifurcation of critical periodic orbits from the outer boundary of an isochronous center. With regard to its proof let us advance that, after a convenient division in the space of coefficients, we proceed as in the proofs of Bautin \cite[\S 3]{Bautin} and Chicone and Jacobs \cite[Theorem~2.2]{Chicone} for the analogous results about limit cycles and critical periods, respectively, bifurcating from the center. Here we tackle the bifurcation from the polycycle, which is more challenging because, contrary to the center, the period function cannot be analytically extended there. To overcome this difficulty it is crucial the fact that the flatness of the remainder in the asymptotic expansion is preserved after the derivation with respect to the parameters. 

\begin{prop}\label{I2}
If $\np_0=(-\frac{1}{2},2)$ then  $\mathrm{Crit}\bigl((\out_{\np_0},X_{\np_0}),X_{\np}\bigr)= 1$.
\end{prop}

\begin{prova}
We show first the upper bound $\mathrm{Crit}\bigl((\out_{\np_0},X_{\np_0}),X_{\np}\bigr)\leqslant 1$, which constitutes the difficult part of the proof. To this end, following the notation introduced in \secc{subsec_2.2}, we define $P(s;\np)$ to be the period of the periodic orbit of $X_\np$ passing through the point $(p_1-s,0)$. Thanks to the reversibility of~$X_\np$ with respect to $\{y=0\}$ it turns out that $P(s;\np)=2T(s;\np)$ where $T(\,\cdot\,;\np)$ is the Dulac time that we consider in \propc{prop52}. Thus, by applying $(d)$ in that result
and setting $\lambda(\np)=\frac{1}{2(F-1)}$, we can assert that
\[
 T(s;\np)=T_{00}(\np)+T_{01}(\np)s^\lambda+T^{\frac{1}{2}}_{101}(\np)s\omega_{1-2\lambda}(s)
        +T_{100}^{\frac{1}{2}}(\np)s+r_1(s;\np),
\]
where $r_1\in\F_{3/2-\upsilon}^\infty(\np_0)$ for all $\upsilon>0$ small enough, the coefficients are analytic in a neighbourhood of $\np_0=(-\frac{1}{2},2)$ and, moreover, the gradients of $T_{01}$ and $T_{101}^{\frac{1}{2}}$ are linearly independent at $\np_0$. Since one can verify that $\partial_ss\omega_\alpha(s)=(1-\alpha)\omega_\alpha(s)-1$, the derivation of the above equality yields 
\[
 s^{1-\lambda}T'(s;\np)=\lambda T_{01}(\np)+2\lambda T_{101}^{\frac{1}{2}}(\np)s^{1-\lambda}\omega_{1-2\lambda}(s)+\big(T_{100}^{\frac{1}{2}}(\np)-T_{101}^{\frac{1}{2}}(\np)\big)s^{1-\lambda}+r_2(s;\np)
\]
where, by using Lemmas~A.3 and~A.4 in \cite{MV20a}, the remainder $r_2\!:=s^{1-\lambda}\partial_sr_1$ belongs to $\F_{1-\upsilon}^\infty(\np_0)$ because $\partial_sr_1\in\F_{1/2-\upsilon}^\infty(\np_0)$ and, on the other hand, $s^{1-\lambda}\in\F_{1/2-\upsilon}^\infty(\np_0)$ due to $\lambda(\np_0)=1/2.$ Note furthermore that $\hat\np=\Psi(\np)\!:=\big(\lambda(\np) T_{01}(\np),2\lambda(\np) T_{101}^{\frac{1}{2}}(\np)\big)$ is local analytic change of coordinates at $\np_0=(-\frac{1}{2},2)$ such that $\Psi(\np_0)=(0,0)$. We can thus write
\begin{equation}\label{p32eq1}
 \mathscr R_1(s;\hat\np)\!:=\left.s^{1-\lambda(\np)}T'(s;\np)\right|_{\np=\Psi^{-1}(\hat\np)}=\hat\np_1+\hat\np_2s^{1-\hat\lambda}\omega_{1-2\hat\lambda}(s)+a(\hat\np)s^{1-\hat\lambda}+h(s;\hat\np),
\end{equation}
where we set $\hat\lambda(\hat\np)\!:=\lambda(\Psi^{-1}(\hat\np))$ for shortness and define
\[
 a(\hat\np)\!:=\big(T_{100}^{\frac{1}{2}}-T_{101}^{\frac{1}{2}}\big)(\Psi^{-1}(\hat\np))
 \text{ and }
 h(s;\hat\np)\!:=r_2(s;\Psi^{-1}(\hat\np))\in\F^\infty_{1-\upsilon}(0_2).
\]

Recall at this point that the center at the origin of $X_{\np_0}$ is isochronous, so that $T'(s;\np_0)\equiv 0.$ Consequently, due to $\Psi(\np_0)=(0,0)$, 
\[
a(0,0)=0\text{  and }h(s;0,0)\equiv 0.
\]
By the Weierstrass Division Theorem (see for instance \cite[Theorem 1.8]{Greuel}), the first equality implies that $a(\hat\np)=\hat\np_1a_1(\hat\np)+\hat\np_2a_2(\hat\np)$ with~$a_1$ and~$a_2$ analytic functions at $(0,0).$ On the other hand, by \lemc{Fdiv}, $h(s;\hat\np)=\hat\np_1h_1(s;\hat\np)+\hat\np_2h_2(s;\hat\np)$ with $h_i\in\F^\infty_{1-\nu}(0_2).$ Therefore, from \refc{p32eq1},
\[
 \mathscr R_1(s;\hat\np)
 =\hat\np_1\Big(1+a_1(\hat\np)s^{1-\hat\lambda}+h_1(s;\hat\np)\Big)
 +\hat\np_2\Big(s^{1-\hat\lambda}\omega_{1-2\hat\lambda}(s)+a_2(\hat\np)s^{1-\hat\lambda}+h_2(s;\hat\np)\Big)
\]
Since $h_i\in\F^\infty_{1-\upsilon}(0_2),$ $h_i(s;\hat\np)$ tends to zero uniformly for $\hat\np\approx (0,0)$ as $s\to 0^+$ (see \defic{defi2}). Due to $\lambda(\np_0)=1/2,$ this is also the case of $s^{1-\hat\lambda}$ and $s^{1-\hat\lambda}\omega_{1-2\hat\lambda}(s)$ by $(c)$ of Lemma~A.4 in \cite{MV20a}. Hence there exists a neighbourhood $U$ of $(0,0)$ such that $\lim_{s\to 0^+}(1+a_1(\hat\np)s^{1-\hat\lambda}+h_1(s;\hat\np))=1$ uniformly on $U$. Accordingly, the function
\begin{align}\label{p32eq2}
 \mathscr R_2(s;\hat\np)\!:=&\,\frac{\mathscr R_1(s;\hat\np)}{1+a_1(\hat\np)s^{1-\hat\lambda}+h_1(s;\hat\np)}
 =\hat\np_1+\hat\np_2\ell(s;\hat\np)\\
\intertext{ with }\notag
\ell(s;\np)\!:=&\,\frac{s^{1-\hat\lambda}\omega_{1-2\hat\lambda}(s)+a_2(\hat\np)s^{1-\hat\lambda}+h_2(s;\hat\np)}{1+a_1(\hat\np)s^{1-\hat\lambda}+h_1(s;\hat\np)}
\end{align}
belongs to the class $\mathscr C^\infty_{s>0}(U)$, see \defic{defi_fun}. 

We claim that, by shrinking $U$, there exists $s_0>0$ such that $\mathscr R_2(s;\hat\np)$ has at most one zero on $(0,s_0)$, counted with multiplicities, for all $\hat\np=(\hat\np_1,\hat\np_2)\in U\setminus\{(0,0)\}.$ Indeed, to show this note first that  if $\hat\np_2=0$ then $\mathscr R_2(s;\hat\np)=\hat\np_1\neq 0,$ so that there is nothing to be proved in this case. 
Let us study consequently the case $\hat\np_2\neq 0$. To this end we observe that $\mathscr R_2'(s;\hat\np)=\hat\np_2\ell'(s;\hat\np)$ where, using a more compact notation,
\begin{align*}
 \ell'(s;\hat\np)
 &
 =\partial_s\left(
 \frac{s^{1-\hat\lambda}\omega_{1-2\hat\lambda} +a_2s^{1-\hat\lambda}+\F^\infty_{1-\upsilon}}{1+a_1 s^{1-\hat\lambda}+\F^\infty_{1-\upsilon}}
 \right)
 \\[2pt]
 &\hspace{-0.5truecm}=\frac{\omega_{1-2\hat\lambda}}{s^{\hat\lambda}}
 \left(
 \frac{\hat\lambda+\frac{a_2-a_2\hat\lambda-1}{\omega_{1-2\hat\lambda}}+\F^\infty_{\frac{1}{2}-\upsilon'}}{1+a_1s^{1-\hat\lambda}+\F^\infty_{1-\upsilon}}
 -\frac{\big(s^{1-\hat\lambda}+a_2\frac{s^{1-\hat\lambda}}{\omega_{1-2\hat\lambda}}+\F^\infty_{1-\upsilon'}\big)\big((1-\hat\lambda)a_1+\F^\infty_{\frac{1}{2}-\upsilon'}\big)}{(1+a_1s^{1-\hat\lambda}+\F^\infty_{1-\upsilon})^2}
 \right).
 \end{align*}
Here we use the identity $\partial_ss^b\omega_\alpha(s)=s^{b-1}((b-\alpha)\omega_\alpha(s)-1)$ and that, by Lemmas A.3 and A.4 in~\cite{MV20a}, we have $1/\omega_{1-2\hat\lambda}\in\F^\infty_{-\upsilon}(0_2)$ and $s^{-\hat\lambda}\in\F_{-1/2-\upsilon}^\infty(0_2)$ for all $\upsilon>0$ small enough due to $\hat\lambda(0,0)=\frac{1}{2}$ and, moreover, that the inclusion $\F_{L}^\infty\F_{L'}^\infty\subset \F_{L+L'}^\infty$ holds. We also remark that, by $(a)$ of Lemma~A.4 in \cite{MV20a},
\[
 \lim_{s\to 0^+}\frac{1}{\omega_{1-2\hat\lambda}(s)}=\frac{|1-2\hat\lambda|-(1-2\hat\lambda)}{2}\text{ uniformly for $\hat\np\approx (0,0)$.}
\]
On account of this, from the above expression of $\ell'$ we obtain that 
\[
\lim_{s\to 0^+}\frac{s^{\hat\lambda}\ell'(s;\hat\np)}{\omega_{1-2\hat\lambda}(s)}
=b(\hat\np)
\text{ uniformly for $\hat\np\approx (0,0)$,}
\]
where $b(\hat\np)\!:=\hat\lambda+\frac{1}{2}(a_2-a_2\hat\lambda-1)(|1-2\hat\lambda|-(1-2\hat\lambda))$. Since $\hat\lambda(0,0)=1/2$, it is clear that $b(\hat\np)$ is a non-vanishing continuous function in a neighbourhood of $(0,0)$. Accordingly, due to $\mathscr R'_2(s;\hat\np)=\hat\np_2\ell'(s;\hat\np)$, we can assert that
\[
 \lim_{s\to 0^+}\frac{s^{\hat\lambda}\mathscr R'_2(s;\hat\np)}{\omega_{1-2\hat\lambda}(s)}=\hat\np_2b(\hat\np) \text{ uniformly for $\hat\np\approx (0,0)$.}
\] 
Since $\omega_\alpha(s)$ only vanishes at $s=1,$ by shrinking $U$ if necessary, we can assert the existence of some $s_0\in (0,1)$ such that $\mathscr R_2'(s;\hat\np)\neq 0$ for all $s\in (0,s_0)$ and $\hat\np=(\hat\np_1,\hat\np_2)\in U$ with $\hat\np_2\neq 0$. Therefore, by Rolle's Theorem, $\mathscr R_2(s;\hat\np)$ has at most one zero on $(0,s_0)$ counted with multiplicities. This shows the validity of the claim for the case $\np_2\neq 0.$ 

Recall finally that the period function $P(s;\np)$ is twice the Dulac time $T(s;\np)$. 
Thus, taking the claim into account, from \refc{p32eq1} and \refc{p32eq2} it turns out that $V\!:=\Psi^{-1}(U)$ is an open neighbourhood of $\np_0=(-\frac{1}{2},2)$ verifying that $P'(s;\np)$ has at most one isolated zero on $(0,s_0)$, counted with multiplicities, for all $\np\in V.$ (To be more precise, the claim applies for the punctured neighbourhood $V\setminus\{\np_0\}$ and, on the other hand, $P'(s;\np_0)\equiv 0,$ so that it has not any isolated zero.) Hence, see \defic{Z}, $\mathcal Z_0(P'(\,\cdot\,;\np),\np_0)\leqslant 1$. Therefore the upper bound $\mathrm{Crit}\bigl((\out_{\np_0},X_{\np_0}),X_{\np}\bigr)\leqslant 1$ follows from assertion $(2a)$ in \lemc{CZ}
since, using the notation in that result, $P(s;\nu)=\hat P(\sigma(s;\nu);\nu)$ with $\sigma(s;\nu)=(p_1-s,0)$ for $s\in [0,\delta)$. Thus it only remains to show that this upper bound is achieved. To this end we recall that, by \cite[Theorem A]{MMV2}, $\np_0=(-\frac{1}{2},2)$ is a local bifurcation value of the period function at the outer boundary, see \defic{def2}. Then, since the period annulus of the centers under consideration varies continuously, see \obsc{+CZ}, by applying $(a)$ in \lemc{BZ} we get that $\mathrm{Crit}\bigl((\out_{\np_1},X_{\nu_1}),X_{\np}\bigr)\geqslant 1$. This completes the proof of the result.
\end{prova}

The following is our second result about the criticality of the quadratic isochrones. 

\begin{prop}\label{I4}
If $\np_0=(-\frac{1}{2},\frac{1}{2})$ then $\mathrm{Crit}\bigl((\out_{\np_0},X_{\np_0}),X_{\np}\bigr)= 1$.
\end{prop}

\begin{prova}
We prove $\mathrm{Crit}\bigl((\out_{\np_0},X_{\np_0}),X_{\np}\bigr)\leqslant 1$ first, which is the most complicated part of the proof. To this end, for each $s\in (0,1)$ we denote by $P(s;\np)$ the period of the periodic orbit of $X_\np$ passing thought the point $(1-s,0)\in\R^2,$ see \figc{fig1}. Then, on account of the reversibility of the vector field with respect to $\{y=0\}$, it follows that $P(s;\np)=2T(s;\np)$, where $T(\,\cdot\,;\np)$ is the Dulac time introduced before \propc{coeficients_Loud1}. Thanks to that result we have thus the asymptotic expansion of $P(s;\np)$ near the polycycle, which corresponds to $s=0.$ On the other hand, it is well known that the period function can be analytically extended to the center (which corresponds to $s=1$ with this parametrization) because it is non-degenerated. The coefficients of the Taylor's series of $P'(s;\np)$ at $s=1$ belong to the polynomial ring $\R[D,F]$. Chicone and Jacobs show (see Lemma 3.1 and Theorem 3.9 in \cite{Chicone}) that these coefficients are in the ideal generated by 
\begin{align*}
 p_2(D,F)=&\,10D^2+10 D F-D+4F^2-5F+1\\
 \intertext{and}
 p_4(D,F)=&\,1540 D^4+4040 D^3 F+1180 D^3+4692 D^2
           F^2+1992 D^2 F+ 453 D^2\\
         &+2768 D F^3+228 D F^2+318 D F-2 D+784 F^4-616 F^3-63
           F^2-154 F+49
\end{align*}
over the local ring $\R\{D,F\}_{\np_i}$ of convergent power series at $\np_i$ localized at any of the of the four quadratic isochrones $\np_0\!:=(-\frac{1}{2},\frac{1}{2}),$ $\np_1\!:=(0,1)$ $\np_2\!:=(0,\frac{1}{4})$ and $\np_3\!=(-\frac{1}{2},2).$ With regard to the first one, we claim that the ideal $\mf B\!:=(p_2,p_4)$ is equal to $\left(D+F,(2F-1)^2\right)$ over the local ring $\R\{D,F\}_{\np_0}$. Indeed, to prove this we use that
\begin{equation}\label{p35eq1}
 \left(\begin{array}{c}
  p_2 \\ p_4
 \end{array}\right)
 =
 \left(\begin{array}{cc}
  q_{11} & q_{12} \\ q_{21} & q_{22}
 \end{array}\right)
 \left(\begin{array}{c}
  (2F-1)^2 \\ D+F
 \end{array}\right)
\end{equation}
with $q_{11}=1,$ $q_{12}=10D-1$, $q_{21}=52{F}^{2}+44F+49$ and 
\[
q_{22}=576{F}^{3}+ ( 2192D-584 ) {F}^{2}+ ( 2500{D}^{2}
+812D-135 ) F+1540{D}^{3}+1180{D}^{2}+453D-2.
\]
(The idea to obtain this is that the zero of $p_{2i}|_{D=-F}$ at $F=1/2$ has multiplicity two.) 
From~\refc{p35eq1} we get that $p_{2i}\in\left(D+F,(2F-1)^2\right)$ over the polynomial ring $\R[D,F].$ Conversely, since one can verify that the determinant $q_{11}q_{22}-q_{21}q_{12}$ is different from zero at $\np_0=(-\frac{1}{2},\frac{1}{2}),$ by inverting the matrix in~\refc{p35eq1} it follows that $(2F-1)^2\in\mf B$ and $D+F\in\mf B$ over the local ring $\R\{D,F\}_{\np_0}$. This proves the validity of the claim. Consequently, thanks to the result of Chicone and Jacobs mentioned above, we have the following equality between ideals over the local ring $\R\{D,F\}_{\np_0}:$
\[
 \mf B=\big(D+F,(2F-1)^2\big)=\big(P^{(i)}(1;\np),i\in\N\big).
\]
Now the crucial point is that the ideal $\big(P^{(i)}(s_0;\np),i\in\N\big)$
does not depend on the point $s_0\in (0,1].$ Indeed, this follows verbatim the argument that R. Roussarie gives in \cite[pp. 76--78]{Roussarie89} or \cite[\S 4.3.1]{Roussarie} to justify the same property about the ideal of the displacement map, the so-called Bautin ideal.  
Here we also use that, such as the displacement map, the period function $P(s;\np)$ extends analytically to the non-degenerate center (i.e., $s=1$). Accordingly,
\begin{equation}\label{p35eq4}
 \mf B=\big(D+F,(2F-1)^2\big)=\big(P^{(i)}(s_0;\np),i\in\N\big)\text{ for all $s_0\in (0,1].$}
\end{equation}

We turn next to the study of the period function near the polycycle (i.e., $s=0).$ In this regard by applying $(d)$ in \propc{coeficients_Loud1} we can assert that, for all $\upsilon>0$ small enough, 
\[
 P(s;\np)=2T_{00}(\np)+2T_{101}^1(\np)s\omega_{1-\lambda}(s)+2T_{100}^1(\np)s+\F_{2-\upsilon}^\infty(\np_0),
 \]
where $\lambda(\np)=\frac{F}{1-F}$ and
\begin{equation}\label{p35eq2}
T_{101}^1(\np)=-\rho_4(\np)(F-1/2)^2\text{ and }T_{100}^1(\np)=\rho_5(\np)(D+1/2)+\rho_6(\np)(F-1/2)
\end{equation}
for some analytic positive functions $\rho_4,$ $\rho_5$ and $\rho_6$ in a neighbourhood of $\np_0=(-\frac{1}{2},\frac{1}{2})$. Consequently, on account of the identity $\partial_ss\omega_\alpha(s)=(1-\alpha)\omega_\alpha(s)-1$ and assertion $(f)$ of Lemma A.3 in~\cite{MV20a},
\[
 P'(s;\np)=2\lambda T_{101}^1(\np)\omega_{1-\lambda}(s)+2\big(T_{100}^1-T_{101}^1\big)(\np)+\F_{1-\upsilon}^\infty(\np_0).
\]
Furthermore, from \refc{p35eq2} it follows that 
\[
 \hat\np=\Psi(\np)\!:=\left(
  (F-1/2)\sqrt{2\lambda\rho_4(\np)},2\big(T_{100}^1-T_{101}^1\big)(\np)
 \right)
\]
is an analytic local change of coordinates in a neighbourhood of $\np=\np_0$ because one can verify that its Jacobian at $\np_0=(-\frac{1}{2},\frac{1}{2})$ is equal to $-2\rho_5(\np_0)\sqrt{2\rho_4(\np_0)}\neq 0.$ Setting $\hat\np=(\hat\np_1,\hat\np_2)$, observe that then
\begin{equation}\label{p35eq5}
 P'(s;\Psi^{-1}(\hat\np))=-\hat\np_1^2\omega_{1-\hat\lambda}(s)+\hat\np_2+f(s;\hat\np),
\end{equation} 
where $f\in\F_{1-\upsilon}^\infty(0_2)$ and we denote $\hat\lambda\!:=\lambda(\Psi^{-1}(\hat\np))$ for shortness. 

We claim that $\mf B=(\hat\np_1^2,\hat\np_2)$ over the local ring $\R\{D,F\}_{\np_0}$. To show this we note that
\[
 \hat\np_2\big|_{D=-F}=(2F-1)\big(\rho_6-\rho_5\big)(-F,F)+2(F-1/2)^2\rho_4(-F,F).
\]
Since $\rho_5(\np_0)=\rho_6(\np_0)$ by $(d)$ in \propc{coeficients_Loud1}, it follows that 
$\big(\rho_6-\rho_5\big)(-F,F)=(F-1/2)r_1(F)$ for some analytic function  $r_1$ at $F=1/2.$ Consequently $\hat\np_2|_{D=-F}=(F-1/2)^2r_2(F)$ with $r_2$ analytic at $F=1/2.$ Taking this into account, the Weierstrass Division Theorem (see \cite[Theorem 1.8]{Greuel}) shows that
\[
 \hat\np_2=(D+F)q(\hat\np)+(F-1/2)^2r_2(F)
\]
for some analytic function $q$  at $\hat\np=(0,0)$ which, from \refc{p35eq2}, verifies $q(0,0)=2\rho_5(\np_0)\neq 0.$ Hence we can write 
\begin{equation*}
 \left(\begin{array}{c}
  \hat\np_1^2 \\ \hat\np_2
 \end{array}\right)
 =
 \left(\begin{array}{cc}
  0 & 2\lambda\rho_4(\np) \\  q(\Psi(\np)) & r_2(F) 
 \end{array}\right)
 \left(\begin{array}{c}
  D+F \\ (F-1/2)^2
 \end{array}\right),
\end{equation*}
where the matrix has an analytic inverse at $\np=\np_0.$ Taking \refc{p35eq4} into account this shows that $\mf B=(\hat\np_1^2,\hat\np_2)$ over the local ring $\R\{D,F\}_{\np_0}$, as desired.

Recall at this point that the center of $X_{\np_0}$ is isochronous. Hence $P'(s;\np_0)\equiv 0.$ Thus, taking $\Psi(\np_0)=(0,0)$ into account, from \refc{p35eq5} we get that $f(s;0,0)\equiv 0.$ Having this in mind we write the remainder in \refc{p35eq5} as 
 \[
  f(s;\hat\np)=f_1(s;\hat\np)+f_2(s;\hat\np_1)
 \]   
with $f_1(s;\hat\np)\!:=f(s;\hat\np_1,\hat\np_2)-f(s;\hat\np_1,0)$ and $f_2(s;\hat\np_1)\!:=f(s;\hat\np_1,0)$. Since $f_1(s;\hat\np_1,0)\equiv 0$, the application of \lemc{Fdiv} shows the existence of $g_1\in\F_{1-\upsilon}^\infty(0_2)$ such that $f_1(s;\hat\np)=\hat\np_2 g_1(s;\hat\np)$. Due to $f_2(s;0)\equiv 0$ and again by \lemc{Fdiv},  $f_2(s;\hat\np_1)=\hat\np_1 g_2(s;\hat\np_1)$ with $g_2\in\F_{1-\upsilon}^\infty(0_2)$. We also have $g_2(s;0)\equiv 0$ because, otherwise, it would exist $s_0>0$ such that $g_2(s_0;\hat\np_1)\neq 0$ for all $\hat\np_1\approx 0$. In this case, taking the claim into account together with \refc{p35eq4} and \refc{p35eq5},
\[
 P'(s_0;\Psi^{-1}(\hat\np))=-\hat\np_1^2\omega_{1-\hat\lambda}(s_0)+\hat\np_2+\hat\np_2 g_1(s_0;\hat\np)+\hat\np_1 g_2(s_0;\hat\np_1)\in\mf B=(\hat\np_1^2,\hat\np_2).
\]
From here, since each $g_i(s_0;\hat\np)$ is analytic at $\hat\np=(0,0)$ and $g_2(s_0;\hat\np_1)\neq 0$ for $\hat\np_1\approx 0$, we would get that $\hat\np_1\in (\hat\np_1^2,\hat\np_2)$ over the local ring $\R\{D,F\}_{\np_0}$, which is clearly a contradiction. Concerning the analyticity of $g_i(s_0;\hat\np)$, let us remark that it follows by applying the Weierstrass Division Theorem thanks to the analyticity of $f(s_0;\hat\np)$ at $\hat\np=(0,0),$ which in its turn follows from \refc{p35eq5} noting that:
\begin{itemize}
\item $P'(s_0;\np)$ is analytic at $\np=\np_0$ because $\{X_\np\}_{\np\in\R^2}$ is an analytic family of the vector fields and hence, by \lemc{obert}, $(s,\np)\mapsto P(s;\np)=\hat P((1-s,0);\np)$ is analytic on $(0,1)\times\R^2,$ 
\item the change of coordinates $\hat\np=\Psi(\np)$ is analytic at $\np=\np_0,$ and
\item $\omega_{\alpha}(s_0)$ is analytic at $\alpha=0$ because we can write it as 
         $\omega_{\alpha}(s_0)=F(\alpha\ln s_0)\alpha$ with $F(x)=\frac{e^{-x}-1}{x}.$
\end{itemize}
Hence $g_2(s;0)\equiv 0$ and, by \lemc{Fdiv} once again, $f_2(s;\hat\np_1)=\hat\np_1^2g_3(s;\hat\np_1)$ with $g_3\in\F_{1-\upsilon}^\infty(0_2)$. Summing-up all this information about the remainder, from \refc{p35eq5} we get that
\[
P'(s;\Psi^{-1}(\hat\np))=-\hat\np_1^2\big(\omega_{1-\hat\lambda}(s)+\F_{1-\upsilon}^\infty(0_2)\big)+\hat\np_2\big(1+\F_{1-\upsilon}^\infty(0_2)\big).
\]
We are now in position to complete the proof by showing that there exist $s_0>0$ and an open neighbourhood~$U$ of $\hat\np=(0,0)$ such that 
\[
 G(s;\hat\np)\!:=\frac{P'(s;\Psi^{-1}(\hat\np))}{\omega_{1-\hat\lambda}(s)+\F_{1-\upsilon}^\infty(0_2)}=-\hat\np_1^2+\hat\np_2\frac{1+\F_{1-\upsilon}^\infty(0_2)}{\omega_{1-\hat\lambda}(s)+\F_{1-\upsilon}^\infty(0_2)}
\] 
has at most one zero on $(0,s_0)$, counted with multiplicities, for all $\hat\np=(\hat\np_1,\hat\np_2)\in U\setminus\{(0,0)\}.$ This is clear in case that $\hat\np_2=0.$ To tackle the case $\hat\np_2\neq 0$ we compute the derivative with respect to $s$ to obtain that
\begin{align*}
 G'(s;\hat\np)&=\hat\np_2\partial_s\left(
 \frac{1+\F_{1-\upsilon}^\infty}{\omega_{1-\hat\lambda}(s)+\F_{1-\upsilon}^\infty}
 \right)
 =\hat\np_2\partial_s\left(
 \frac{1+\F_{1-\upsilon}^\infty}{\omega_{1-\hat\lambda}(s)(1+\F_{1-2\upsilon}^\infty)}
 \right)\\[6pt]
 &=\hat\np_2\partial_s\left(
 \frac{1+\F_{1-2\upsilon}^\infty}{\omega_{1-\hat\lambda}(s)}
 \right)
 =\frac{\hat\np_2}{s^{2-\hat\lambda}\omega^2_{1-\hat\lambda}(s)}(1+\F_{1-2\upsilon}^\infty)
  +\frac{\hat\np_2}{\omega_{1-\hat\lambda}(s)}\F_{-2\upsilon}^\infty\\[6pt]
  &=\frac{\hat\np_2}{s^{2-\hat\lambda}\omega^2_{1-\hat\lambda}(s)}\big(1+\F_{1-2\upsilon}^\infty+s^{2-\hat\lambda}\omega_{1-\hat\lambda}(s)\F_{-2\upsilon}^\infty\big)
  =\frac{\hat\np_2}{s^{2-\hat\lambda}\omega^2_{1-\hat\lambda}(s)}\big(1+\F_{1-3\upsilon}^\infty\big).
\end{align*}
Here, in the second equality we apply first assertion $(c)$ of Lemma A.4 in~\cite{MV20a} to get that $1/\omega_{1-\hat\lambda}(s)\in\F^\infty_{-\upsilon}$ for all $\upsilon>0$ small enough, due to $\hat\lambda(0,0)=1$, and use next that $\F_{-\upsilon}^\infty\F_{1-\upsilon}^\infty\subset\F_{1-2\upsilon}^\infty$ from $(g)$ of Lemma~A.3 in~\cite{MV20a}. In the third equality,
on account of $\frac{1}{1+s}-1\in\F_1^\infty$ and by $(h)$ of Lemma~A.3 in~\cite{MV20a}, we use first the inclusion $\frac{1}{1+\F^\infty_{1-2\upsilon}}\subset 1+\F^\infty_{1-2\upsilon}$. Then, by using $(d)$ and $(g)$ of Lemma~A.3 in~\cite{MV20a}, we expand the numerator to get that $(1+\F_{1-\upsilon}^\infty)(1+\F_{1-2\upsilon}^\infty)\subset 1+\F_{1-2\upsilon}^\infty$. Next, in the fourth equality we use that $\partial_s\omega_\alpha(s)=s^{-\alpha-1}$ and assertion $(f)$ of Lemma A.3 in \cite{MV20a} to deduce that $\partial_s\F_{1-2\upsilon}^\infty\subset\F_{-2\upsilon}^\infty$. Finally in the last equality we apply $(c)$ of Lemma~A.4 in \cite{MV20a} to get that $s^{2-\hat\lambda}\omega_{1-\hat\lambda}(s)\in\F^\infty_{1-\upsilon}$ and we use again that $\F_{1-\upsilon}^\infty\F_{-2\upsilon}^\infty\subset\F_{1-3\upsilon}^\infty$.
On account of \defic{defi2} we can assert the existence of some $s_0\in (0,1)$ and a neighbourhood $U$ of~$(0,0)$ such that $G'(s;\hat\np)\neq 0$ for all $s\in (0,s_0)$ and $\hat\np\in U$ with $\hat\np_2\neq 0.$ Consequently $P'(s;\Psi^{-1}(\hat\np))$ has at most one isolated zero on $(0,s_0),$ counted with multiplicities, for all $\hat\np\in U\setminus\{(0,0)\}.$ Thus, on account of \defic{Z} and the fact that $\Psi(\np_0)=(0,0)$, we get $\mathcal Z_0(P'(\,\cdot\,;\np),\np_0)\leqslant 1$. Finally the upper bound $\mathrm{Crit}\bigl((\out_{\np_0},X_{\np_0}),X_{\np}\bigr)\leqslant 1$ follows from assertion $(2a)$ in \lemc{CZ}
since, using the notation in that result, $P(s;\nu)=\hat P(\sigma(s;\nu);\nu)$ with $\sigma(s;\nu)=(1-s,0)$ for $s\in [0,\delta)$. Therefore it only remains to show that this upper bound is attained. To this end we recall that, by \cite[Theorem A]{MMV2}, $\np_0=(-\frac{1}{2},\frac{1}{2})$ is a local bifurcation value of the period function at the outer boundary, see \defic{def2}. Then, since the period annulus of the centers under consideration varies continuously, see \obsc{+CZ}, by applying $(a)$ in \lemc{BZ} we get that $\mathrm{Crit}\bigl((\out_{\np_1},X_{\nu_1}),X_{\np}\bigr)\geqslant 1$. This finishes the proof of the result.
\end{prova}

As we explain at the beginning of this section, the maximum criticality of the period function at the inner boundary is~2 and it is achieved at the three Loud points $\nu=L_i$, see \refc{Li}. We refer the interested reader to the paper of Chicone and Jacobs~\cite{Chicone} for a proof of this result. In a joint paper with P.~Marde\v si\'c, see \cite[Theorem 4.3]{MMV2}, we prove that at each $\nu=L_i$ there exists a germ of analytic curve that consists of local bifurcation values of the period function at the interior, see \defic{def2}. Since the period function extends analytically to the center, this follows readily by applying the Weierstrass Preparation Theorem.
In our next result we identify a parameter $\nu=\nu_\star$ for which the criticality at the outer boundary is 2. Furthermore we prove that at $\nu=\nu_\star$ there exists a $\cc^1$ germ of curve of local bifurcation values of the period function at the interior. Hence, roughly speaking, this parameter is the mirror image at the outer boundary of one of the Loud points, see \figc{figdoble} and \obsc{figdobe_com}. In the statement, following the notation introduced at the beginning of \secc{subsec_2.2}, for each $s\in (0,p_1)$ and $\np\approx\np_\star$ we denote by $P(s;\np)$ the period of the periodic orbit of $X_\np$ passing through the point $(p_1-s,0)\in\R^2$. We also remark that $T_{10}$ and $T_{01}$ are the coefficients given in \propc{prop52}, which vanish at 
$\np_\star=(\mathcal G(4/3),4/3)$.

\begin{prop}\label{doble}
Let us consider $\np_\star=(\mathcal G(4/3),4/3)$. Then the following holds:
\begin{enumerate}[$(a)$]
\item $\mathrm{Crit}\bigl((\out_{\np_\star},X_{\np_\star}),X_{\np}\bigr)=2$.
\item There exist an open neighbourhood $U$ of $\np_\star$ and $s_0>0$ such that the set
\[
 \Delta\!:=\{\np\in U;\text{ there exists $s\in(0,s_0)$ such that } P'(s;\np)=P''(s;\np)=0\}
 \]
satisfies the following conditions:
\begin{enumerate}[$(b1)$]
\item Each $\np\in\Delta$ is a local bifurcation value of the period function at the interior,

\item there exist $\varepsilon>0$ and a 
 $\mathscr C^1$ injective curve $\delta:(-\varepsilon,\varepsilon)\to U$ with $\delta(0)=\np_\star$, $\delta((0,\varepsilon))=\Delta$ and such that $\delta'(0)\neq(0,0)$ is tangent to $\{\np\in U;\;T_{10}(\np)=0\}$, 
 
\item for each $\np\in\Delta$ there exists a unique $s_\np\in (0,s_0)$ such that $P'(s_\np;\np)=P''(s_\np;\np)=0$ and, moreover, $\lim_{\np\to\np_\star}s_\np=0^+$,

\item $\Delta\subset\{\np\in U; \text{ $T_{10}(\np)<0$ and $T_{01}(\np)>0$}\}$, and
\item for any $\np_0\in\Delta$ and any neighbourhood $V$ of $\np_0$ there exist $\bar\np\in V$ and different $s_1,s_2\in (0,s_0)$ such that $P'(s_1;\bar\np)=P'(s_2;\bar\np)=0$.

\end{enumerate}
\end{enumerate}
\end{prop}

\begin{prova}
We observe first of all that $\sigma(s;\nu)\!:=(p_1-s,0)$ is a parametrization of the outer boundary of the period annulus verifying the hypothesis in \lemc{CZ}. This will enable us to relate $\mathrm{Crit}\bigl((\out_{\np_\star},X_{\np_\star}),X_{\np}\bigr)$ with $\mathcal Z_0(P'(\,\cdot\,;\np),\np_\star).$ That being said, thanks to the reversibility of the vector field with respect to $\{y = 0\}$, we note that $P(s;\np) = 2T(s;\np)$, where $T(\,\cdot\,;\np)$ is the Dulac time considered in \propc{prop52}. From point~$(b)$ in that result we can assert that, for all $\upsilon>0$ small enough,
\begin{equation*}
P(s;\np)=2T_{00}(\np)+2T_{10}(\np)s+2T_{01}(\np)s^{\lambda}+2T_{20}(\np)s^2+\F_{5/2-\upsilon}^\infty(\np_\star),
\end{equation*}
where $\lambda(\np)=\frac{1}{2(F-1)},$ $T_{10}(\np_\star)=T_{01}(\np_\star)=0$, $T_{20}(\np_\star)< 0$ and the gradients $\nabla T_{10}(\np_\star)$ and $\nabla T_{01}(\np_\star)$ are linearly independent. Due to $T_{20}(\np_\star)\neq 0$, by applying  \cite[Theorem C]{MV20b} we get that $\mathcal Z_0(P'(\,\cdot\,;\np),\np_\star)\leqslant 2$.
(For readers convenience, let us explain that \cite[Theorem C]{MV20b} is a general result addressed to the Dulac time which, 
by using the well-known derivation-division algorithm, gives a bound for $\mathcal Z_0(T'(\,\cdot\,;\nu),\nu_0)$ in terms of the position of the first non-vanishing coefficient in the asymptotic expansion of $T(s;\nu)$ at $s=0$.) Consequently, by assertion $(2a)$ in \lemc{CZ}, $\mathrm{Crit}\bigl((\out_{\np_\star},X_{\np_\star}),X_{\np}\bigr)\leqslant 2$. In addition, since 
\begin{equation}\label{F1}
F_1(s;\np)\!:=P'(s;\np)=2T_{10}(\np)+2\lambda T_{01}(\np)s^{\lambda-1}+4T_{20}(\np)s+\F_{3/2-\upsilon}^\infty(\np_\star)
\end{equation}
and the gradients $\nabla T_{10}(\np_\star)$ and $\nabla T_{01}(\np_\star)$ are linearly independent, by \cite[Proposition~4.2]{MV20b} it turns out that $\mathcal Z_0(P'(\,\cdot\,;\np),\np_\star)\geqslant 2$. As a matter of fact, from the proof of that result, this lower bound is achieved by means of two different sequences of zeros of $P'(\,\cdot\,;\np)$ and, therefore, by assertion $(2b)$ in \lemc{CZ}, $\mathrm{Crit}\bigl((\out_{\np_\star},X_{\np_\star}),X_{\np}\bigr)\geqslant 2$. Accordingly $\mathrm{Crit}\bigl((\out_{\np_\star},X_{\np_\star}),X_{\np}\bigr)= 2$
and this proves $(a).$ 

Let us turn next to the proof of the  assertions in $(b).$ For this purpose, from \refc{F1}
and by applying Lemmas A.3 and A.4 in \cite{MV20a}, we get
\begin{equation}\label{F2}
F_2(s;\np)\!:=s^{2-\lambda}P''(s;\np)=2\lambda(\lambda-1)T_{01}(\np)+4T_{20}(\np)s^{2-\lambda}+\F_{1-\upsilon}^\infty(\np_\star).
\end{equation}
Setting $U_\varepsilon\!:=\{\np\in\R^2:\|\np-\np_\star\|<\varepsilon\}$, the map
$F\!:=(F_1,F_2)$ is well-defined for $(s,\np)\in (0,\varepsilon)\times U_\varepsilon$ taking $\varepsilon>0$ small enough. Since $T_{10}(\np_\star)=T_{01}(\np_\star)=0$, $T_{20}(\np_\star)\neq 0$ and the gradients $\nabla T_{10}(\np_\star)$ and $\nabla T_{01}(\np_\star)$ are linearly independent, we can assume by reducing $\varepsilon>0$ if necessary that $\hat\np=\Psi(\np)$, defined by means of
\begin{equation}\label{Psi}
 \Psi(\np)\!:=\left(\frac{T_{10}(\np)}{2T_{20}(\np)},\frac{\lambda(\np) T_{01}(\np)}{2T_{20}(\np)}\right),
\end{equation}
is an analytic change of coordinates from $U_\varepsilon$ to the neighbourhood $\hat U_{\hat\varepsilon}\!:=(-\hat\varepsilon,\hat\varepsilon)^2$ of $(0,0)=\Psi(\np_\star).$ 
Recall that our aim is to study the solutions of the system of equations $\{P'=0,P''=0\}$ which, on account of~\refc{F1} and \refc{F2}, is equivalent to $\{F_1=0,F_2=0\}.$ In order to study the latter we first lift $\Psi$ to an analytic change of 
variables ${\Phi}$ given by 
\[
 (\hat s,\hat\np)=\Phi(s,\np)\!:=(s^{2-\lambda(\np)},\Psi(\np)),
\]
which (after diminishing $\varepsilon$ and $\hat\varepsilon$ if necessary) is defined from $\mathscr U_\varepsilon\!:=(0,\varepsilon)\times U_{\varepsilon}$ to $\hat{\mathscr U}_{\hat\varepsilon}\!:=(0,\hat\varepsilon)\times \hat U_{\hat\varepsilon},$ and then we consider the map $\hat F:\hat{\mathscr U}_{\hat\varepsilon}\to\R^2$ defined by $\hat F(\hat s,\hat\np)=(\hat F_1(\hat s,\hat\np),\hat F_2(\hat s,\hat\np))$ with
\[
 \hat F_1(\hat s,\hat\np)\!:=\frac{F_1(\Phi^{-1}(\hat s,\hat\np))}{4T_{20}(\Psi^{-1}(\hat\np))}
 \text{ and }
 \hat F_2(\hat s,\hat\np)\!:=\frac{F_2(\Phi^{-1}(\hat s,\hat\np))}{4T_{20}(\Psi^{-1}(\hat\np))}.
\]
By assertions $(h)$ and $(c)$ of Lemmas A.3 and A.4 in \cite{MV20a}, respectively, it follows that
\[
 \hat F_1(\hat s,\hat\np)=\hat\np_1+\hat\np_2\hat s^{\frac{\hat\lambda-1}{2-\hat\lambda}}+
f_1(\hat s;\hat\np)
\text{ and }
\hat F_2(\hat s,\hat\np)=(\hat\lambda-1)\hat\np_2+\hat s+f_2(\hat s;\hat\np),
\]
where $f_1,f_2\in\F^\infty_{2-\upsilon}(0_2)$ for some $\upsilon>0$ small enough and we set 
$\hat\lambda(\hat\np)\!:=\lambda (\Psi^{-1}(\hat\np))$ for shortness. Here we also use that $\F^\infty_{3-\upsilon}(0_2)\subset \F^\infty_{2-\upsilon}(0_2)$ and $s=\hat s^{1/(2-\hat\lambda)}\in \F_{2-\upsilon}^\infty(0_2)$ due to $\lambda(\np_\star)=3/2$. Observe on the other hand that, via the diffeomorphism $\Phi,$ the system $\{P'(s;\np)=0,P''(s;\np)=0\}$ on~$\mathscr U_\varepsilon$ is equivalent to the system $\{\hat F_1(\hat s,\hat\np)=0,\hat F_2(\hat s,\hat\np)=0\}$ on~$\hat{\mathscr U}_{\hat\varepsilon}$. With regard to the latter note that, by  \cite[Lemma A.1]{MV20a}, the remainders $f_1$ and $f_2$ extend to $\mathscr C^1$ functions in a neighbourhood of $(0,0,0)$ satisfying that $\nabla f_1(0,0,0)=\nabla f_2(0,0,0)=(0,0,0)$. Observe in particular that $\hat F_2(\hat s,\hat\np)$ extends to a $\cc^1$ function in a neighbourhood of $(0,0,0)$. Hence, taking $\hat\lambda(0,0)=3/2$ into account, by the Implicit Function Theorem there exists a $\cc^1$ function $h(\hat s,\hat\np_1)$ in a neighbourhood of $(0,0)$ such that, by shrinking $\hat\varepsilon>0$ if necessary,  
\[
 \text{$\hat F_2(\hat s,\hat\np)=0$ with $(\hat s,\hat\np)\in\hat{\mathscr U}_{\hat\varepsilon}$} \Leftrightarrow \hat\np_2=h(\hat s,\hat\np_1).
\]
Furthermore $h$ satisfies $h(0,\hat\np_1)\equiv 0$ and $\nabla h(0,0)=(-2,0)$. Our next task is to substitute $\hat\np_2=h(\hat s,\hat\np_1)$ in $\hat F_1(\hat s,\hat\np_1,\hat\np_2)=0$ and analyze the resulting equation. To this end we extend $\hat F_1(\hat s,\hat\np_1,\hat\np_2)|_{\hat\np_2=h(\hat s,\hat\np_1)}$ on a neighbourhood of $(\hat s,\hat\np_1)=(0,0)$ by means of
\[
(\hat s,\hat\np_1)\mapsto\hat\np_1+h(\hat s,\hat\np_1)|\hat s|^{{e}(\hat s,\hat\np_1)}+\hat f_1(\hat s,\hat\np_1),
\]
where ${e}(\hat s,\hat\np_1)\!:=\frac{\hat\lambda-1}{2-\hat\lambda}\big|_{\hat\lambda=\hat\lambda(\hat\np_1,h(\hat s,\hat\np_1))}$ and $\hat f_1(\hat s,\hat\np_1)=f_1(\hat s;\hat\np_1,h(\hat s,\hat\np_1))$ are clearly $\mathscr C^1$ in a neighbourhood of $(0,0)$. We claim that the function $g(\hat s,\hat\np_1)\!:= h(\hat s,\hat\np_1)|\hat s|^{{e}(\hat s,\hat\np_1)}$ is $\mathscr C^1$ in a neighbourhood of $(0,0)$ as well and that its gradient vanishes at $(0,0)$. To show this notice first that $g(0,\hat\np_1)=0$ and, consequently, $\partial_{\hat\np_1}g(0,\hat\np_1)=0.$ Moreover, using that 
$h(0,\hat\np_1)=0$, we get
\[
 \partial_{\hat s} g(0,\hat\np_1)=\lim\limits_{\hat s\to 0}\frac{h(\hat s,\hat\np_1)|\hat s|^{{e}(\hat s,\hat\np_1)}}{\hat s}=\lim\limits_{\hat s\to 0}\frac{h(\hat s,\hat\np_1)-h(0,\hat\np_1)}{\hat s}\lim\limits_{\hat s\to 0}|\hat s|^{{e}(\hat s,\hat\np_1)}=\partial_{\hat s}h(0,\hat\np_1)\cdot 0=0
\]
because $h$ is $\cc^1$ and ${e}(0,0)=1$ implies ${e}(0,\hat\np_1)>0$ for $\hat\np_1\approx 0.$ Similarly, if $\hat s\neq 0$ then
\begin{align*}
\partial_{\hat s}g(\hat s,\hat\np_1)&=(\partial_{\hat s}h(\hat s,\hat\np_1))|\hat s|^{{e}(\hat s,\hat\np_1)}+h(\hat s,\hat\np_1)|\hat s|^{{e}(\hat s,\hat\np_1)}\left(\log|\hat s|\partial_{\hat s}{e}(\hat s,\hat\np_1)+\frac{{e}(\hat s,\hat\np_1)}{\hat s}\right)\\
&=(\partial_{\hat s}h(\hat s,\hat\np_1))|\hat s|^{{e}(\hat s,\hat\np_1)}+\frac{h(\hat s,\hat\np_1)-h(0,\hat\np_1)}{\hat s}|\hat s|^{{e}(\hat s,\hat\np_1)}\Big(\hat s\log|\hat s|\partial_{\hat s}{e}(\hat s,\hat\np_1)+{e}(\hat s,\hat\np_1)\Big)
\intertext{and}
\partial_{\hat\np_1}g(\hat s,\hat\np_1)&=(\partial_{\hat\np_1}h(\hat s,\hat\np_1))|\hat s|^{{e}(\hat s,\hat\np_1)}+h(\hat s,\hat\np_1)|\hat s|^{{e}(\hat s,\hat\np_1)}\log|\hat s|\partial_{\hat\np_1}{e}(\hat s,\hat\np_1)
\end{align*}
tend to zero as $\hat s\to 0$ uniformly on $\hat\np_1\approx 0$. This clearly implies that $g$ is $\mathscr C^1$ in a neighbourhood of $(0,0)$ and $\nabla g(0,0)=(0,0)$, so that the claim is true. Thus, by applying the Implicit Function Theorem to the ``extended equation''
\[
\hat\np_1+h(\hat s,\hat\np_1)|\hat s|^{{e}(\hat s,\hat\np_1)}+\hat f_1(\hat s,\hat\np_1)=0
\]
and reducing $\hat\varepsilon>0$ once again, we obtain a $\cc^1$ function $\ell(\hat s)$ on $(-\hat\varepsilon,\hat\varepsilon)$ such that  $\hat F_1(\hat s,\hat\np_1,h(\hat s,\hat\np_1))=0$ with $(\hat s,\hat\np_1)\in (0,\hat\varepsilon)\times (-\hat\varepsilon,\hat\varepsilon)$ if, and only if, $\hat\np_1=\ell(\hat s)$. Moreover $\ell(0)=0$ and $\ell'(0)=0$.  Accordingly, after shrinking $\hat{\varepsilon}>0$ once again if necessary, we can assert that 
\[
 \text{$\hat F_1(\hat s,\hat\np)=\hat F_2(\hat s,\hat\np)=0$ with $(\hat s,\hat\np)\in\hat{\mathscr U}_{\hat\varepsilon}$} \Leftrightarrow \hat\np=(\hat\np_1,\hat\np_2)=\big(\ell(\hat s),h(\hat s,\ell(\hat s)\big).
\]
At this point, since we reduced the original $\hat\varepsilon>0$,
we also diminish $\varepsilon>0$ so that $\Phi(s,\np)=(s^{2-\lambda(\np)},\Psi(\np))$ is still a diffeomorphism from $\mathscr U_\varepsilon$ into $\hat{\mathscr U}_{\hat\varepsilon}$. Then, from \refc{F1} and \refc{F2}, the following assertions are equivalent:
\begin{enumerate}[(1)]
\item $\np_0\in\Delta\!:=\{\np\in U_\varepsilon; \text{ there exists $s\in (0,\varepsilon)$ such that $P'(s;\np)=P''(s;\np)=0$}\}$,
\item $\Psi(\np_0)=\big(\ell(\hat s_0),h(\hat s_0,\ell(\hat s_0))\big)$ for some $\hat s_0\in (0,\hat\varepsilon)$,
\item $\np_0=\delta(t_0)$ for some $t_0\in (0,\hat\varepsilon)$, where 
$\delta(t)\!:=\Psi^{-1}\big(\ell(t),h(t,\ell(t))\big)$. 
\end{enumerate}
It is clear from these equivalences that \map{\delta}{(-\hat\varepsilon,\hat\varepsilon)}{U_\varepsilon\subset\R^2} is a $\cc^1$ parametrized curve with $\delta(0)=\np_\star$ satisfying that $\delta\big((0,\hat\varepsilon)\big)=\Delta$. One can easily verify, taking $\ell'(0)=0$ and $\partial_1h(0,0)\neq 0$ into account, together with the definition of $\Psi$ in \refc{Psi}, that $\delta'(0)$ is a non-zero vector tangent to $\{\np\in U_\varepsilon;\;T_{10}(\np)=0\}.$ In particular, on account of $\delta'(0)\neq (0,0)$ and by reducing $\hat\varepsilon>0,$ we have that $\delta$ is one-to-one. This proves the assertion $(b2)$ in the statement.
Due to $\mathcal Z_0(P'(\,\cdot\,;\np),\np_\star)\leqslant 2$, and after shrinking $\varepsilon>0$ if necessary, note also that the zeros of $P'(\,\cdot\,;\np)$ on $(0,\varepsilon)$ can have at most multiplicity two. Therefore, since the interior of $\Delta=\delta\big((0,\hat\varepsilon)\big)$ is empty (as a subset of $\R^2$), by applying \lemc{lema-interior} we can assert that each $\np_0\in\Delta$ is a local bifurcation value of the period function at the interior, which shows the validity of~$(b1)$ in the statement. 
With regard to the assertions in $(b3)$, we note that the uniqueness of~$s_\np$ and $\lim_{\np\to\np_\star}s_\np=0^+$ follow from the point $(2)$ above using that $\hat s\mapsto h(\hat s,\ell(\hat s))$ is invertible at $\hat s=0$ and that, by definition, $\hat s=s^{2-\lambda(\np)}.$
On the other hand, since $P''(s_\np;\np)=0$ for all $\np\in\Delta$, from \refc{F2} we get that $T_{01}(\np)T_{20}(\np)<0$ for all $\np\in\Delta$. Here we also use that $\lim_{\np\to\np_\star}s_\np=0^+$ to take advantage of the properties of the remainder and the fact that $\lambda(\np_\star)=3/2.$ By arguing similarly, on account of $P'(s_\np;\np)=0$ for all $\np\in\Delta$, from~\refc{F1} it follows that $T_{01}(\np)T_{10}(\np)<0$ for all $\np\in\Delta$. Taking this into account the assertion in $(b4)$ is a consequence of $T_{20}(\np_\star)<0,$ see $(b)$ in \propc{prop52}. 
Finally, in order to prove $(b5)$, let us consider $\np\in\Delta$ and note that from \refc{F1} we obtain
\[
 \lim_{\np\to\np_\star}\left.\big(\partial_{\np_1}P'(s;\np),\partial_{\np_2}P'(s;\np)\big)\right|_{s=s_\np}=2\nabla T_{10}(\np_\star),
\] 
where we use that the flatness of the remainder $\F_{3/2-\upsilon}^\infty(\np_\star)$ is preserved after derivation with respect to parameters, see \defic{defi2}. Similarly, in this case from \refc{F2} and using also $P''(s_\np;\np)\equiv 0$, we get
\[
 \lim_{\np\to\np_\star}\left.s^{2-\lambda(\np)}\big(\partial_{\np_1}P''(s;\np),\partial_{\np_2}P''(s;\np)\big)\right|_{s=s_\np}=\frac{3}{2}\nabla T_{01}(\np_\star).
\] 
Thus, since the vectors $\nabla T_{10}(\np_\star)$ and $\nabla T_{01}(\np_\star)$ are linearly independent, so they are 
$\nabla P'(s;\np)|_{s=s_\np}$ and $\nabla P''(s;\np)|_{s=s_\np}$ for all $\np\in U_\varepsilon$ (after shrinking $\varepsilon>0$ if necessary). That being said, we fix any $\np_0\in\Delta$ and compute the second order Taylor's expansion of $P'(s;\np)$ at $s=s_{\np_0},$
\[
 P'(s;\np)=P'(s_{\np_0};\np)+P''(s_{\np_0};\np)(s-s_{\np_0})+\op(s-s_{\np_0}).
\]
Then, due to $P'(s;\np_0)\not\equiv 0,$ $P'(s_{\np_0};\np_0)=P''(s_{\np_0};\np_0)=0$ and the fact that the gradients $\nabla P'(s_{\np_0};\np)$ and $\nabla P''(s_{\np_0};\np)$ are linearly independent at $\np=\np_0,$ the application of \cite[Proposition 4.2]{MV20b} shows that for each open neighbourhood $V$ of $\np_0$ there exist $\bar\np\in V$ and two $s_1,s_2\in (0,\varepsilon)$ such that $P'(s_1;\bar\np)=P'(s_2;\bar\np)=0.$ This proves the validity of the assertion in $(b5)$ and completes the proof of the result.
\end{prova}

\usetikzlibrary{calc}

\begin{figure}[t]
\begin{center}
\hfill
\begin{tikzpicture}
\node[inner sep=0pt] at (7,0) {\includegraphics[width=14cm]{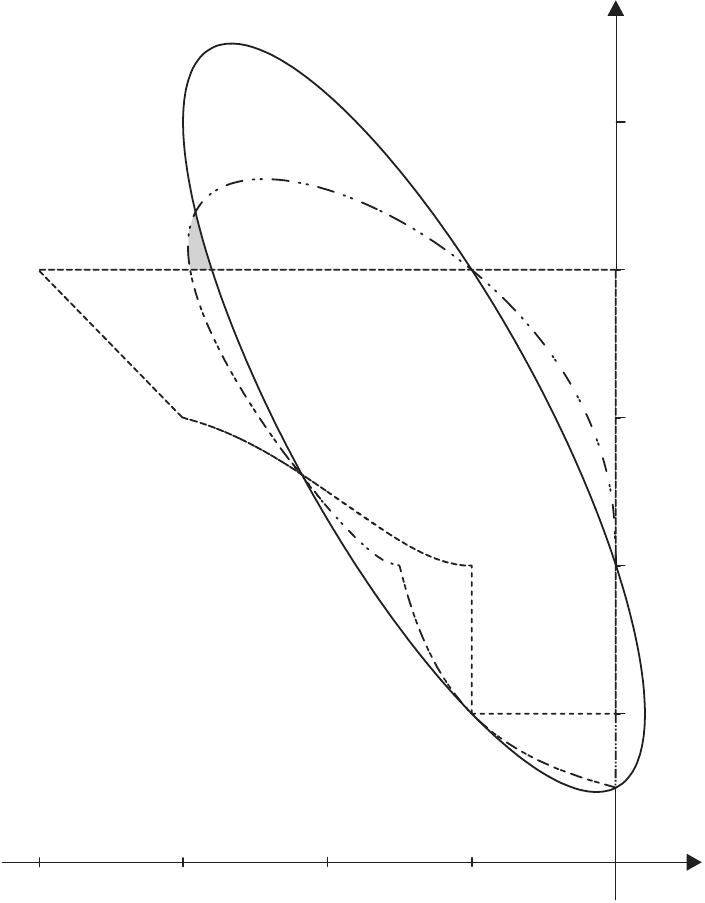}};
\node at (4.5,0.75) {$\Gamma_B$};
\node at (7.5,3) {$\Gamma_C$};
\node at (6.25,2.75) {$\Gamma_0$};
\node at (6.75,-4) {$-1$};
\node at (4.25,-4) {$-2$};
\node at (10,-4) {$D$};
\node at (9.75,4) {$F$};
\node at (9.75,-1) {$1$};
\node at (9.75,1.55) {$2$};
\begin{scope}[shift={(-4cm,0cm)}]
\draw[shift={(10cm,0cm)}] (-5.5,0) to (-1.85,0) to (-1.85,-3.5) to (-5.5,-3.5) to (-5.5,0);
\draw[shift={(10cm,0cm)}] (-5.5,-0.5) to [out=-20,in=135] (-1.85,-3);
\draw[shift={(10cm,0cm)}] (-4,-0) to [out=-80,in=120] (-2.7,-3.5);
\draw[shift={(10cm,0cm)}]  [dashed] (-4.5,0) to [out=-70,in=135] (-2.05,-3.5);
\draw[shift={(10cm,0cm)},fill] (-5,-0.7) circle [radius=0.05];
\draw[shift={(10cm,0cm)},fill] (-2.98,-3) circle [radius=0.05];
\draw[shift={(10cm,0cm)},thick] (-5,-0.7) to [out=-20,in=120] (-2.98,-3);
\draw[shift={(10cm,0cm)},fill,white] (-3.95,-1.75) circle [radius=0.25];
\draw[shift={(10cm,0cm)},fill, lightgray] (-3.925,-1.26) to (-3.68,-1.42) to (-3.51,-1.89);
\node[shift={(10cm,0cm)}]  at (-3.65,-0.25) {$\Gamma_C$};
\node[shift={(10cm,0cm)}]  at (-2.1,-2.4) {$\Gamma_B$};
\node[shift={(10cm,0cm)}]  at (-2.1,-3.15) {$\Gamma_0$};
\node[shift={(10cm,0cm)}]  at (-5,-1) {$\nu_\star$};
\node[shift={(10cm,0cm)}]  at (-4.5,-1.25) {$\Delta$};
\node[shift={(10cm,0cm)}]  at (-3.8,-2.2) {$\delta_3$};
\node[shift={(10cm,0cm)}]  at (-3.3,-3) {$L_3$};
\draw[shift={(11.35cm,-0.05cm)}] (-1.25,0.25) to (-0.5,0.25) to (-0.5,-0.47) to (-1.25,-0.47) to (-1.25,0.25);
\draw [thick] (8.15,-2) to (9,-1) to (8.8,-1.75) to (10.05,-0.4);
\draw [thick] (10.05,-0.4) to (9.95,-0.455) to (10.05,-0.4) to (10,-0.5);
\end{scope}
\end{tikzpicture}
\end{center}
\caption{Arrangement of the three types of local bifurcation curve (inner boundary, interior and outer boundary) near the parameter values $\nu=L_3$ and $\nu=\nu_\star$, see \refc{Li} and \propc{doble}, respectively. We refer the reader to \obsc{figdobe_com} for a detailed explication.}\label{figdoble}
\end{figure}

\begin{obs}\label{figdobe_com}
Let us finish this section contextualizing the results in \propc{doble}. In \figc{figdoble} we display the ellipse $\Gamma_C$ that consists of local bifurcation values of the period function at the inner boundary (i.e., the center) of $\PA$. It corresponds, see \cite[Lemma 3.1]{Chicone}, to the vanishing of the first period constant
\[
p_2(\np)=10D^2+10DF-D+4F^2-5F+1.
\]
Moreover the curve $\Gamma_B$ consists of local bifurcation values at the outer boundary (i.e., the polycycle) of $\PA$, see \cite[Theorem A]{MMV2}. It is made of the arc $\{D=\mathcal G(F): F\in (1,\frac{3}{2})\}$ joining $(-\frac{3}{2},\frac{3}{2})$ and $(-\frac{1}{2},1)$ together with several straight segments. According to \propc{doble} and \cite[Theorem 4.3]{MMV2}, respectively, the germs of curve $\Delta$ and $\delta_3$ are inside the set of local bifurcation values of the period function at the interior of $\PA$. At this moment we do not have any analytical tool to fully characterize this set. We conjecture that $\Delta$ and $\delta_3$ connect with each other to delimit a region of parameters for which the corresponding center has exactly two critical periodic orbits. With regard to this conjecture
it is proved in \cite[Theorem 5.2]{MMV2} that the center of any parameter inside one of the two light gray regions has at least two critical periodic orbits. The boundary of these regions is inside $\Gamma_C,$ $\Gamma_B$ and $\Gamma_0.$ For completeness let us explain that the curve $\Gamma_0$ consists of those parameters such that the period function tends to $2\pi$ as the periodic orbits tend to the outer boundary. 
\end{obs}

\section{Proof of Theorem~\ref{Loud}}\label{provaA}

\begin{prooftext}{Proof of Theorem~\ref{Loud}.}
The statement covers all the parameters $\nu_0\in\R^2$ outside the vertical segments 
$\ell_0\!:=\{D=-1,F\in [0,1]\}\cup\{D=0, F\in [0,\frac{1}{4})\}.$ For simplicity in the exposition, instead of proving the five assertions in the statement separately, we split $\R^2\setminus\ell_0$ depending on the result and tool applied to study the corresponding criticality. For reader's convenience we enumerate the different cases that we obtain in this way. 

\begin{enumerate}[1.]

\item Let us consider first of all the set $\ell_1\!:=\R^2\setminus (\Gamma_B\cup\Gamma_U)$, 
         where recall (see \figc{diagrama}) that $\Gamma_U$ is the 
         union of the dotted straight lines, whatever its colour is, 
         and~$\Gamma_B$ is the Jordan curve in boldface type. Then, by \cite[Theorem A]{MMV2}, 
         we know that any $\nu_0\in\ell_1$ is a local regular value of the period function at
        the outer boundary, see $(c)$ in \defic{def2}. On account of this, by~$(b)$ in \lemc{BZ} we get that
        $\mathrm{Crit}\bigl((\out_{\np_0},X_{\np_0}),X_{\np}\bigr)=0$. Here we also use that the period annuli of the 
        Loud's centers vary continuously, see \obsc{+CZ}, and that the outer boundary of  
        $\PA_\nu$ for $\nu\notin \Gamma_B\cup\Gamma_U$ is a hyperbolic polycycle, 
        see for instance \cite[\S 3.1]{MMV2}. 
                 
\item The criticality at  
         $\ell_2\!:=\{D=-\frac{1}{2},F\in(\frac{1}{2},1)\}\cup\{F=\frac{1}{2},D\in(-\frac{1}{2},0)\}$ 
         and $\ell_3\!:=\{F=\frac{1}{2}, D\in(-1,-\frac{1}{2})\}$ follows 
         from the results in \secc{subsec_2.1}. In this case $\sigma(s;\nu)=(1-s,0)$ 
         is a parametrization of the outer boundary of the 
         period annulus verifying the hypothesis 
         $(a)$, $(b)$ and $(c)$ in \lemc{CZ}. Moreover denoting by $P(s;\nu)$ the period of the periodic orbit of $X_\nu$ 
         passing through $\sigma(s;\nu)$, we have that $P(s;\nu)=2T(s;\nu)$, where $T$ is the Dulac map considered 
         in \propc{coeficients_Loud1}. By applying this result we know that the first non-vanishing coefficient 
         in the asymptotic expansion of $P(s;\np)$ at $s=0$ is the third one for all $\np\in \ell_2.$ Therefore
         \cite[Theorem C]{MV20b} implies that
         $\mathcal Z_0(P'(\,\cdot\,;\nu),\nu_0)\leqslant 1$ for all $\nu_0\in \ell_2.$ 
         On account
         of this, by assertion $(2a)$ in \lemc{CZ} it follows that 
         $\mathrm{Crit}\bigl((\out_{\np_0},X_{\nu_0}),X_{\np}\bigr)\leqslant 1$ for all $\nu_0\in \ell_2.$ 
         On the other hand, due to
         $\ell_2\subset\Gamma_B$, we know by \cite[Theorem A]{MMV2} that these parameters are local
         bifurcation values of the period function at the outer boundary. Thus, since the period annuli of
         the Loud's centers vary continuously (see \obsc{+CZ}), by applying $(a)$ in \lemc{BZ} we get that
         $\mathrm{Crit}\bigl((\out_{\np_0},X_{\nu_0}),X_{\np}\bigr)\geqslant 1$ for all $\nu_0\in \ell_2.$ 
         Hence $\mathrm{Crit}\bigl((\out_{\np_0},X_{\nu_0}),X_{\np}\bigr)=1$ for all $\nu_0\in \ell_2.$
         
         We turn now to the criticality in the segment $\ell_3.$ So let us fix any 
         $\nu_0=(D_0,\frac{1}{2})$ with $D_0\in (-1,-\frac{1}{2})$ and note that then, 
         by $(d)$ in \propc{coeficients_Loud1},
         \[
         T'(s;\np)=-\rho_4(\np)(F-1/2)^2\big(\lambda\omega_{1-\lambda}(s)-1\big)
                         +\rho_5(\np)(D+1/2)+\rho_6(\np)(F-1/2)+\mathscr R(s;\nu),
         \]    
        where $\mathscr R\in\F_{1-\upsilon}^\infty(\np_0)$ for all $\upsilon>0$ small enough. To obtain
        the derivative of the Dulac time,
        we use that $\partial_s(s\omega_\alpha(s))=(1-\alpha)\omega_\alpha(s)-1$ and that, 
        by $(f)$ in Lemma~A.3 in \cite{MV20a}, 
        $\partial_s\F_{2-\upsilon}^\infty(\np_0)\subset\F_{1-\upsilon}^\infty(\np_0).$ From this equality, 
        since $\lambda\omega_{1-\lambda}(s)-1$ tends to $+\infty$ as $(s,\nu)\to (0,\nu_0)$ due 
        to $\lambda(\nu_0)=1$ (see \defic{defi_comp}), $\mathscr R(s;\nu)$ tends to 0 as $(s,\nu)\to (0,\nu_0)$, 
        $\rho_i(\nu_0)>0$ and $D_0+\frac{1}{2}<0$, we can 
        assert the existence of an open neighbourhood~$V$ of~$\nu_0$ and $\varepsilon>0$ such that 
        $P'(s;\nu)=2T'(s;\nu)<0$ for all $\nu\in V$ and $s\in (0,\varepsilon).$ Consequently 
        $\mathcal Z_0(P'(\,\cdot\,;\nu),\nu_0)=0$ and so, by applying $(2c)$ in \lemc{CZ}, we conclude that
        $\mathrm{Crit}\bigl((\out_{\np_0},X_{\np_0}),X_{\np}\bigr)=0.$
        
\item We turn next to study the horizontal segments 
         $\ell_4\!:=\big\{F=2,D\in(-2,0)\setminus\{-\frac{1}{2}\}\big\}$ 
         and the curve 
         $\ell_5\!:=\big\{D=\mathcal G(F): F\in(1,\frac{3}{2})\big\}.$ Here we set 
         $\nu_\star\!:=(\mathcal G(4/3),4/3)$ because this parameter yields to a distinguished case.
        
         We begin by noting (see the first paragraph in \secc{subsec_2.2}) that $\sigma(s;\nu)=(p_1-s,0)$ 
         is a parametrization of the 
         outer boundary of the period annulus verifying the assumptions in \lemc{CZ} and that if we denote the period of 
         the periodic orbit of $X_\nu$ 
         passing through $\sigma(s;\nu)$ by $P(s;\nu)$ then $P(s;\nu)=2T(s;\nu)$, 
         where $T$ is the Dulac map considered 
         in \propc{prop52}. Thus, by applying first that result and then \cite[Theorem C]{MV20b} we obtain that
         $\mathcal Z_0(P'(\,\cdot\,;\nu),\nu_0)\leqslant 1$ for all $\nu_0\in\ell_4\cup\ell_5\setminus\{\nu_\star\}$.         
         \mbox{Moreover~\cite[Theorem A]{MMV2}} shows that these parameters are local
         bifurcation values of the period function at the outer boundary because
         $\ell_4\cup\ell_5\subset\Gamma_B$. Since the period annuli of
         the Loud's centers vary continuously (see \obsc{+CZ}), by applying $(a)$ in \lemc{BZ} we have
         $\mathrm{Crit}\bigl((\out_{\np_0},X_{\nu_0}),X_{\np}\bigr)\geqslant 1$ for all 
         $\nu_0\in\ell_4\cup\ell_5\setminus\{\nu_\star\}.$ 
         Therefore $\mathrm{Crit}\bigl((\out_{\np_0},X_{\nu_0}),X_{\np}\bigr)=1$ for all 
         $\nu_0\in\ell_4\cup\ell_5\setminus\{\nu_\star\}.$ 
        
         On the other hand we have that
         $\mathrm{Crit}\bigl((\out_{\np_\star},X_{\nu_\star}),X_{\np}\bigr)=2$ by assertion~$(a)$ in \propc{doble}. 
         Finally the fact that there is a curve of local 
         bifurcation values of the period function at the interior of $\PA$ arriving at $\nu=\nu_\star$ 
         tangent to $\Gamma_B$ follows from assertion $(b)$ in the same result. 
         
\item Next we analyze the parameters in the segment $\ell_6\!:=\{F=1, D\in(-1,0)\}$, that corresponds to a case in 
         which there is a saddle-node singularity at the outer boundary of the period annulus. 
         This is treated in \secc{subsec_2.3}, where we introduce the map 
         \[
         \sigma(s;\np)\!:=\left\{
          \begin{array}{cl}
            (1-s,0) & \text{ if $F\leqslant 1,$} \\[2pt]
            (p_1-s,0) & \text{ if $F> 1,$}
         \end{array}
        \right.
        \]
        that provides a parametrization of the outer boundary of the period annulus verifying the assumptions
        in \lemc{CZ}. In addition if we denote by $P(s;\nu)$ the period of the periodic orbit of $X_\nu$ passing 
        through $\sigma(s;\nu)$ then $P(s;\nu)=2T(s;\nu)$, where $T$ is the Dulac map considered in \propc{prop53}. 
        From that result we get the existence of an open neighbourhood $\mathscr U$ of $\ell_6=(-1,0)\times\{1\}$ 
        such that
        \[
         P'(s;\np)=2T_{1}(\np)+4T_{2}(\np)s+s\hat h(s;\np)
        \]
        with $T_1,T_2\in\cc^0(\mathscr U)$ and where $\hat h(s;\np)$ and $s\partial_s\hat h(s;\np)$ 
        tend to zero as $s\to 0^+$
        uniformly on compact subsets of $\mathscr U$. We know         
        moreover that $T_{1}(\nu)=0$ if, and only if, $\nu=\nu_\star\!:=(-\frac{1}{2},1)$ 
        and that $T_{2}(\nu_\star)\neq 0.$
        
        If we take any $\nu_0\in\ell_6\setminus\{\nu_\star\}$ then, thanks to the good properties of the remainder,
        we get that $\lim_{(s,\nu)\to (0,\nu_0)}P'(s;\nu)=2T_1(\nu_0)\neq 0$ and this easily implies 
        $\mathcal Z_0(P'(\,\cdot\,;\nu),\nu_0)=0.$ Hence, by applying assertion $(2c)$ in \lemc{CZ}, 
        $\mathrm{Crit}\bigl((\out_{\np_0},X_{\nu_0}),X_{\np}\bigr)=0$ for all $\nu_0\in\ell_6\setminus\{\nu_\star\}$.
        
        In order to study the criticality of $\nu_\star$ we observe that 
        $\lim_{(s,\nu)\to (0,\nu_\star)}P''(s;\nu)=4T_2(\nu_\star)\neq 0$ and, consequently, 
        $\mathcal Z_0(P'(\,\cdot\,;\nu),\nu_0)\leqslant 1$ by Rolle's Theorem. Therefore, by assertion $(2a)$ in \lemc{CZ}, 
        $\mathrm{Crit}\bigl((\out_{\np_\star},X_{\nu_\star}),X_{\np}\bigr)\leqslant 1.$ On the other hand, the application
        of \cite[Theorem A]{MMV2} together with~$(a)$ in \lemc{BZ} shows that 
        $\mathrm{Crit}\bigl((\out_{\np_\star},X_{\nu_\star}),X_{\np}\bigr)\geqslant 1$ due to $\nu_\star\in\Gamma_B.$
        Hence $\mathrm{Crit}\bigl((\out_{\np_\star},X_{\nu_\star}),X_{\np}\bigr)=1$.
        
\item We proceed with the study of the segment $\ell_7\!:=\{F=0, D\in(-1,0)\}$ 
         which, as in the previous case, corresponds to 
         period annuli having a saddle-node singularity at the outer boundary. In order to compute the
         criticality of any $\nu_0\in\ell_7$ we apply the results obtained in \cite{MSV}. In that paper it is proved that
         for each $\nu_0\in \ell_7$ there exist $\delta>0$, an open neighbourhood $V$ 
         of $\nu_0$ and a continuous function
         \map{\sigma}{[0,\delta)\times V}{\RP^2} verifying the hypothesis in \lemc{CZ}. Moreover, denoting 
         the period of the periodic orbit of $X_\nu$ passing through $\sigma(s;\nu)$ by $P(s;\nu)$, the proof of
         \cite[Theorem B]{MSV} shows that $P'(s;\nu)$ tends to $-\infty$ as $(s,\nu)\to (0,\nu_0).$ Consequently 
         $\mathcal Z_0(P'(\,\cdot\,;\nu),\nu_0)=0$ and hence, by applying $(c)$ in \lemc{CZ}, we get that     
         $\mathrm{Crit}\bigl((\out_{\np_0},X_{\nu_0}),X_{\np}\bigr)=0$.   
         
\item We analyze next the parameters inside the segment $\ell_8\!:=\big\{D=0,F\in [\frac{1}{4},2]\big\}.$ So let us
         fix any $\nu_0=(0,F_0)$ with $F_0\in [\frac{1}{4},2].$ By \cite[Theorem A]{MMV2} we can assert
         that if $F_0\in [\frac{1}{2},2]$ then $\nu_0$ is a local bifurcation value of the period function at the outer
         boundary. On the other hand, if $F_0\in [\frac{1}{4},\frac{1}{2}]$ then we can conclude the same  by applying
         \cite[Theorem B]{MV}. Hence, since the period annuli of the 
        Loud's centers vary continuously, see \obsc{+CZ}, assertion $(a)$ in \lemc{BZ} shows that 
        $\mathrm{Crit}\bigl((\out_{\np_0},X_{\nu_0}),X_{\np}\bigr)\geqslant 1$ for all $\nu_0\in\ell_8.$
        
\item From the results in \cite{TopaV,Jordi1} it follows that 
         $\mathrm{Crit}\bigl((\out_{\np_0},X_{\np_0}),X_{\np}\bigr)=0$ for any parameter $\nu_0$ inside the  
         set
         $\ell_9\!:=\{D=0,F\notin[0,2]\}\cup\{D=-1,F<0\}\cup\{D+F=0,F<0\}$. 
         Indeed, in those papers the authors determine 
         a region $M$ in the parameter plane for which the corresponding centers have a globally monotonic period 
         function. Taking this into account the assertion follows easily from the fact that 
         $\ell_9$ is contained in the interior of $M$, see \defic{criticality}.
         
\item We consider now the half-line $\ell_{10}\!:=\{D+F=0,F>1\},$ 
         so let us take a parameter $\nu_0=(-F_0,F_0)$ with $F_0>1.$
         In this case the assertions with regard to its criticality 
         follow from the results in~\cite{MMV20}.            
         It is proved there that there exists a function $\xi=\xi(\nu)$ in a neighbourhood $U$ of~$\nu_0$ 
         such that $\sigma(s;\nu)=\big(0,\xi(\nu)(1-s)\big)$ is a $\cc^0$ map on 
         $[0,\delta)\times U$ verifying the hypothesis 
         $(a)$, $(b)$ and $(c)$ in \lemc{CZ}. Moreover if we denote by $P(s;\nu)$ the period of 
         the periodic orbit of $X_\nu$ 
         passing through $\sigma(s;\nu)$ then \cite[Theorem B]{MMV20} shows that 
         \begin{itemize}
         \item $\mathcal Z_0(P'(\,\cdot\,;\nu),\nu_0)=0$ if $F_0\notin [3/2,2]$,                  
         \item $\mathcal Z_0(P'(\,\cdot\,;\nu),\nu_0)=1$ if $F_0\in [3/2,2)$ and
         \item $\mathcal Z_0(P'(\,\cdot\,;\nu),\nu_0)=2$ if $F_0=2$.  
         \end{itemize}
         In the first case 
         $\mathrm{Crit}\bigl((\out_{\np_0},X_{\np_0}),X_{\np}\bigr)=0$ by $(2c)$ in \lemc{CZ}, whereas 
         in the second case the combination of $(2a)$ and $(2b)$
         implies $\mathrm{Crit}\bigl((\out_{\np_0},X_{\np_0}),X_{\np}\bigr)=1.$ 
         In the third case, by applying~$(2a)$ 
         we get $\mathrm{Crit}\bigl((\out_{\np_0},X_{\np_0}),X_{\np}\bigr)\leqslant 2.$
         To show that this upper bound is attained we also apply $(2b)$ in \lemc{CZ} but to this end we must 
         check the assumption that for each open neighbourhood~$V$ of $\nu_0=(-2,2)$ and $\delta>0$ there exist
         distinct $s_1,s_2\in (0,\delta)$ and $\hat\nu\in V$ such that $P'(s_i;\hat\nu)=0$ for $i=1,2.$ To verify this
         we note first, see \cite[\S 4]{MMV20}, that we can write 
         \[
          P'(s;\nu)=\delta_1(\nu)f_1(s;\nu)+\delta_2(\nu)f_2(s;\nu)+f_3(s;\nu)
         \]
         where the coefficients $\delta_1$ and $\delta_2$ are independent at $\nu_0$ in the sense of 
         \cite[Definition 4.1]{MV20b} and, for $i=1,2,$ $\lim_{s\to 0^+}\frac{f_{i+1}(s;\nu)}{f_i(s;\nu)}=0$.
         On account of this and $P'(s;\nu_0)\not\equiv 0,$ the fact that the mentioned assumption is verified 
         follows from the proof of \cite[Proposition 4.2]{MV20b}. Related with this let us also mention, 
         see again \cite[\S 4]{MMV20}, that the ordered set
         $(f_1,f_2,f_3)$ is an extended complete Chebyshev system 
         on $(0,\varepsilon)$ 
         for $\varepsilon>0$ sufficiently small (see \cite{KS} for a definition). 
               
         On the other hand, by assertion $(a)$ in \cite[Theorem C]{MMV20}, there exist a neighbourhood $U$
         of $\nu_0=(-2,2)$, $s_0>0$ and a injective $\cc^0$ curve \map{\rho}{(-\varepsilon,\varepsilon)}{U} satisfying 
         $\rho(0)=(-2,2)$ and
         \[
          \rho((0,\varepsilon))=\Delta\!:=\{\nu\in U;\text{ there exists $s\in (0,s_0)$ such that $P'(s;\nu)=P''(s;\nu)=0$}\}.
         \] 
         Furthermore, assertion $(b)$ in that result shows that the curve $\Delta$ has an exponentially flat contact 
         with the straight line $\{F=2\}$ at $\nu_0=(-2,2)$, see \figc{doblebis}. 
         Since the interior of $\Delta$ is clearly empty and 
         the Chebyshev property explained above prevents the zeros of 
         $P'(\,\cdot\,;\nu)$ to have multiplicity greater than 2, the application of \lemc{lema-interior} shows 
         that~$\Delta$ consists of local bifurcation values of the period function at the interior. 
         
\item Finally the fact that the criticality at the outer boundary of the isochrones 
         $\nu_1\!:=(-\frac{1}{2},2)$ and $\nu_2\!:=(-\frac{1}{2},\frac{1}{2})$ is 1 follows from Propositions~\ref{I2} 
         and~\ref{I4}, respectively.
         
\end{enumerate}
Since $\R^2\setminus\ell_0=\left(\cup_{i=1}^{10}\ell_i\right)\cup\{\nu_1\}\cup\{\nu_2\},$ this concludes the proof of the result.
\end{prooftext}

\begin{figure}[t]
\begin{center}
\begin{tikzpicture}
\node[inner sep=0pt] at (7,0) {\includegraphics[width=14cm]{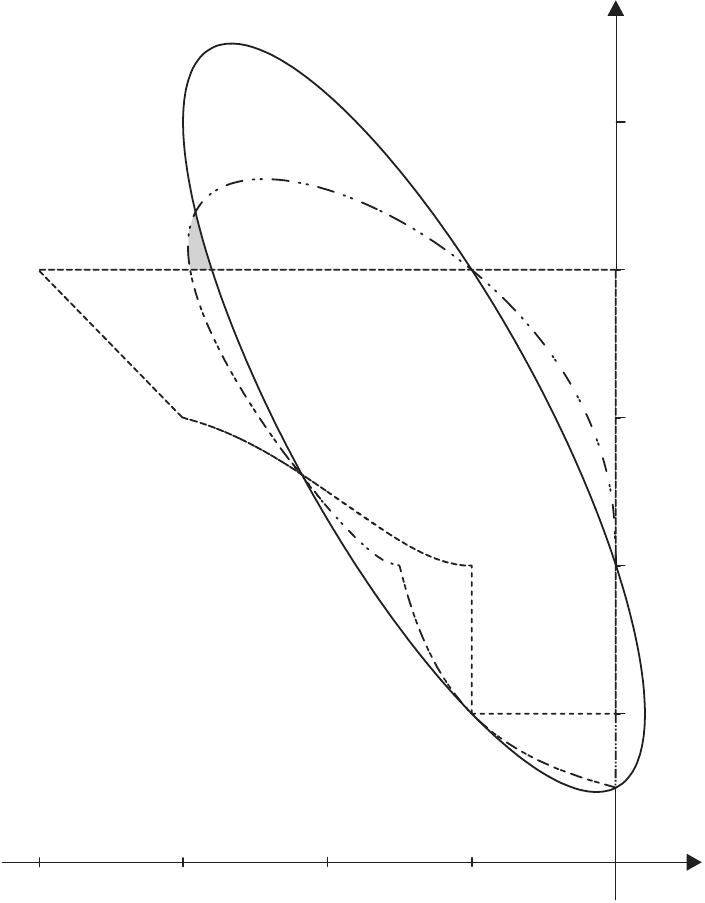}};
\node at (4.6,0.75) {$\Gamma_B$};
\node at (7.5,3) {$\Gamma_C$};
\node at (6.75,2.55) {$\Gamma_0$};
\node at (6.75,-4) {$-1$};
\node at (4.25,-4) {$-2$};
\node at (10,-4) {$D$};
\node at (9.75,4) {$F$};
\node at (9.75,-1) {$1$};
\node at (9.75,1.55) {$2$};
\node at (4.6,1.8) {$\Delta$};
\node at (5.3,2.7) {$\delta_1$};
\node at (5.8,2.88) {$L_1$};
\draw[thick]  (4.279,1.58) to [out=0,in=190] 
(5.1,1.62)
to [out=10,in=-90] (5.53,2.88);
\draw[fill] (5.53,2.88) circle [radius=0.04];
\draw[fill]  (4.279,1.58) circle [radius=0.04];
\draw[fill,white] (5.3,2.1) circle [radius=0.25];
\end{tikzpicture}
\end{center}
\caption{Arrangement of the three types of local bifurcation curve (inner boundary, interior and outer boundary) near the parameter values $\nu=L_1$, see \refc{Li}, and $\nu=(-2,2)$. We refer the reader to \obsc{figdoblebis_com} for a detailed explanation.}\label{doblebis}
\end{figure}

\begin{obs}\label{figdoblebis_com} 
We conclude this section by making further comments about \figc{doblebis}. It follows from~$(e)$ in \teoc{Loud} and \cite[Theorem 4.3]{MMV2}, respectively, that the germs of curve $\Delta$ at $\nu=(-2,2)$ and $\delta_1$ at $\nu=L_1$ are inside the set of local bifurcation values of the period function at the interior of~$\PA.$ Exactly as we explain in \obsc{figdobe_com}, we conjecture that both curves connect each other to delimit 
a region of parameters for which the period function has exactly two critical periodic orbits. In this regard
 \cite[Theorem 5.2]{MMV2} shows that the center of any parameter inside the light gray sector has at least two critical periodic orbits. We know now, see point 8 in the proof of \teoc{Loud}, that $\nu=(-2,2)$ is at the boundary of this region with exactly two critical periodic orbits. The numerical visualization of this fact is a challenging problem because~$\Delta$ has a 
exponential flat contact with $\{F=2\}$ at $\nu=(-2,2)$.
\end{obs}


\appendix


\section{Coefficient formulas}\label{appA}

\subsection{Previous results about the Dulac time}\label{appA1}

This appendix is entirely devoted to the proof of Propositions~\ref{coeficients_Loud1} and~\ref{prop52} in \secc{sec2}. For the parameter values under consideration in both results, and thanks to the symmetry of the vector field $X_\nu$ in~\refc{sist_loud}, it turns out that the period function is twice the Dulac time associated to the passage through a hyperbolic saddle at infinity. The asymptotic expansion of this type of passage is the subject of our recent papers \cite{MV20a,MV20b,MV21} and in order to prove the results in \secc{sec2} we strongly rely on the tools developed there. For this reason we first summarize for reader's convenience the definitions and results from those papers that are indispensable here. We recap the results in three theorems. In short, \teoc{5punts} will provide us with the monomial scale needed in each asymptotic expansion, which only depends on the hyperbolicity ratio of the saddle, whereas \teoc{oldA} will give the explicit expression of their coefficients in terms of a sort of Mellin transform that is introduced in \teoc{L8}. 

\begin{figure}[t]
   \centering
  \begin{lpic}[l(0mm),r(0mm),t(0mm),b(5mm)]{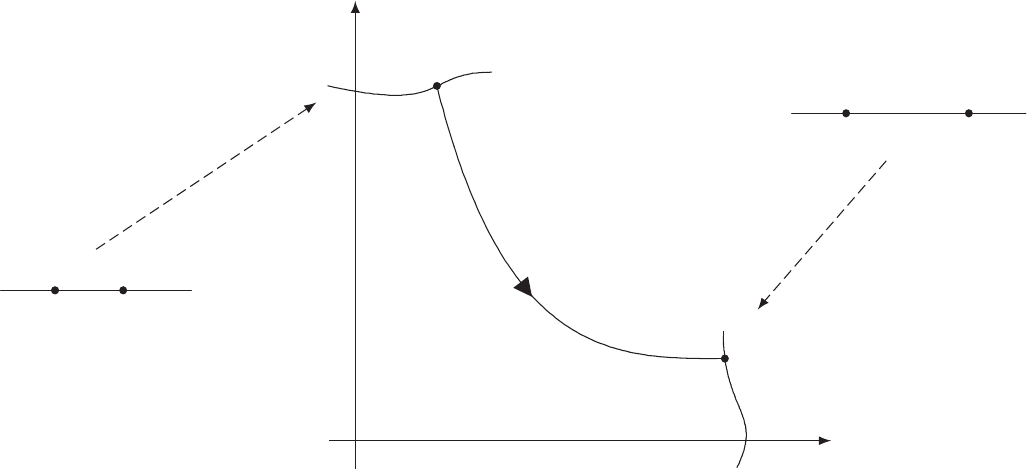}
   \lbl[l]{5,15.5;$0$}   
   \lbl[l]{12,16;$s$}   
   \lbl[l]{18,32;$\sigma_1$}
   \lbl[l]{28,40;$\Sigma_1$}      
   \lbl[l]{31.5,45;$x_2$}   
   \lbl[l]{38,42;$\sigma_1(s)$}   
   \lbl[l]{52,25;$\varphi(\,\cdot\,,\sigma_1(s))$}
   \lbl[l]{75.5,12;$\varphi(T(s),\sigma_1(s))$}  
   \lbl[l]{75,-2;$\Sigma_2$}   
   \lbl[l]{82,1;$x_1$}
   \lbl[l]{80,26;$\sigma_2$}
   \lbl[l]{85,33.5;$0$}
   \lbl[l]{98,33.5;$s$}                   
   \end{lpic}
  \caption{Definition of $T(\,\cdot\,;\no)$, where $\varphi(t,p;\no)$ is the solution of~$X_\no$ passing through the point $p\in\mathscr U$ at time $t=0.$}
  \label{DefTyR}
 \end{figure}

In order to facilitate the application of the above-mentioned results we particularize them to fit in the context needed to prove Propositions~\ref{coeficients_Loud1} and~\ref{prop52}. Thus, following the notation that we use in \cite{MV20b}, let us consider the parameter $\no\!:=(\lambda,\mu)\in\hat W\!:=(0,+\infty)\times W$, where $W$ is an open set of $\R^N$, and the family of vector fields $\{X_{\no}\}_{\no\in\hat W}$ with
\begin{equation}\label{X}
 X_\no({x_1},{x_2})\!:=\frac{1}{x_2}\Big({x_1}P_1({x_1},{x_2};\no)\partial_{x_1}+{x_2}P_2({x_1},{x_2};\no)\partial_{x_2}\Big),
\end{equation}
where
\begin{itemize}
\item $P_1$ and $P_2$ belong to $\mathscr C^{\omega}(\mathscr U\!\times\!\hat W)$ for some open set $\mathscr U$ of $\R^2$ containing the origin, 
\item $P_1({x_1},0;\no)>0$ and $P_2(0,{x_2};\no)<0$ for all $({x_1},0),(0,{x_2})\in\mathscr U$ and $\no\in \hat W,$
\item $\lambda=-\frac{P_2(0,0;\no)}{P_1(0,0;\no)}$.
\end{itemize}
Moreover, for $i=1,2,$ let \map{\sigma_i}{(-\varepsilon,\varepsilon)\times \hat W}{\Sigma_i} be a $\mathscr C^{\omega}$~transverse section to~$X_{\no}$ at $x_i=0$ defined by
 \[
  \sigma_i(s;\no)=\bigl(\sigma_{i1}(s;\no),\sigma_{i2}(s;\no)\bigr)
 \]
such that $\sigma_1(0,\no)\in\{(0,x_2);x_2>0\}$ and $\sigma_2(0,\no)\in\{(x_1,0);x_1>0\}$
for all $\no\in \hat W.$ Then \teoc{5punts} is concerned with the time $T(s;\no)$ that spends the solution of $X_\no$ passing through the point $\sigma_1(s;\no)\in\Sigma_1$ to arrive at~$\Sigma_2,$ see \figc{DefTyR}. More concretely it  shows that $T(s;\no)$ has an asymptotic expansion at $s=0$ with the remainder having good flatness properties with respect to the parameters. We specify these properties in the following two definitions. 

\begin{defi}\label{defi_fun}
Given $K\in\Z_{\geq 0}\cup\{\infty\}$ and an open subset $U\subset\hat W\subset \R^{N+1},$ we say that a function $\psi(s;\no)$ belongs to the class $\mathscr C^K_{s>0}(U)$ if there 
exists an open neighbourhood $\Omega$ of 
\[ 
 \{(s,\no)\in\R^{N+2};s=0,\no\in U\}=\{0\}\times U
\]  
in $\R^{N+2}$ such that $(s,\no)\mapsto \psi(s;\no)$ is $\mathscr C^K$ on $\Omega\cap\big((0,+\infty)\times U\big)$.     
\end{defi}

\begin{defi}\label{defi2} 
Consider $K\in\Z_{\geq 0}\cup\{\infty\}$ and an open subset $U\subset\hat W\subset\R^{N+1}.$ Given $L\in\R$ and $\no_0\in U$, we say that $\psi(s;\no)\in\mathscr C^K_{s>0}(U)$ is \emph{$(L,K)$-flat with respect to $s$ at $\no_0$}, and we write $\psi\in\F_L^K(\no_0)$, if for each $\nu=(\nu_0,\ldots,\nu_{N+1})\in\Z_{\geq 0}^{N+2}$  with $|\nu|=\nu_0+\cdots+\nu_{N+1}\leqslant K$ there exist a neighbourhood~$V$ of~$\no_0$ and $C,s_0>0$ such that
\begin{equation*}
 \left|\frac{\partial^{|\nu|}\psi(s;\no)}{\partial s^{\nu_0}
 \partial\no_1^{\nu_1}\cdots\partial\no_{N+1}^{\nu_{N+1}}}\right|\leqslant C s^{L-\nu_0}
 \text{ for all $s\in(0,s_0)$ and $\no\in V$.} 
\end{equation*}
If $W$ is a (not necessarily open) subset of $U$ then define $\F_L^K(W)\!:=\bigcap_{\no_0\in W}\F_L^K(\no_0).$
\end{defi}

Next result merges the statements of Theorem~1.6, Theorem~4.3 and Corollary~B in \cite{MV21}. Following the notation in that paper, we particularise them to the case $(n_1,n_2)=(0,1)$ for simplicity. Moreover, for the sake of shortness, we only include those items that will be used in the present paper. 

\begin{theo}\label{5punts}

Let $T (s;\no)$ be the Dulac time of the hyperbolic saddle \refc{X} from $\Sigma_1$ and $\Sigma_2$. Then, setting $D_{00}=\emptyset$, $D_{10}=\frac{1}{\N}$, $D_{01}=\N$, $D_{20}=\frac{2}{\N}$ and $D_{02}=\frac{\N}{2},$ for each $(i,j)\in\{(0,0),(1,0),(0,1),(2,0),(0,2)\}$ there exists a meromorphic function $T_{ij}(\no)$ on $\hat W=(0,+\infty)\times W$, having poles only along $D_{ij}\times W$, such that the following assertions hold:

\begin{enumerate}[$(1)$]

\item If $\lambda_0\in(1,2)$ then $T(s;\no)=T_{00}(\no)+T_{10}(\no)s+T_{01}(\no)s^{\lambda}+T_{20}(\no)s^2+\F_{L}^\infty(\{\lambda_0\}\times W)$ for any $L\in\big[2,\lambda_0+1\big)$.

\item If $\lambda_0>2$ then $T(s;\no)=T_{00}(\no)+T_{10}(\no)s+T_{20}(\no)s^2+\F_{L}^\infty(\{\lambda_0\}\times W)$ for any $L\in \big[2,\min(3,\lambda_0)\big)$.

\item If $\lambda_0=\frac{1}{2}$ then $T(s;\no)=T_{00}(\no)+T_{01}(\no)s^{\lambda}+s\T_{10}^{\lambda_0}(\omega;\no)+\F_L^\infty(\{\lambda_0\}\times W)$ for any $L\in [1,\frac{3}{2})$, where $\omega=\omega(s;\alpha)$, $\alpha=1-2\lambda$ and $\T_{10}^{\lambda_0}(w;\no)\in\cc^\infty(\hat U)[w]$ for some open neighbourhood $\hat U$ of $\{\lambda_0\}\times W$. Moreover 
\[
 \T_{10}^{\lambda_0}(\omega;\no)=T_{10}(\no)+T_{02}(\no)(1+\alpha\omega)\,\text{ for $\lambda\neq\lambda_0$}.
\]  

\item If $\lambda_0=1$ then
$T(s;\no)=T_{00}(\no)+s\T_{10}^{\lambda_0}(\omega;\no)+\F_{L}^\infty(\{\lambda_0\}\times W)$ for any $L\in (1,2)$, where $\omega=\omega(s;\alpha)$, $\alpha=1-\lambda$ and $\T_{10}^{\lambda_0}(w;\no)\in\cc^\infty(\hat U)[w]$ for some open neighbourhood $\hat U$ of $\{\lambda_0\}\times W$. Moreover 
\[
 \T_{10}^{\lambda_0}(\omega;\no)=T_{10}(\no)+T_{01}(\no)(1+\alpha\omega)\,\text{ for $\lambda\neq\lambda_0$}.
\]

\item If $\lambda_0=2$ then $T(s;\no)=T_{00}(\no)+T_{10}(\no)s+s^2\T_{20}^{\lambda_0}(\omega;\no)+\F_L^\infty(\{\lambda_0\}\times W)$ for any $L\in \big[2,3\big)$, where $\omega=\omega(s;\alpha)$, $\alpha=2-\lambda$ and $\T_{20}^{\lambda_0}(w;\no)\in\cc^\infty(\hat U)[w]$ for some open neighbourhood $\hat U$ of $\{\lambda_0\}\times W$. Moreover 
\[
\T_{20}^{\lambda_0}(\omega;\no)=T_{20}(\no)+T_{01}(\no)\left(1+\alpha\omega\right)\,\text{ for $\lambda\neq\lambda_0$}.
\]

\end{enumerate}
\end{theo}

We focus next on the expression of the coefficients $T_{ij}$ and the result that we state below in this regard follows from assertion~$(c)$ in \cite[Theorem A]{MV21} particularized to $(n_1,n_2)=(0,1)$. In its statement we use the following functions:
\begin{equation}\label{def_fun}
\begin{array}{ll}
\dsp L_1(u)\!:=\exp\int_0^u\left(\frac{P_1(0,z)}{P_2(0,z)}+\frac{1}{\lambda}\right)\frac{dz}{z} & 
\dsp L_2(u)\!:=\exp\int_0^u\left(\frac{P_2(z,0)}{P_1(z,0)}+{\lambda}\right)\frac{dz}{z} 
\\[15pt]
\dsp A_1(u)\!:=\frac{1}{P_2(0,u)}&  
\dsp A_2(u)\!:=\frac{L_2(u)}{P_1(u,0)}
\\[15pt]
\dsp M_1(u)\!:=L_1(u)\partial_1\!\left(\frac{P_1}{P_2}\right)(0,u)& 
\dsp B_1(u)\!:=L_1(u)\partial_1P_2^{-1}(0,u)
\\[15pt]
\dsp C_1(u)\!:=L_1^2(u)\partial_1^2P_2^{-1}(0,u)
& 
\dsp \hspace{-1.8truecm}+2L_1(u)\gorro{M}_1(1/\lambda,u)\partial_1P_2^{-1}(0,u)
\end{array}
\end{equation}
Here, given $\alpha\in\R\setminus\Z_{\ge 0}$ and a real valued function $f(x)$ that is $\cc^\infty$ in an open interval containing $x=0$, $\hat f(\alpha,x)$ is a sort of incomplete Mellin transform (see \teoc{L8} below). Moreover, for the sake of shortness, in the following statement  we use the compact notation $\sigma_{ijk}$ for the $k$th derivative at $s=0$ of the $j$th component of $\sigma_i(s;\no)$, i.e., 
\[
 \sigma_{ijk}(\no)\!:=\partial^k_s\sigma_{ij}(0;\no).
\]
Also with regard to the statement, note that $D_{ij}$ refers to the discrete sets introduced in \teoc{5punts}.

\begin{theo}\label{oldA}
For each $(i,j)\in\{(0,0),(1,0),(0,1),(2,0)\}$, the following expression of $T_{ij}(\no)$ is valid provided that $\lambda\notin D_{ij}$:
\begin{align*}
T_{00}(\no)
\!=&-\sigma_{120}\gorro{A}_1(-1,\sigma_{120}),\\[8pt]
T_{01}(\no)
\!=&\frac{\sigma_{120}\sigma_{111}^\lambda}{\sigma_{210}^\lambda L_1^\lambda(\sigma_{120})}\gorro A_2( \lambda,\sigma_{210}),\\[8pt]
T_{10}(\no)
\!=&-\frac{\sigma_{121}}{P_2(0,\sigma_{120})}-\frac{\sigma_{120}\sigma_{111}}{L_1(\sigma_{120})}\gorro{B}_1\big(1/\lambda-1,\sigma_{120}\big),\\[8pt]
%
\intertext{and}
T_{20}(\no)
\!=&- \frac{\sigma_{120}\sigma_{122}}{2\sigma_{120}P_2(0,\sigma_{120})}-\frac{1}{2}\sigma_{121}^2\partial_2P_2^{-1}(0,\sigma_{120})
-\sigma_{121}\sigma_{111}
\partial_1P_2^{-1}(0,\sigma_{120})
\\
&
-\frac{\sigma_{120}\sigma_{111}^2}{2L_1^2(\sigma_{120})}\gorro{C}_1(2/\lambda-1,\sigma_{120})-S_1\frac{\sigma_{120}\sigma_{111}}{L_1(\sigma_{120})}\gorro{B}_1(1/\lambda-1,\sigma_{120}),
\intertext{where}
S_1
\!=&\frac{\sigma_{112}}{2\sigma_{111}}-\frac{\sigma_{121}}{\sigma_{120}}\left(\frac{P_1}{P_2}\right)\!(0,\sigma_{120})-\frac{\sigma_{111}}{L_1(\sigma_{120})}\gorro{M}_1(1/\lambda,\sigma_{120}).
\end{align*}
\end{theo}

As we already explained, the following result (that merges Theorem~B.1 and Corollary~B.3 in \cite{MV21}) is the third ingredient needed in the proof of Propositions~\ref{coeficients_Loud1} and~\ref{prop52}.

\begin{theo}\label{L8}
Let us consider an open interval $I$ of $\R$ containing $x=0$ and an open subset $U$ of $\R^N$.
\begin{enumerate}[$(a)$]

\item Given $f(x;\nu)\in\mathscr C^{\infty}(I\times U)$, there exits a unique $\hat f(\alpha,x;\nu)\in\mathscr C^{\infty}((\R\setminus\Z_{\ge 0})\times I\times U)$ such that 
\begin{equation*}
 x\partial_x\hat f({\alpha},x;\nu)-\alpha\hat f({\alpha},x;\nu)=f(x;\nu).
\end{equation*}

\item If $x\in I\setminus\{0\}$ then $\partial_x(\hat f({\alpha},x;\nu)|x|^{-\alpha})=f(x;\nu)\frac{|x|^{-\alpha}}{x}$ and, taking any $k\in\Z_{\ge0}$ with $k>\alpha$,
\begin{equation*}
\hat f(\alpha,x;\nu)=
\sum_{i=0}^{k-1}\frac{\partial_x^if(0;\nu)}{i!(i-\alpha)}x^i+|x|^{\alpha}\int_0^x\!\left(f(s;\nu)-T_0^{k-1}f(s;\nu)\right)|s|^{-\alpha}\frac{ds}{s},
\end{equation*}
where $T_0^kf(x;\nu)=\sum_{i=0}^{k}\frac{1}{i!}\partial_x^if(0;\nu)x^i$ is the $k$-th degree Taylor polynomial of $f(x;\nu)$ at $x=0$.


\item If $f(x;\nu)$ is analytic on $I\times U$ then $\hat f(\alpha,x;\nu)$ is analytic on $(\R\setminus\Z_{\ge 0})\times I\times U$. Finally,
for each $(\alpha_0,x_0,\nu_0)\in\Z_{\ge 0}\times I\times U$ the function
$(\alpha,x,\nu)\mapsto(\alpha_0-\alpha)\hat f(\alpha,x;\nu)$ extends analytically to $(\alpha_0,x_0,\nu_0)$.

\item If $f(x;\nu)=x^ng(x;\nu)$ with $g\in \cc^{\infty}(I\times U)$ and $n\in\N$ then 
         $\hat f(\alpha,x;\nu)=x^n\hat g(\alpha-n,x;\nu).$ 

\end{enumerate}
\end{theo}

The following simple observation will be useful in order to study the coefficients of the asymptotic expansions that we shall deal with.

\begin{obs}\label{rm1}
If $\sum_{i=1}^ma_ix^{\lambda_i}+\psi(x)=0$ for all $x\in (0,\varepsilon),$ where $\lambda_i\in\R$ with $\lambda_1<\lambda_2<\cdots<\lambda_m$, $a_1,a_2,\ldots,a_m\in\R$ and $\psi(x)=\mathrm{o}(x^{\lambda_m})$ then $a_1=a_2=\cdots=a_m=0$.
\end{obs}

We are now in position to begin the proof of the two first results in \secc{sec2}.

\subsection{Proof of \propc{coeficients_Loud1}}

\begin{prooftext}{Proof of \propc{coeficients_Loud1}.}
We follow the approach in \cite[\S 5]{MMV03} to take advantage of the general setting developed in \cite{MV20b}. To this end we will work on an extended parameter space $\bar\np\in\bar V$ that we specify as follows. Firstly, introducing an auxiliary parameter $\eta\approx 0,$ we consider two local transverse sections $\Sigma_1^\eta$ and $\Sigma_2^\eta$ parametrized respectively by $s\mapsto(1-s,\eta)$ and $s\mapsto(-\frac{1}{s},\frac{\eta}{s})$, for $s> 0$, cf. \figc{fig1}. 
Secondly, taking any $\alpha,\beta\in\R$ such the straight line $y=\alpha x+\beta$ does not intersect any solution of $X_\np$ while traveling from $\Sigma_1^\eta$ to $\Sigma_2^\eta.$ One can readily see that a sufficient condition for this to hold is that
\[
 \alpha+\beta<\eta\text{ and }\eta>-\alpha.
\]
Then, setting $\bar\np\!:=(D,F,\alpha,\beta,\eta)$, we will work on the extended parameter space
\[ 
 \bar V\!:=\left\{\bar\np\in\R^5:D\in (-1,0), F\in (0,1), \alpha+\beta<\eta,-\alpha<\eta\right\}.
\] 
Taking this into account, we consider the projective change of coordinates, see \figc{recta},
\begin{figure}[t]
 \centering
 \begin{lpic}[l(0mm),r(0mm),t(0mm),b(0mm)]{recta(1)}
  \lbl[l]{41,50;$x_1$}  
  \lbl[l]{48.5,45;$x_2$}
  \lbl[l]{43,22;$y$} 
  \lbl[l]{64,10;$x$}
  \lbl[l]{53,30;$\{x=1\}$}
  \lbl[l]{35,2.5;$L^\infty_{(x_1,x_2)}=\{y=\alpha x+\beta\}$}
  \lbl[l]{0,40;$L^\infty_{(x,y)}=\{x_2=0\}$}           
 \end{lpic}
  \caption{Projective coordinate change in the proof of \propc{coeficients_Loud1}.}\label{recta}
 \end{figure}
\begin{equation*}
 (x_1,x_2)=\left(\frac{1-x}{y-\alpha x-\beta},\frac{1}{y-\alpha x-\beta}\right).
\end{equation*}
One can verify that in these coordinates the parametrizations of $\Sigma_1^\eta$ and $\Sigma_2^\eta$ become
\begin{equation}\label{t1eq3}
 \sigma_1(s)\!:=\left(\frac{s}{\eta+\alpha s-\alpha-\beta},\frac{1}{\eta+\alpha s-\alpha-\beta}\right)
 \text{ and } 
 \sigma_2(s)\!:=\left(\frac{1+s}{\alpha+\eta-\beta s},\frac{s}{\alpha+\eta-\beta s}\right),
\end{equation}
respectively, whereas the vector field \refc{sist_loud} is brought to
\begin{equation}\label{t1eq0}
 \bar X_{\bar\np}\!:=\frac{1}{x_2}\big(x_1P_1(x_1,x_2;\bar\np)\partial_{x_1}+x_2P_2(x_1,x_2;\bar\np)\partial_{x_2}\big)
\end{equation}
with
\begin{align*}
 P_1(x_1,x_2;\bar\np)&=(1-\alpha x_1)^2+(\alpha+\beta)x_2-x_2^2+(1-\alpha^2-\alpha\beta)x_1x_2
 \\&\hspace{2truecm}-F(1+\alpha(x_2-x_1)+\beta x_2)^2-D(x_1-x_2)^2\\
\intertext{and}
 P_2(x_1,x_2;\bar\np)&=\alpha^2x_1^2-\alpha x_1-x_2^2+(1-\alpha^2-\alpha\beta)x_1x_2\\
 &\hspace{2truecm}-F(1+\alpha(x_2-x_1)+\beta x_2)^2-D(x_1-x_2)^2.
\end{align*} 
The reason why we introduce the auxiliary parameters $\alpha,$ $\beta$ and $\eta$ is because the computations are much easier taking the projective change of coordinates that sends $y=0$ to infinity (i.e., with $\alpha=\beta=0)$, which is not compatible with the placement of the original transverse sections (i.e., with $\eta=0).$ Since the parameters $\bar\np$ with $\alpha=\beta=\eta=0$ are in the boundary of the admissible set $\bar V$, we will work in the interior and then make a limit argument. By introducing these auxiliary parameters we end up in a setting where the assumptions to apply the results in Section~\ref{appA1} are fulfilled. Observe in particular that $P_1$ and $P_2$ are analytic on $\R^2\times\bar V$. Following the notation introduced there, note that the hyperbolicity ratio of the saddle
\[
\lambda=-\frac{P_2(0,0)}{P_1(0,0)}=\frac{F}{1-F}
\]
depends on $\mu=\bar\np$. We denote the Dulac time of $\bar X_{\bar\np}$ from $\Sigma_1^\eta$ to $\Sigma_2^\eta$ by $\bar T(s;\bar\np)$. Note that, by construction,  it does not depend on $\alpha$ and $\beta$ as long as $\alpha+\beta<\eta$ and $\eta>-\alpha$ holds. Moreover, and this is the key point for our purposes, the Dulac time $T(s;\np)$ in the statement is precisely $\bar T(s;\bar\np)$ for $\eta=0.$ 

Let us fix any $\bar\np_0=(D_0,F_0,\alpha_0,\beta_0,\eta_0)\in\bar V$ with $F_0\in [\frac{1}{2},1)$. Observe that $\lambda(\bar\nu_0)>2$ if $F_0\in (\frac{2}{3},1),$ $\lambda(\bar\nu_0)\in (1,2)$ if $F_0\in (\frac{1}{2},\frac{2}{3})$, $\lambda(\bar\nu_0)=2$ if $F_0=\frac{2}{3}$ and $\lambda(\bar\nu_0)=1$ if $F_0=\frac{1}{2}$. Consequently, by applying (2), (1), (5) and (4) in \teoc{5punts}, respectively, we get that
\begin{enumerate}[$(a')$]

\item $\bar T(s;\bar\np)=\bar T_{00}(\bar\np)+\bar T_{10}(\bar\np)s+\bar T_{20}(\bar\np)s^2+\F_{L_0-\upsilon}^\infty(\bar\np_0)$ if $F_0\in (\frac{2}{3},1),$ where $L_0=\min(3,\lambda(\bar\np_0)),$ 

\item $\bar T(s;\bar\np)=\bar T_{00}(\bar\np)+\bar T_{10}(\bar\np)s+\bar T_{01}(\bar\np)s^\lambda+\F^\infty_{2-\upsilon}(\bar\np_0)$ if $F_0\in (\frac{1}{2},\frac{2}{3})$,

\item $\bar T(s;\bar\np)=\bar T_{00}(\bar\np)+\bar T_{10}(\bar\np)s+\bar T_{201}^2(\bar\np)s^2\omega_{2-\lambda}(s)+\bar T_{200}^2(\bar\np)s^2+\F_{3-\upsilon}^\infty(\bar\np_0)$ if $F_0=\frac{3}{2}$, where
$\bar T_{201}^2(\bar\np)$ and $\bar T_{200}^2(\bar\np)$ are smooth in a neighbourhood of $\{\bar\np\in\bar V:\lambda(\bar\np)=2\}$ and, moreover,
\[
 \bar T_{201}^2(\bar\np)=(2-\lambda)\bar T_{01}(\bar\np)
 \text{ and }
 \bar T_{200}^2(\bar\np)=\bar T_{20}(\bar\np)+\bar T_{01}(\bar\np)
 \text{ for $\lambda(\bar\np)\neq 2,$}
\]

\item $ \bar T(s;\bar\np)=\bar T_{00}(\bar\np)+\bar T_{101}^1(\bar\np)s\omega_{1-\lambda}(s)+\bar T_{100}^1(\bar\np)s+\F_{2-\upsilon}^\infty(\bar\np_0)$ if $F_0=\frac{1}{2}$,  where
$\bar T_{101}^1(\bar\np)$ and $\bar T_{100}^1(\bar\np)$ are smooth in a neighbourhood of $\{\bar\nu\in\bar V:\lambda(\bar\nu)=1\}$ and, moreover, 
\[
 \bar T_{101}^1(\bar\np)=(1-\lambda)\bar T_{01}(\bar\np)
 \text{ and }
 \bar T_{100}^1(\bar\np)=\bar T_{10}(\bar\np)+\bar T_{01}(\bar\np)
  \text{ for $\lambda(\bar\np)\neq 1.$}
\]

\end{enumerate}
Here $\upsilon$ is a small enough positive number depending on $\bar\np_0.$ Furthermore the coefficients $\bar T_{ij}(\bar\np)$ are meromorphic functions on $\bar V$ with poles only at those $F\in (0,1)$ such that $\lambda(\bar\np)=\frac{F}{1-F}\in D_{ij},$ where $D_{00}=\emptyset$, $D_{10}=\frac{1}{\N}$, $D_{01}=\N$ and $D_{20}=\frac{2}{\N}$.

We claim that the coefficients $\bar T_{ij}(\bar\np)$ do not depend on $\alpha$ and $\beta$. Indeed, to see this recall that the Dulac time $\bar T(s;\bar\np)$ does not depend on $\alpha$ and $\beta$ provided that $\alpha+\beta<\eta$ and $\eta>-\alpha$, which is verified for $\bar\np\in\bar V.$ Hence $\partial_\alpha \bar T(s;\bar\np)\equiv 0$ and $\partial_\beta \bar T(s;\bar\np)\equiv 0.$ Thus from $(b')$ we get that, for each fixed $\bar\np_\star\in \bar V\cap\{\frac{1}{2}<F<\frac{2}{3}\}$,
\[
 \partial_\alpha\bar T_{00}(\bar\np_\star)+\partial_\alpha\bar T_{10}(\bar\np_\star)s+\partial_\alpha\bar T_{01}(\bar\np_\star)s^{\lambda(\bar\np_\star)}+\op(s^{2-2\upsilon})=0\text{ for all $s>0$,}
\] 
where we use that the flatness order of the remainder is preserved when derived with respect to parameters (see \defic{defi2}). Then, since $1<\lambda<2$ for $F\in (\frac{1}{2},\frac{2}{3})$ and we can choose $\upsilon>0$ arbitrary small (depending on $\bar\np_\star)$, by taking \obsc{rm1} into account we can assert that 
 \[
 \partial_\alpha\bar T_{00}(\bar\np_\star)= 0\text{, }\partial_\alpha\bar T_{10}(\bar\np_\star)= 0\text{ and }\partial_\alpha\bar T_{01}(\bar\np_\star)= 0.
\] 
Since $\bar V\cap\{\frac{1}{2}<F<\frac{2}{3}\}$ is open and the coefficients are meromorphic on $\bar V$, by \lemc{continuacio} it follows that $\partial_\alpha\bar T_{00}$, $\partial_\alpha\bar T_{10}$ and  $\partial_\alpha\bar T_{01}$ are identically zero. The claim for $\partial_\beta$ and the other coefficients follows verbatim. 

On account of the claim and the fact that $\bar T(s;\bar\np)|_{\eta=0}=T(s;\np)$, the assertions in $(a)$--$(d)$ with regard to the asymptotic expansion of $T(s;\np)$ follow from $(a')$--$(d')$, respectively, by setting
\[
 T_{ij}(\np)\!:=\left.\bar T_{ij}(\bar\np)\right|_{\eta=0}.
\]
Next we proceed with the computation of the expression of each coefficient. To this end observe that
\begin{equation}\label{t1eq1}
T_{ij}(\np)=\left.\bar T_{ij}(D,F,\alpha,\beta,\eta)\right|_{\eta=0}=\lim_{\eta\to 0^+}\bar T_{ij}(D,F,\alpha,\beta,\eta)=\lim\limits_{\eta\to 0^+}\bar T_{ij}(D,F,0,0,\eta),
\end{equation}
where the second follows by the continuity of $\bar\np\mapsto\bar T_{ij}(\bar\np)$ on $\bar V$ and the last one on account of the previous claim. In view of this the plan is to compute $\bar T_{ij}(D,F,0,0,\eta)$ with $\eta>0$ by applying \teoc{oldA} and then make $\eta\to 0.$ With this aim it is first necessary to obtain the functions in \refc{def_fun}. In doing so, and setting 
\[
 \text{$\delta_1\!=\frac{1}{2F},$ $\kappa_1\!=\frac{D+1}{F},$ $\delta_2\!=\frac{1}{2(1-F)}$ and $\kappa_2\!=\frac{D}{F-1},$}
\]
for shortness, one can verify that 
\[
 L_i(u)=h(u;\delta_i,\kappa_i)\text{ for $i=1,2$, where $h(u;\delta,\kappa)\!:=(1+\kappa u^2)^\delta$.}
\]
We stress that $L_1$ and $L_2$ are defined in \refc{def_fun} in terms of the functions $P_1(x_1,x_2;\bar\np)$ and $P_1(x_1,x_2;\bar\np)$ given in \refc{t1eq0} and that, as we already explained, we take $\alpha=\beta=0$ here and in what follows. Similarly, setting 
\[
 \gamma_1=-\frac{2D+1}{F^2}
\]
for the sake of shortness again, some computations show that $M_1(u)=\gamma_1 u h(u;\delta_1-2,\kappa_1),$
\begin{equation}\label{t1eq4}
 A_1(u)=-2\delta_1h(u;-1,\kappa_1),\;
 B_1(u)=M_1(u)
 \text{ and }
 \left.C_1(u)\right|_{D=-\frac{1}{2}}=-4\delta_1^2h(u;2\delta_1-2,\delta_1).
\end{equation}
Moreover $A_2(u)=2\delta_2h(u;\delta_2-1,\kappa_2).$ With regard to the parametrization of the transverse sections in the expression of the coefficients, see \refc{t1eq3}, using the compact notation $\sigma_{ijk}(\bar\np)\!:=\partial^k_s\sigma_{ij}(0;\bar\np)$ we get that 
\begin{equation}\label{t1eq2}
\sigma_{112}=\sigma_{121}=\sigma_{122}=0
\text{ and }
\sigma_{111}=\sigma_{120}=\sigma_{210}=\sigma_{221}=1/\eta.
\end{equation}
We are now in position to apply \teoc{oldA} to obtain the coefficients $\bar T_{ij}(\bar\np)$. (We omit the computation leading to $T_{00}(\np)$ because it is given in \cite[Proposition 5.2]{MMV03}.) In doing so we obtain that
\[
\bar T_{01}(\np,0,0,\eta)=\frac{\sigma_{120}\sigma_{111}^\lambda}{\sigma_{210}^\lambda L_1^\lambda(\sigma_{120})}\gorro A_2( \lambda,\sigma_{210})=2\delta_2\frac{\eta^\lambda}{(\eta^2+\kappa_1)^{\delta_2}}\hat h(\lambda,1/\eta;\delta_2-1,\kappa_2).
\]
Therefore
\begin{align*}
T_{01}(\np)
&=\lim_{\eta\to 0^+}\bar T_{01}(\np,0,0,\eta)
=
\frac{2\delta_2}{\kappa_1^{\delta_2}}\frac{\kappa_2^{\frac{\lambda}{2}}}{2}\mathrm{B}\left(-\frac{\lambda}{2},-\delta_2+1+\frac{\lambda}{2}\right)
\\[5pt]
&
=
\frac{\delta_2}{\sqrt{\kappa_2}}\left(\frac{\kappa_2}{\kappa_1}\right)^{\delta_2} \mathrm{B}\left(-\frac{\lambda}{2},\frac{1}{2}\right)
=
\frac{\sqrt{\pi}}{2(1-F)}\left(\frac{F}{D+1}\right)^{\frac{\lambda+1}{2}}\left(\frac{D}{F-1}\right)^{\frac{\lambda}{2}}\frac{\Gamma(-\frac{\lambda}{2})}{\Gamma(\frac{1-\lambda}{2})},
\end{align*}
where the first equality follows from \refc{t1eq1}, the second one by $(a)$ in \propc{B2F1} (provided that $\lambda\notin\N$), and the last one by using \refc{eqB0} and that $\Gamma(\frac{1}{2})=\sqrt{\pi}.$
Since $\rho_1(\np)=\frac{\sqrt{\pi}}{2(1-F)}\big(\frac{F}{D+1}\big)^{\frac{\lambda+1}{2}}\big(\frac{D}{F-1}\big)^{\frac{\lambda}{2}}$ is an analytic positive function on $V$, this proves the equality in the statement because we know that $T_{01}(\np)$ is meromorphic on $V$ with poles only at those $F\in (0,1)$ such that $\lambda(\np)\in\N.$ 

Let us turn next to the computation of $T_{10}$. With this aim we note that
\begin{align*}
 \bar T_{10}(\np,0,0,\eta)&=
 - \sigma_{120} \left(\frac{\sigma_{121}}{\sigma_{120}P_2(0,\sigma_{120})}+\frac{\sigma_{111}}{L_1(\sigma_{120})}\gorro{B}_1\big(1/\lambda-1,\sigma_{120}\big)\right)
 \\[5pt]
 &=
  -\frac{\eta^{-2}}{L_1(\eta^{-1})}\hat B_1(1/\lambda-1,\eta^{-1})=-\gamma_1\frac{\eta^{-3+2\delta_1}}{(\eta^2+\kappa_1)^{\delta_1}}\hat h(1/\lambda-2,\eta^{-1};\delta_1-2,\kappa_1).
\end{align*}
Here the first equality follows by \teoc{oldA}, the second one from \refc{t1eq2} and the last one by applying~$(d)$ in \teoc{L8} to the function $B_1(u)=\gamma_1 u h(u;\delta_1-2,\kappa_1)$, see \refc{t1eq4}. Since $-3+2\delta_1=\frac{1}{\lambda}-2$, by applying~$(a)$ in \propc{B2F1} we get
\[
 T_{10}(\np)=\lim_{\eta\to 0^+}T_{10}(\np,0,0,\eta)=-\frac{\gamma_1\kappa_1^{\frac{1}{2\lambda}-1}}{2\kappa_1^{\delta_1}}\mathrm{B}\left(1-\frac{1}{2\lambda},\frac{1}{2}\right)
 =\frac{\sqrt{\pi}(2D+1)}{2\sqrt{F(1+D)^3}}\frac{\Gamma(1-\frac{1}{2\lambda})}{\Gamma(\frac{3}{2}-\frac{1}{2\lambda})}
 \]
and, due to $\rho_2(\np)=\frac{\sqrt{\pi}}{2\sqrt{F(1+D)^3}}$, this proves the validity of the expression for $T_{10}$ given in the statement. Let us finally compute the coefficient $T_{20}.$ In this case, on account of \refc{t1eq2}, by \teoc{oldA} we get
\[
 \bar T_{20}(\np,0,0,\eta)=-\eta^{-1}\left(\frac{\eta^{-2}}{2L_1^2(\eta^{-1})}\hat C_1(2/\lambda-1,\eta^{-1})+\frac{\eta^{-1}S_1}{L_1(\eta^{-1})}\hat B_1(1/\lambda-1,\eta^{-1})\right).
\]
Thus, since $\gamma_1=0$ for $D=-\frac{1}{2},$ from \refc{t1eq4} it turns out that
\[
 \bar T_{20}(-1/2,F,0,0,\eta)=\frac{2\delta_1^2\eta^{-3+4\delta_1}}{(\eta^2+\kappa_1)^{2\delta_1}}\hat h(2/\lambda-1,\eta^{-1};2\delta_1-2,\delta_1)\Big|_{D=-\frac{1}{2}}
\]
Hence, due to $-3+4\delta_1=\frac{2}{\lambda}-1$, by applying $(a)$ in \propc{B2F1} once again, 
\begin{equation*}
 T_{20}(-1/2,F)=\lim_{\eta\to 0^+}\bar T_{20}(-1/2,F,0,0,\eta)=\frac{\kappa_1^{\frac{1}{\lambda}-\frac{1}{2}}}{4F^2\kappa_1^{2\delta_1}}\Big|_{D=-\frac{1}{2}}\mathrm{B}\left(\frac{1}{2}-\frac{1}{\lambda},\frac{1}{2}\right)=\frac{\sqrt{\pi}}{\sqrt{2F}}\frac{\Gamma(\frac{1}{2}-\frac{1}{\lambda})}{\Gamma(1-\frac{1}{\lambda})}.
\end{equation*}
Since $T_{20}(\np)$ is a meromorphic function having poles only at those $\np_0\in V$ such that $\lambda(\np_0)\in D_{20}=\frac{2}{\N}$ and, on the other hand, $\lambda(\np)>2$ for all $\np\in V\cap\{\frac{2}{3}<F<1\}$, by applying the Weierstrass Division Theorem (see for instance \cite[Theorem 1.8]{Greuel}), we can assert the existence of an analytic function $\rho_3$ on $V\cap\{\frac{2}{3}<F<1\}$ such that
\[
T_{20}(\np)=\frac{\sqrt{\pi}}{\sqrt{2F}}\frac{\Gamma(\frac{1}{2}-\frac{1}{\lambda})}{\Gamma(1-\frac{1}{\lambda})}+\rho_3(\np)(2D+1).
\]

It only remains to prove the assertions with regard to the properties of the coefficients in the respective asymptotic expansions. Being the ones in $(a)$ and $(b)$ an easy consequence of well-known properties of the gamma function (see for instance \cite{stegun}), we proceed with the other two:

\begin{enumerate}

\item[$(c)$] Let us take any $\np_0\in V\cap\{F=\frac{2}{3}\}$ and note that, from $(c')$, the functions $T_{201}^2(\np)\!:=\left.\bar T_{201}^2(\bar\np)\right|_{\eta=0}$ and $T_{200}^2(\np)\!:=\left.\bar T_{200}^2(\bar\np)\right|_{\eta=0}$ are smooth in a neighbourhood of $\{\np\in V:\lambda(\np)=2\}=V\cap\{F=\frac{2}{3}\}$ and
\[
 T_{201}^2(\np)=\big(2-\lambda(\np)\big)T_{01}(\np)
 \text{ and }
 T_{200}^2(\np)=T_{20}(\np)+T_{01}(\np)\text{ for $\lambda(\np)\neq 2.$}
\]
Recall that $T_{01}(\np)$ and $T_{20}(\np)$ are meromorphic with a pole 
at those $\np$ such that $\lambda(\np)=2\in D_{01}\cap D_{20}$.
What is more, by Propositions 3.2 and 3.6 in \cite{MV21}, respectively, we know that in both cases the pole is simple. Consequently by the Weierstrass Division Theorem (or, more directly, by \cite[Lemma 2.8]{MV21}) it follows that $T_{201}^2(\np)$ and $T_{200}^2(\np)$ are analytic in a neighbourhood of $V\cap\{F=\frac{2}{3}\}$. On the other hand, from the already proved part of the statement, if $\np_\star=(-\frac{1}{2},\frac{2}{3})$ then $T_{10}(\np_\star)=0$ and 
\[
 T_{201}^2(\np_\star)=\lim_{\np\to\np_\star}(2-\lambda)T_{01}(\np)=\frac{\rho_1(\np_\star)}{\Gamma(-\frac{1}{2})}\lim_{\np\to\np_\star}\Gamma(-\lambda/2)(2-\lambda)=-2\frac{\rho_1(\np_\star)}{\Gamma(-\frac{1}{2})}\neq 0
 \]
because $\lim_{x\to -1}(x+1)\Gamma(x)=-1$ (see \cite[\S 6]{stegun} for instance). 

\item[$(d)$] Consider finally any $\np_0\in V\cap\{F=\frac{1}{2}\}$. Then, from assertion
$(d')$, $T_{101}^1(\np)\!:=\left.\bar T_{101}^1(\bar\np)\right|_{\eta=0}$ and $T_{100}^1(\np)\!:=\left.\bar T_{100}^1(\bar\np)\right|_{\eta=0}$ are smooth functions in a neighbourhood of $\{\np\in V:\lambda(\np)=1\}=V\cap\{F=\frac{1}{2}\}$ and, in addition,
\[
 T_{101}^1(\np)=\big(1-\lambda(\np)\big)T_{01}(\np)
 \text{ and }
 T_{100}^1(\np)=T_{10}(\np)+T_{01}(\np)\text{ for $\lambda(\np)\neq1.$}
\]
Since $T_{10}(\np)$ and $T_{01}(\np)$ are meromorphic with a pole at those $\np$ such that $\lambda(\np)=1\in D_{10}\cap D_{01},$ the above equality shows exactly as before that $T_{101}^1(\np)$ and $T_{100}^1(\np)$ are analytic in a neighbourhood of $V\cap\{F=\frac{1}{2}\}$. Moreover, from the expression for $T_{01}$ in the statement that we already proved and using that $1-\lambda=2\frac{F-1/2}{F-1}$, we can write
\[
 T_{101}^1(\np)=(1-\lambda)T_{01}(\np)=-\rho_4(\np)(F-1/2)^2\text{ with $\rho_4(\np)\!:=\frac{\rho_1(\np)}{(F-1)^2}\frac{4\Gamma(-\frac{\lambda}{2})}{\Gamma(\frac{1-\lambda}{2})(\lambda-1)}.$}
\]
The function $\rho_4$ is analytic at $\nu_0=(D_0,\frac{1}{2})$ because $\lambda(\np_0)=1$ and 
$\Gamma(z)$ has simple pole at $z=0$. In addition, since $\lim_{z\to 0}z\Gamma(z)=1$ and $\Gamma(-\frac{1}{2})=-2\sqrt{\pi}$, we get that $\rho_4(\np_0)=16\sqrt{\pi}\rho_1(\nu_0)>0.$ From the expressions in the statement as well, we get 
\[
 T_{100}^1(\np)=\rho_5(\np)(D+1/2)+\rho_6(\np)(F-1/2)
\]  
with $\rho_5(\np)\!:=2\rho_2(\np)\frac{\Gamma(1-\frac{1}{2\lambda})}{\Gamma(\frac{3}{2}-\frac{1}{2\lambda})}$ and $\rho_6(\np)\!:=\frac{1-F}{2}\rho_4(\np)$, that are analytic and positive at $\nu_0=(D_0,\frac{1}{2})$ due to $\lambda(\np_0)=1$. Finally a computation shows that $\rho_5(-\frac{1}{2},\frac{1}{2})=\rho_6(-\frac{1}{2},\frac{1}{2}).$ 

\end{enumerate}
This concludes the proof of the result.
\end{prooftext}

\subsection{Proof of \propc{prop52}}

\begin{prooftext}{Proof of \propc{prop52}.} We will adapt the arguments in \cite[\S3.2.1]{MMV2} to take advantage of the general setting developed in~\cite{MV21}. To this end, as we did in the proof of the previous result, we will work in an extended parameter space $\bar W$ to be specified. In this case the computations are a little bit more involved because we also need to straighten the separatrices of the saddle, see \figc{fig4}. With this aim in view we first take $\varepsilon\in\R$ and consider the local change of coordinates given by 
\[
 (x_1,x_2)=\phi_{\varepsilon}(x,y)\!:=
 \left(\frac{q(x)-\frac{1}{2}y^2}{a(p_2-x+\varepsilon y)^2},\frac{p_2-p_1}{p_2-x+\varepsilon y}\right),
\]
where recall that $q(x)=a(x-p_1)(x-p_2)$ with $a=\frac{D}{2(1-F)}>0$ and $p_1,p_2\in\R,$ $p_1<p_2$, for all $\np\in W.$ In what follows, for the sake of shortness we set 
\[
 \kappa_1\!:=p_2-p_1\text{ and }\kappa_2\!:=1/\sqrt{2a}. 
\]
One can check that the Jacobian determinant of $\phi_\varepsilon$ vanishes at $(x,y)$ if and only if $y-\varepsilon q'(x)=0,$ where
\[
 q'(x)=2ax-a(p_1+p_2)x,
\]
and that this straight line is mapped by $\phi_\varepsilon$ to 
\[
 \mathscr D_{\varepsilon}(x_1,x_2)\!:=2a(1-x_1-x_2)+a^2\varepsilon^2(4x_1+x_2^2)=0.
\]
We claim that $\phi_\varepsilon$ is an analytic map from
\[
\Omega_{\varepsilon}\!:=\{(x,y)\in\R^2:p_2-x+\varepsilon y>0,\,y-\varepsilon q'(x)>0\}
\]
to $\mathscr U_{\varepsilon}\!:=\{(x_1,x_2)\in\R^2:x_2>0,\,\mathscr D_{\varepsilon}(x_1,x_2)>0\}$ with a well defined analytic inverse given by 
 \[
  \psi_{\varepsilon}(x_1,x_2)\!:=
  \left(\frac{\kappa_1(\varepsilon\mathscr D_{\varepsilon}(x_1,x_2)^{\frac{1}{2}}-1)+(p_2-\varepsilon^2a(p_1+p_2))x_2}{(1-2a\varepsilon^2)x_2},\frac{\kappa_1(\mathscr D_{\varepsilon}(x_1,x_2)^{\frac{1}{2}}+a\varepsilon(x_2-2))}{(1-2a\varepsilon^2)x_2}\right).
 \]
Indeed, the claim follows by checking that $\phi_\varepsilon\circ\psi_\varepsilon=\text{Id}$ on $\{\mathscr D_{\varepsilon}(x_1,x_2)>0,x_2\neq 0\}$ and that $\psi_\varepsilon\circ\phi_\varepsilon=\text{Id}$ on $\{\frac{y-\varepsilon q'(x)}{p_2-x+\varepsilon y}>0\}.$ To show the second identity we use that $(\mathscr D_\varepsilon\circ\phi_\varepsilon)(x,y)=\left(\frac{y-\varepsilon q'(x)}{p_2-x+\varepsilon y}\right)^2.$ Consequently $(x_1,x_2)=\phi_\varepsilon(x,y)$ is an analytic global change of variables from $\Omega_\varepsilon$ to $\mathscr U_\varepsilon$ for all $\varepsilon.$ In what follows we require $|\varepsilon|<\frac{1}{\sqrt{2a}}$ in order that the straight line 
$\{\varepsilon y+p_2-x=0\}$ does not intersect the left branch of the hyperbola $\mathcal C=\{\frac{1}{2}y^2-q(x)=0\}$, see \figc{fig5}.
\begin{figure}[t]
 \centering
 \begin{lpic}[l(0mm),r(0mm),t(0mm),b(0mm)]{dib6c(1)}    
  \lbl[l]{0,33;$\Sigma_2^\eta$}  
  \lbl[l]{65,14;$\Sigma_1^\eta$} 
  \lbl[l]{67,5;$\frac{p_1+p_2}{2}$}  
  \lbl[l]{79,5;$p_2$}
  \lbl[l]{60,27;$\mathcal C$}
  \lbl[l]{53,52;$\{\varepsilon y+p_2-x=0\}$}
  \lbl[l]{20,29;$\ell_\eta$}
  \lbl[l]{10,18;$\{y-\varepsilon q'(x)=0\}$}  
  \lbl[l]{88,8;$x$}
  \lbl[l]{49,57;$y$}                   
 \end{lpic}
 \caption{Auxiliary transverse sections in the proof of \propc{prop52}.}\label{fig5}
 \end{figure}

Next we introduce a second auxiliary parameter $\eta\in\R$ and consider two additional transverse sections~$\Sigma_1^\eta$ and $\Sigma_2^\eta$ laying on the straight line $\ell_\eta\!:=\{y=\eta(p_2-x)\},$ see \figc{fig5}. 
Observe that~$\ell_\eta$ intersects the right branch of the hyperbola $\mathcal C$ at $(p_2,0)$ for all $\eta$. We require $|\eta|<\sqrt{2a}$ additionally in order that~$\ell_\eta$ intersects the hyperbola at a point $(x_\eta,y_\eta)$ in the left branch. Then we parametrize $\Sigma_1^\eta$ and $\Sigma_2^\eta$, respectively, by 
 \begin{equation}\label{p2eq1}
  s\mapsto \big(x_\eta-s,\eta(p_2-x_\eta+s)\big)
  \text{ and }
  s\mapsto \big(-1/s,\eta(p_2+1/s)\big), 
 \end{equation}
for $s>0$ small enough. We also require $\varepsilon(\eta^2+2a)+2\eta>0$ so that $(x_\eta,y_\eta)\in\Omega_\varepsilon.$ Summing up, the admissible conditions 
\[
 |\varepsilon|<\frac{1}{\sqrt{2a}},\;|\eta|<\sqrt{2a}\text{ and }\varepsilon(\eta^2+2a)+2\eta>0
\]
guarantee that any solution of $X_\np$ going from $\Sigma_1^\eta$ to $\Sigma_2^\eta$ is inside the domain $\Omega_\varepsilon$ of the coordinate change $(x_1,x_2)=\phi_\varepsilon(x,y)$. Thus, setting $\bar\np\!:=(\np,\varepsilon,\eta)$, we will work on the extended parameter space
\[\textstyle
 \bar W\!:=\big\{\bar\np\in\R^4: F+D>0,D<0,F>1,|\varepsilon|<\frac{1}{\sqrt{2a}},\;|\eta|<\sqrt{2a}\text{ and }\varepsilon(\eta^2+2a)+2\eta>0\big\}.
\]
Clearly the sets $\{(\np,\varepsilon,0): \np\in W\text{ and }\varepsilon\in(0,\frac{1}{\sqrt{2a}})\}$ and $\{(\np,0,\eta): \np\in W\text{ and }\eta\in (0,\sqrt{2a})\}$ are inside~$\bar W$, that will be crucial in the forthcoming steps. 

At this point we define $\sigma_1(\,\cdot\,;\bar\np)$ and $\sigma_2(\,\cdot\,;\bar\np)$ to be, respectively, the composition with $\phi_\varepsilon$ of the parametrization of~$\Sigma_1^\eta$ and~$\Sigma_2^\eta$ given in \refc{p2eq1}. In its regard one can check that
\begin{align}\label{p2eq2}
 \left.\sigma_1(s;\bar\np)\right|_{\varepsilon=0}&=\left(\frac{s(1-\kappa_2^2\eta^2)^2}{\kappa_1+s(1-\kappa_2^2\eta^2)},\frac{\kappa_1(1-\kappa_2^2\eta^2)}{\kappa_1+s(1-\kappa_2^2\eta^2)}\right)\\
\intertext{and}\label{p2eq3}
\left.\sigma_2(s;\bar\np)\right|_{\varepsilon=0}&=\left(\frac{1+p_1s}{1+p_2s}-\kappa_2^2\eta^2,\frac{\kappa_1s}{1+p_2s}\right).
\end{align}
One can also verify that the coordinate change $(x_1,x_2)=\phi_{\varepsilon}(x,y)$ brings the vector field $X_\np$ in \refc{sist_loud} to 
\[
 \bar X_{\bar\np}(x_1,x_2)=\frac{1}{x_2}\big(x_1P_1(x_1,x_2;\bar\np)\partial_{x_1}
      +x_2P_2(x_1,x_2;\bar\np)\partial_{x_2}\big),
\]
where $P_1$ and $P_2$ analytic functions on $\{(x_1,x_2,\bar\np)\in\R^2\times\bar W: \mathscr D_\varepsilon(x_1,x_2)>0\}.$ The hyperbolicity ratio of the saddle at the origin is
\[
\lambda=-\frac{P_2(0,0)}{P_1(0,0)}=\frac{1}{2(F-1)}.
\]
Moreover $P_1|_{\varepsilon=0}=R\bar P_1$ and $P_2|_{\varepsilon=0}=R\bar P_2$ where $R(x_1,x_2)=\frac{1}{\kappa_2}\sqrt{1-x_1-x_2}$,
\begin{equation}\label{p2eq8}
 \bar P_1(x_1,x_2)=2\kappa_1(F-1)+2(p_2-1)x_2\text{ and }
 \bar P_2(x_1,x_2)=-\kappa_1+(p_2-1)x_2.
\end{equation}
(It will be clear in a moment the reason why it suffices to restrict to $\varepsilon=0.$) For each $\bar\np\in\bar W,$ we define $\bar T(s;\bar\np)$ to be the Dulac time of $\bar X_{\bar\np}$ between the transverse sections $\phi_{\varepsilon}(\Sigma^\eta_1)$ and $\phi_{\varepsilon}(\Sigma^\eta_2)$ parametrized by $\sigma_1(\,\cdot\,;\bar\np)$ and $\sigma_2(\,\cdot\,;\bar\np)$, respectively. We point out that, by construction, $\bar T(s;\bar\np)$ does not depend on $\varepsilon$ and that, furthermore, $\left.\bar T(s;\bar\np)\right|_{\eta=0}=T(s;\np).$ 

Next we will apply \teoc{5punts} to obtain the asymptotic expansion of $\bar T(s;\bar\np)$ at $s=0$. Note to this end that, by construction, given any $\bar\np_0\in\bar W$ there exists a relatively compact neighbourhood~$\mathcal V_0$ of 
\[
 \big\{(x_1,0): x_1\in [0,\bar\sigma_{21}(0;\bar\np_0)]\big\}
 \cup 
 \big\{(0,x_2): x_2\in [0,\bar\sigma_{12}(0;\bar\np_0)]\big\}
\]
in $\R^2$ and a neighbourhood $\bar W_0$ of $\bar\np_0$ in $\bar W$ such that $\phi_\varepsilon\big(\Sigma_1^\eta\cup\Sigma_2^\eta\big)\subset\mathcal V_0$ for all $\bar\np\in\bar W_0$ and 
\[
 \mathcal V_0\times\bar W_0\subset\{(x_1,x_2,\bar\np)\in\R^2\times\bar W: \mathscr D_\varepsilon(x_1,x_2)>0\}.
\] 
Here we use (see also \figc{fig5}) that $\phi_\varepsilon$ maps the straight line $\{y-\varepsilon q'(x)=0\}$ to $\{\mathscr D_\varepsilon(x_1,x_2)=0\}$. The above inclusion guarantees that $P_1$ and $P_2$ are analytic on $\mathcal V_0\times\bar W_0,$ so that we can apply \teoc{5punts} to study the Dulac time of $\bar X_{\bar\np}$ for $\bar\np\approx\bar\np_0$. Accordingly, with this aim, let us fix any $\bar\np_0=(D_0,F_0,\varepsilon_0,\eta_0)\in \bar W$ with $F_0\in (1,\frac{3}{2})\cup\{2\}.$ Observe that 
$\lambda(\bar\np_0)>2$ if $F_0\in (1,\frac{5}{4})$, $\lambda(\bar\np_0)\in (1,2)$ if $F_0\in (\frac{5}{4},\frac{3}{2})$, $\lambda(\bar\np_0)=2$ if $F_0=\frac{5}{4}$ and $\lambda(\bar\np_0)=\frac{1}{2}$ if $F_0=2.$ Then, by applying $(2),$ $(1)$, $(5)$ and $(3)$ in \teoc{5punts}, respectively, we can assert that
\begin{enumerate}[$(a')$]

\item $\bar T(s;\bar\np)=\bar T_{00}(\bar\np)+\bar T_{10}(\bar\np)s+\bar T_{20}(\bar\np)s^2+\F_{L_0-\upsilon}^\infty(\bar\np_0)$ if $F_0\in (1,\frac{5}{4})$, where $L_0=\min(3,\lambda(\bar\np_0)),$ 

\item $\bar T(s;\bar\np)=\bar T_{00}(\bar\np)+\bar T_{10}(\bar\np)s+\bar T_{01}(\bar\np)s^\lambda+\bar T_{20}(\bar\np)s^2+\F^\infty_{L_0-\upsilon}(\bar\np_0)$ if $F_0\in (\frac{5}{4},\frac{3}{2})$, where $L_0=\lambda(\bar\np_0)+1,$ 

\item $\bar T(s;\bar\np)=\bar T_{00}(\bar\np)+\bar T_{10}(\bar\np)s+\bar T_{201}^2(\bar\np)s^2\omega_{2-\lambda}(s)+\bar T_{200}^2(\bar\np)s^2+\F_{3-\upsilon}^\infty(\bar\np_0)$ if $F_0=\frac{5}{4}$, where
$\bar T_{201}^2(\bar\np)$ and $\bar T_{200}^2(\bar\np)$ are smooth in a neighbourhood of 
$\{\bar\np\in\bar W:\lambda(\bar\np)=2\}$ and, moreover,
\[
 \bar T_{201}^2(\bar\np)=(2-\lambda)\bar T_{01}(\bar\np)
 \text{ and }
 \bar T_{200}^2(\bar\np)=\bar T_{20}(\bar\np)+\bar T_{01}(\bar\np)
 \text{ for $\lambda(\bar\np)\neq 2,$}
\]

\item $\bar T(s;\bar \np)=\bar T_{00}(\bar \np)+\bar T_{01}(\bar \np)s^\lambda
        +\bar T^{\frac{1}{2}}_{101}(\bar \np)s\omega_{1-2\lambda}(s)
        +\bar T_{100}^{\frac{1}{2}}(\bar \np)s+\F_{3/2-\upsilon}^\infty(\bar \np_0)$ if $F_0=2$, 
         where $\bar T_{101}^1(\bar\np)$ and $\bar T_{100}^1(\bar\np)$ are smooth in a neighbourhood of $\{\bar\np\in\bar W:\lambda(\bar\np)=1/2\}$ 
         and, moreover,
         \[
        \bar T_{101}^{\frac{1}{2}}(\bar\np)=(1-2\lambda)\bar T_{02}(\bar\np)
        \text{ and }
       \bar T_{100}^{\frac{1}{2}}(\bar\np)=\bar T_{10}(\bar\np)+\bar T_{02}(\bar\np)
       \text{ for $\lambda(\bar\np)\neq 1/2.$}
\]

\end{enumerate}
Here $\upsilon$ is a small enough positive number depending on $\bar\np_0.$ Furthermore by applying locally \teoc{5punts} we know that the coefficients $\bar T_{ij}(\bar\np)$ are meromorphic functions on $\bar W$ with poles only at those $F>1$ such that $\lambda(\bar\np)=\frac{1}{2(F-1)}\in D_{ij},$ where $D_{00}=\emptyset$, $D_{10}=\frac{1}{\N}$, $D_{01}=\N$, $D_{20}=\frac{2}{\N}$ and $D_{02}=\frac{\N}{2}$.

We claim that the coefficients $\bar T_{ij}(\bar\np)$ do not depend on $\varepsilon.$ To prove this observe that the Dulac time $\bar T(s;\bar\np)$ does not depend on $\varepsilon$ as long as $\bar\np\in\bar W.$ Accordingly $\partial_\varepsilon\bar T(s;\bar\np)\equiv 0$. Thus from $(b')$ we get that, for each fixed $\bar\np_\star\in \bar W\cap\{\frac{5}{4}<F<\frac{3}{2}\}$,
\[
 \partial_\varepsilon\bar T_{00}(\bar\np_\star)+\partial_\varepsilon\bar T_{10}(\bar\np_\star)s+\partial_\varepsilon\bar T_{01}(\bar\np_\star)s^{\lambda(\bar\np_\star)}+\partial_\varepsilon\bar T_{20}(\bar\np_\star)s^2+\op(s^{L_0-2\upsilon})=0\text{ for all $s>0$,}
\] 
where $L_\star=\lambda(\bar\np_\star)+1$ and we use that the flatness order of the remainder is preserved when derived with respect to parameters (see \defic{defi2}). Then, since $1<\lambda<2$ for $F\in (\frac{5}{4},\frac{3}{2})$ and we can choose $\upsilon>0$ arbitrary small (depending on $\bar\np_\star)$, the application of \obsc{rm1} shows that 
 \[
 \partial_\varepsilon\bar T_{00}(\bar\np_\star)= 0\text{, }\partial_\varepsilon\bar T_{10}(\bar\np_\star)= 0\text{, }\partial_\varepsilon\bar T_{01}(\bar\np_\star)= 0\text{ and }\partial_\varepsilon\bar T_{20}(\bar\np_\star)= 0.
\] 
Since $\bar W\cap\{\frac{5}{4}<F<\frac{3}{2}\}$ is open and the coefficients are meromorphic on $\bar W$, by \lemc{continuacio} it follows that $\partial_\varepsilon\bar T_{00}$, $\partial_\varepsilon\bar T_{10}$, $\partial_\varepsilon\bar T_{01}$ and  $\partial_\varepsilon\bar T_{20}$ are identically zero. The claim for $\partial_\varepsilon\bar T_{02}$ follows verbatim. 

Thanks to the claim and the fact that $\left.\bar T(s;\bar\np)\right|_{\eta=0}=T(s;\np)$ by construction, the assertions in $(a)$--$(d)$ concerning the asymptotic expansion of $T(s;\np)$ at $s=0$ follow from $(a')$--$(d')$, respectively, by setting 
\[
 T_{ij}(\np)\!:=\left.\bar T_{ij}(\bar\np)\right|_{\eta=0}.
\]
We proceed next with the computation of these coefficients and for this purpose the idea is that if $\np\in W$ and $\varepsilon>0$ then 
 \begin{equation}\label{eq521}
T_{ij}(\np)=\bar T_{ij}(\np,\varepsilon,0)=\lim_{\eta\to 0^+}\bar T_{ij}(\np,\varepsilon,\eta)=\lim_{\eta\to 0^+}\bar T_{ij}(\np,0,\eta),
 \end{equation}
where in the second equality we use the continuity of $\bar T_{ij}(\bar\np)$ at any $\bar\np_0\in\bar W$ with $\lambda(\bar\np_0)=\frac{1}{2(F_0-1)}\notin D_{ij}$ and in the last one the fact that $\bar T_{ij}(\bar\np)$ does not depend on $\varepsilon$. Hence our first goal is to obtain $\bar T_{ij}(\bar\np)$ for $\varepsilon=0$ and to this end we shall apply \teoc{oldA}. (We point out that from now on all the computations are performed taking $\varepsilon=0$ and $\eta>0$.) In doing so, and setting 
\[
\kappa_0\!:=\frac{p_2-1}{p_2-p_1}
\]
for shortness, from \refc{def_fun} we obtain that $L_1(u)=(1-\kappa_0 u)^{2F}$, $L_2(u)\equiv 1$, $M_1(u)\equiv 0$ and 
\begin{equation}\label{p2eq4}
\begin{array}{ll}
  A_1(u)=-\frac{\kappa_2}{\kappa_1}(1-\kappa_0 u)^{-1}(1-u)^{-\frac{1}{2}} 
 & 
  A_2(u)=\frac{\kappa_2}{2\kappa_1(F-1)}(1-u)^{-\frac{1}{2}}
 \\[15pt]
  B_1(u)=-\frac{\kappa_2}{2\kappa_1}(1-\kappa_0 u)^{2F-1}(1-u)^{-\frac{3}{2}}
 &
  C_1(u)=-\frac{3\kappa_2}{4\kappa_1}(1-\kappa_0 u)^{4F-1}(1-u)^{-\frac{5}{2}}.
 \\
\end{array}
\end{equation}
From \refc{p2eq2} and \refc{p2eq3}, the necessary information with regard to the transverse sections is the following:
\begin{equation}\label{p2eq5}
\textstyle
\sigma_{120}=\sigma_{210}=1-\kappa_2^2\eta^2,\,\; 
 \sigma_{111}=-\sigma_{121}=\frac{(1-\kappa_2^2\eta^2)^2}{\kappa_1} \text{ and }
 \sigma_{122}=-\sigma_{112}=\frac{2(1-\kappa_2^2\eta^2)^3}{\kappa_1^2}.
\end{equation}
Taking this into account we obtain that
\begin{align*}
 T_{00}(\np)
 &
 =\lim\limits_{\eta\to 0^+}\bar T_{00}(\np,0,\eta)
 =-\lim\limits_{\eta\to 0^+}(1-\kappa_2^2\eta^2)\hat A_1(-1,1-\kappa_2^2\eta^2) \\
 &\textstyle
 =\frac{\kappa_2}{\kappa_1}\,\mathrm{B}(1,\frac{1}{2})\,{}_2F_1(1,1;\frac{3}{2};\kappa_0)
 =\frac{2\kappa_2}{\kappa_1(1-\kappa_0)}\,{}_2F_1(1,\frac{1}{2};\frac{3}{2};\frac{\kappa_0}{\kappa_0-1})
\end{align*}
where in the second equality we use \teoc{oldA}, in the third one we apply $(b)$ in \propc{B2F1} taking 
$\{\alpha=-1,\ga=1,\delta=\frac{1}{2},x=\kappa_0\}$ and in the last one \cite[\S 15.3]{stegun}. Since $\frac{\kappa_0}{\kappa_0-1}=\frac{1-p_2}{1-p_1}$, this shows the validity of the expression for $T_{00}$ given in the statement. 
Similarly, from \refc{p2eq4} and \refc{p2eq5} again,  
\begin{align*}\textstyle
T_{01}(\np)
 &
 =\lim\limits_{\eta\to 0^+}\bar T_{01}(\np,0,\eta) \\
 &\textstyle
 =\frac{\kappa_1^{-\lambda}}{(1-\kappa_0)^{2\lambda F}}\lim\limits_{\eta\to 0^+}\hat A_2(\lambda,1-k^2_2\eta^2)
 =\frac{\kappa_2(1-\kappa_0)^{-2\lambda F}}{2\kappa_1^{\lambda+1}(F-1)}
  {\mathrm{B}(-\lambda,\frac{1}{2})}.
\end{align*}
Here the last equality follows by applying $(b)$ in \propc{B2F1} taking $\{\alpha=\lambda,\delta=\frac{1}{2},x=0\}$, so that $\alpha=\lambda(\np)=\frac{1}{2(F-1)}\notin\Z_{\geq 0}$, and noting that ${}_2F_1(a,b;c;0)=1$, see \refc{hyper}. Since $\rho_1(\np)\!:=\frac{\kappa_2(1-\kappa_0)^{-2\lambda F}}{2\kappa_1^{\lambda+1}(F-1)}$ is an analytic positive  function on $W$, this proves the the expression for $T_{01}$ given in the statement. 

Let us study next the coefficient $T_{10}$. For this purpose we apply \teoc{oldA}, which on account of~\refc{p2eq5} shows that if $\eta>0$ and $\lambda(\np)\notin D_{10}=\frac{1}{\N}$ then
 \begin{equation}\label{p2eq6}
  \bar T_{10}(\np,0,\eta)=\frac{1}{\kappa_1}\left(\frac{-1+\kappa_2^2\eta^2}{\eta\,(\kappa_1+(p_2-1)(1-\kappa_2^2\eta^2))}-
   \frac{(1-\kappa_2^2\eta^2)^3}{(1-\kappa_0(1-\kappa_2^2\eta^2))^{2F}}\,\hat B_1\!\left(1/\lambda-1,1-\kappa_2^2\eta^2\right)\right).
 \end{equation}
Following the notation in \propc{B2F1}, see \refc{p2eq4}, we can write $B_1(1-\kappa_2^2\eta^2)=-\frac{\kappa_2}{2\kappa_1}g(y;\delta,\ga;x)$ with $\{y=1-\kappa_2^2\eta^2,\delta=-\frac{1}{2},\ga=2F,x=\kappa_0\}$ but we cannot apply it to get the limit of $\hat B_1\!\left(1/\lambda-1,1-\kappa_2^2\eta^2\right)$ as $\eta\to 0^+$ because the condition $\delta>0$ is not satisfied. As a matter of fact this is coherent because, since the first summand in \refc{p2eq6} is divergent as $\eta\to 0^+,$ it happens that $\lim_{\eta\to 0^+}\hat B_1\!\left(1/\lambda-1,1-\kappa_2^2\eta^2\right)$ diverges as well. Hence the approach to compute this coefficient hast to be different. The idea is to take advantage of \cite[Theorem~3.6]{MMV2}, which shows that if $F\in (1,\frac{3}{2})$ then $T_{10}(\np)=\rho_2(\np)\bar\rho_2(\np)$ where
\[
 \rho_2(\np)\!:=\frac{\kappa_2}{2\kappa_1(1-p_1)}
 \text{ and }
 \bar\rho_2(\np)\!:=-2+
 \int_0^1\left((1-s)^{-\frac{1}{\lambda}}(1-\bar\kappa s)^{1+\frac{1}{\lambda}}-1\right)s^{-\frac{1}{2}}\frac{ds}{s}
\]
with $\bar\kappa\!:=\frac{1-p_2}{1-p_1}.$ By applying assertion $(b)$ in \teoc{L8} taking $f(s;\np)=(1-s)^{-\frac{1}{\lambda}}(1-\bar\kappa s)^{1+\frac{1}{\lambda}}$, $k=1$ and $\alpha=\frac{1}{2}$ we can assert that $\bar\rho_2(\np)=\lim_{s\to 1^-}\hat f(\frac{1}{2},s;\np)$. Observe on the other hand that, following the notation in \propc{B2F1}, we can write $f(s;\np)=g(y;\delta,\ga;x)$ with $\{y=s,\delta=1-\frac{1}{\lambda},\ga=-1-\frac{1}{\lambda},x=\bar\kappa\}.$ Thus, since one can verify that $\delta=1-\frac{1}{\lambda}>0$ and $x=\bar\kappa<1$ for all $\np\in W\cap\{1<F<\frac{3}{2}\},$ the application of assertion $(b)$ in that result gives
\[\textstyle
 \bar\rho_2(\np)=\lim_{y\to 1^-}\hat f(\frac{1}{2},y;\np)=\mathrm B\left(-\frac{1}{2},1-\frac{1}{\lambda}\right){}_2F_1\left(-1-\frac{1}{\lambda},-\frac{1}{2};\frac{1}{2}-\frac{1}{\lambda};\bar\kappa\right).
\]
Due to $\mathrm B(-\frac{1}{2},1-\frac{1}{\lambda})=\mathrm B(1-\frac{1}{\lambda},-\frac{1}{2})$, see \refc{eqB0},
this proves the validity of the expression for $T_{10}(\np)$ in the statement for all $\np\in W\cap\{1<F<\frac{3}{2}\}.$ Accordingly, since $T_{10}$ is meromorphic on $W$, this is also the case of~$\bar\rho_2$ thanks to \lemc{2F1}, and $\rho_2$ is analytic on $W$, the application of the real version of \lemc{continuacio} implies the validity of the equality on $W$. Observe moreover that $\rho_2$ is positive on $W.$ 

We proceed with the computation of the coefficient $T_{20}$. In first instance, for the sake of convenience we shall work with $\np\in W\cap\{1<F<\frac{5}{4}\},$ so that $\lambda(\np)>2.$ Due to $M_1\equiv0,$ \teoc{oldA} shows that if $\lambda(\np)\notin D_{20}=\frac{2}{\N}$ then
\begin{align}\label{p2eq7}
\bar T_{20}(\np,0,\eta)
&=- \frac{\sigma_{120}\sigma_{122}}{2\sigma_{120}P_2(0,\sigma_{120})}-\frac{1}{2}\sigma_{121}^2\partial_2P_2^{-1}(0,\sigma_{120})-\sigma_{121}\sigma_{111}\partial_1P_2^{-1}(0,\sigma_{120})
\\[4pt]
&\hspace{-1truecm}
-\frac{\sigma_{120}\sigma_{111}^2}{2L_1^2(\sigma_{120})}\gorro{C}_1(2/\lambda-1,\sigma_{120})-\left(\frac{\sigma_{112}}{2\sigma_{111}}-\frac{\sigma_{121}}{\sigma_{120}}\left(\frac{P_1}{P_2}\right)\!(0,\sigma_{120})\right)
\frac{\sigma_{120}\sigma_{111}}{L_1(\sigma_{120})}\gorro{B}_1(1/\lambda-1,\sigma_{120}).
\notag
\end{align}
On account of $\left.P_2\right|_{\varepsilon=0}=R\bar P_2$, $R(0,\sigma_{120})=\frac{1}{\kappa_2}\left.\sqrt{1-x_2}\right|_{x_2=1-\kappa_2^2\eta^2}=\eta$ and $\bar P_2(0,1)=p_1-1\neq 0$, see~\refc{p2eq8} and~\refc{p2eq5}, it follows that we can write
 \begin{align}\label{eq525}
 &
 - \frac{\sigma_{120}\sigma_{122}}{2\sigma_{120}P_2(0,\sigma_{120})}-\frac{1}{2}\sigma_{121}^2\partial_2P_2^{-1}(0,\sigma_{120})-\sigma_{121}\sigma_{111}\partial_1P_2^{-1}(0,\sigma_{120})=\frac{\phi_{1}(\eta^2)}{\eta^{3}}\\
 \intertext{and}\label{eq526}
 &
 -\frac{\sigma_{112}}{2\sigma_{111}}+\frac{\sigma_{121}}{\sigma_{120}}\left(\frac{P_1}{P_2}\right)\!(0,\sigma_{120})=\phi_{2}(\eta^2),
 \end{align}
where here (and in what follows) $\phi_i(x)$ stands for an analytic function at $x=0.$ (In fact $\phi_i$ depends also on $\np$ and this dependence is analytic on $W$. We omit this dependence for brevity when there is no risk of confusion.) Following this notation, from \refc{p2eq6} and using that $\lim_{\eta\to 0^+}L_1(1-\kappa_2^2\eta^2)=(1-\kappa_0)^{2F}\neq 0,$ we can assert that
\begin{equation}\label{eq523}
\frac{\sigma_{120}\sigma_{111}}{L_1(\sigma_{120})}\gorro{B}_1(1/\lambda-1,\sigma_{120})
=
\phi_{3}(\eta^2)\bar T_{10}(\np,0,\eta)+\frac{\phi_{4}(\eta^2)}{\eta}
\end{equation}
as long as $\lambda(\np)\notin D_{10}=\frac{1}{\N}.$ On the other hand, since $\lambda(\np)>2$ for all $\np\in W\cap\{1<F<\frac{5}{4}\},$ by applying assertion~$(b)$ in \teoc{L8} with $k=0$ and $\alpha=\frac{2}{\lambda}-1$, we get
\begin{align*}
\sigma_{120}^{1-\frac{2}{\lambda}}\hat C_1{(2/\lambda-1},\sigma_{120})&=\int_0^{\sigma_{120}}C_1(u)u^{1-\frac{2}{\lambda}}\frac{du}{u}=\int_{0}^{1-\kappa_2^2\eta^2}C_1(u)u^{-\frac{2}{\lambda}}du\\
 &=-\frac{3\kappa_2}{4\kappa_1}(1-\kappa_0)^{4F-1}\underbrace{\int_{\kappa_2^2\eta^2}^1(1-\bar\kappa x)^{4F-1}x^{-\frac{5}{2}}(1-x)^{-\frac{2}{\lambda}}dx}_{I},
\end{align*}
where in the second equality we perform the change of variable $u=1-x$ and use that $\frac{\kappa_0}{\kappa_0-1}=\frac{1-p_2}{1-p_1}=\bar\kappa$. Next we split the above integral as $I=I_1+I_2$ with
\[
 I_1\!:=\int_{\bar\kappa_2^2\eta^2}^1\left((1-\bar\kappa x)^{4F-1}(1-x)^{-\frac{2}{\lambda}}-(1+\beta x)\right)x^{-\frac{3}{2}}\frac{dx}{x}
 \text{ and }
 I_2\!:=\int_{\bar\kappa_2^2\eta^2}^1(1+\beta x)x^{-\frac{3}{2}}\frac{dx}{x},
 \]
where we take $\beta\!:=\frac{2}{\lambda}-(4F-1)\bar\kappa$ so that $I_1$ converges as $\eta\to 0$. Due to 
$I_2=-\frac{2}{3}-2\beta+\frac{2}{3}\frac{1}{(\kappa_2\eta)^3}+\frac{2\beta}{\kappa_2\eta},$ we can write $I=J+\frac{2}{3}\frac{1}{(\kappa_2\eta)^3}+\frac{2\beta}{\kappa_2\eta}$ with
\[
 J(\np,\eta)\!:=
 \int_{\kappa_2^2\eta^2}^1\left((1-\bar\kappa x)^{4F-1}(1-x)^{-\frac{2}{\lambda}}-(1+\beta x)\right)x^{-\frac{3}{2}}\frac{dx}{x}
 -\frac{2}{3}-2\beta.
\]
Consequently we obtain that
\begin{equation}\label{eq522}
 \hat C_1{(2/\lambda-1},\sigma_{120})=-\frac{3\kappa_2}{4\kappa_1}\frac{(1-\kappa_0)^{4F-1}}{(1-\kappa_2^2\eta^2)^{1-\frac{2}{\lambda}}}\left(J(\np,\eta)+\frac{2}{3(\kappa_2\eta)^3}+\frac{2\beta}{\kappa_2\eta}\right).
 \end{equation}
On the other hand, setting $g(x;\np)\!:=(1-\bar\kappa x)^{4F-1}(1-x)^{-\frac{2}{\lambda}}$ we get that
\[
 \lim_{\eta\to 0}J(\np,\eta)=\lim_{x\to 1^-}\hat g(3/2,x;\np)
 =
 \mathrm B\left(1-\frac{2}{\lambda},-\frac{3}{2}\right){}_2F_1\left(1-4F,-\frac{3}{2};-\frac{1}{2}-\frac{2}{\lambda};\bar\kappa\right)=:\!J_0(\np).
\] 
Here the first equality follows by $(b)$ in \teoc{L8} taking $\{\alpha=\frac{3}{2},k=2\}$ and the second one by~$(b)$ in \propc{B2F1} taking $\{\alpha=\frac{3}{2}, \ga=1-4F,\delta=1-\frac{2}{\lambda},x=\bar\kappa\}$ and using that $\delta>0$ thanks to $\lambda(\np)>2$ for all $\np\in W\cap\{1<F<\frac{5}{4}\}.$ That being said, substituting \refc{eq525}, \refc{eq526}, \refc{eq523} and \refc{eq522} into \refc{p2eq7} and gathering next the analytic functions at $\eta=0$ we obtain that
\[
\bar T_{20}(\np,0,\eta)=\phi_5(\np,\eta^2)J(\np,\eta)+\phi_6(\np,\eta^2)\bar T_{10}(\np,0,\eta)+\phi_7(\np,\eta^2)/\eta^3,
\] 
with $\phi_i(\np,x)$ analytic at $x=0$. (Here we specify again the dependence on~$\np$ for the sake of consistency in the exposition.) Consequently, from \refc{eq521}, 
 \[
  T_{20}(\np)=\lim_{\eta\to 0^+}\bar T_{20}(\np,0,\eta)=\phi_5(\np,0)J_0(\np)
         +\phi_6(\np,0)T_{10}(\np),
 \]
where we use that $\phi_7(\np,\eta^2)=\rm O(\eta^4)$ because the limit must be finite. One can easily check that 
\[
\rho_3(\np)\!:=\phi_5(\np,0)=\frac{3\kappa_2}{8\kappa_1^3(1-\kappa_0)}
\text{ and }
\rho_4(\np)\!:=\phi_6(\np,0)=\frac{p_1-1+2F\kappa_1}{\kappa_1(p_1-1)}.
\]
This proves the validity of the expression for $T_{20}(\np)$ for all $\np\in W\cap\{1<F<\frac{5}{4}\}.$ Similarly as before, by applying \lemc{continuacio} this equality extends to $W$ since $\rho_3$ and $\rho_4$ are analytic on $W$ and, on the other hand, $J_0$ is meromorphic on $W$ by \lemc{2F1} and so are $\bar T_{10}$ and $\bar T_{20}$. Finally an easy computation shows that $\rho_3$ and $\rho_4$ are positive on $W.$

So far we have proved the validity of the expression of the coefficients $T_{ij}(\np)$ that we give in the first part of the statement. Moreover, since $T(s;\np)=\bar T(s;\bar\np)|_{\eta=0}$ and $T_{ij}(\np)=\bar T_{ij}(\bar\np)|_{\eta=0},$ the assertions with regard to the asymptotic expansion of the Dulac map at $s=0$ in $(a)$--$(d)$ follow, respectively, from $(a')$--$(d')$. The only remaining point concerns the behaviour of the coefficients in each one of these cases. This is our final task, that we carry out case by case: 

\begin{enumerate}[$(a)$]

\item Let us consider any $\np_0=(D_0,F_0)\in W\cap\{1<F<\frac{5}{4}\}$ such that $T_{10}(\np_0)=0$. We claim that then $D_0\in (-1,-\frac{1}{2}).$ Indeed, to prove this we first use Proposition 3.11 in \cite{MMV2}, which shows that the set $\big\{\np\in W: T_{10}(\np)=0\text{ with }1<F<\frac{3}{2}\big\}$ is the graphic of an analytic function $D=\mathcal G(F)$ verifying $-F<\mathcal G(F)<-\frac{1}{2}$ and $\lim_{F\to 1^+}\mathcal G(F)=-\frac{1}{2}.$ Therefore it is clear that the claim will follow once we prove that $T_{10}(-1,F)\neq 0$ for all $F\in (1,\frac{5}{4}).$ In order to show this we note that $p_2|_{D=-1}=1$ and, consequently,
\[\textstyle
 T_{10}(-1,F)=\left.\rho_2(\np)\mathrm B\big(1-\frac{1}{\lambda},-\frac{1}{2}\big)
 {}_2F_1\big(-1-\frac{1}{\lambda},-\frac{1}{2};\frac{1}{2}-\frac{1}{\lambda};\frac{1-p_2}{1-p_1}\big)\right|_{D=-1}\neq 0
\] 
because ${}_2F_1(a,b;c;0)=1$ by definition and, on the other hand, $\lambda(\np)=\frac{1}{2(1-F)}>2$ for all $F\in (1,\frac{5}{4})$ and one can check that  
$\mathrm B\big(1-\frac{1}{\lambda},-\frac{1}{2}\big)=\frac{\Gamma(1-\frac{1}{\lambda})\Gamma(-\frac{1}{2})}{\Gamma(\frac{1}{2}-\frac{1}{\lambda})}\neq 0$ for all $\lambda>2$. This proves the claim. 

Recall at this point that
\begin{equation}\label{p2eq9}
T_{20}(\np)\textstyle=\rho_3(\np)\mathrm B\big(1-\frac{2}{\lambda},-\frac{3}{2}\big){}_2F_1\big(-\frac{2}{\lambda}-3,-\frac{3}{2};-\frac{1}{2}-\frac{2}{\lambda};\frac{1-p_2}{1-p_1}\big)+\rho_4(\np)T_{10}(\np).
\end{equation}
Accordingly, since $T_{10}(\np_0)=0$ with $\np_0\in Q\!:=W\cap\{1<F<\frac{5}{4}\}\cap\{-1<D<-\frac{1}{2}\}$ and $\rho_2$ and $\rho_3$ are positive functions, in order to prove that $T_{20}(\np_0)\neq 0$ it suffices to show that the linear combination 
\[\textstyle
\mathrm B\big(1-\frac{2}{\lambda},-\frac{3}{2}\big){}_2F_1\big(-\frac{2}{\lambda}-3,-\frac{3}{2};-\frac{1}{2}-\frac{2}{\lambda};\frac{1-p_2}{1-p_1}\big)-4\mathrm B\big(1-\frac{1}{\lambda},-\frac{1}{2}\big){}_2F_1\big(-1-\frac{1}{\lambda},-\frac{1}{2};\frac{1}{2}-\frac{1}{\lambda};\frac{1-p_2}{1-p_1}\big)
\] 
does not vanish on $Q$. Since one can easily verify that $\frac{2}{\lambda}\in (0,1)$ and $\frac{1-p_2}{1-p_1}\in (-1,0)$ for all $\np\in Q$, this follows directly by applying \propc{propC1} with $\{\al=\frac{2}{\lambda},\be=-\frac{1-p_2}{1-p_1}\}.$ 

\item We already proved that $T_{01}(\np)=\rho_1(\np)\mathrm B\big(-\lambda,\frac{1}{2}\big)$ with $\rho_1$ an analytic positive function on $W$. The function $\mathrm{B}(-\lambda,\frac{1}{2})=\frac{\Gamma(-\lambda)\Gamma(\frac{1}{2})}{\Gamma(\frac{1}{2}-\lambda)}$ vanishes only when $\frac{1}{2}-\lambda(\np)\in\Z_{\leq 0}$ and, for $\np\in W\cap\{\frac{5}{4}<F<\frac{3}{2}\}$, this occurs if and only if $\lambda(\np)=\frac{3}{2},$ i.e., $F=\frac{4}{3}.$ 
Recall moreover that, by Proposition~3.11 in \cite{MMV2}, the set $\big\{\np\in W: T_{10}(\np)=0\text{ with }1<F<\frac{3}{2}\big\}$ is the graphic of an analytic function $D=\mathcal G(F)$ on $(1,\frac{3}{2}).$ Consequently there exists a unique $\np_\star=(D_\star,\frac{4}{3})$ inside 
$W\cap\{\frac{5}{4}<F<\frac{3}{2}\}$ such that $T_{10}(\np_\star)=0$ and one can prove that $D_\star\in (-1.15,-1.10)$. The gradients of $T_{10}$ and $T_{01}$ are linearly independent at $\np_\star$ because $\partial_D T_{10}(\np)\neq 0$ for all $\np\in W\cap\{1<F<\frac{3}{2}\}$ by $(a)$ of Lemma 3.13 in~\cite{MMV2} and, on the other hand, $\partial_D T_{01}(\np_\star)=\mathrm B\big(-\lambda,\frac{1}{2}\big)|_{\np=\np_\star}\partial_D\rho_1(\np_\star)= 0$ and $\partial_F T_{10}(\np_\star)\neq 0$ since the gamma function has simple poles at $\Z_{\leq 0}$ with non-zero residue. Finally $T_{20}(\np_\star)<0$ follows noting that, from \refc{p2eq9} and $\lambda(\np_\star)=\frac{3}{2}$, we get
\[\textstyle
 T_{20}(\np_\star)=\rho_3(\np_\star)\left.\mathrm B\big(-\frac{1}{3},-\frac{3}{2}\big){}_2F_1\big(-\frac{13}{3},-\frac{3}{2};-\frac{11}{6};\frac{1-p_2}{1-p_1}\big)\right|_{\np=\np_\star},
\]
which is negative because one can easily check that $D\mapsto {}_2F_1\big(-\frac{13}{3},-\frac{3}{2};-\frac{11}{6};\frac{1-p_2}{1-p_1}\big)|_{F=\frac{4}{3}}$ is positive on $(-1.15,-1.10)$ and $\mathrm B\big(-\frac{1}{3},-\frac{3}{2}\big)\approx -2.6.$

\item Let us fix any $\np_0\in W\cap\{F=\frac{5}{4}\},$ so that $\lambda(\np_0)=2.$ Then, from $(c'),$ $T_{201}^2(\np)\!:=\bar T_{201}^2(\bar\np)|_{\eta=0}$ and $T_{200}^2(\np)\!:=\bar T_{200}^2(\bar\np)|_{\eta=0}$ are smooth functions in a neighbourhood of $\{\np\in W:\lambda(\np)=2\}$ and, in addition,
\[
  T_{201}^2(\np)=\big(2-\lambda(\np)\big) T_{01}(\np)
 \text{ and }
  T_{200}^2(\np)= T_{20}(\np)+ T_{01}(\np)
 \text{ for $\lambda(\np)\neq2.$}
\]
The functions $T_{01}(\np)$ and $T_{20}(\np)$ are meromorphic with a pole 
at those $\np$ such that $\lambda(\np)=2\in D_{01}\cap D_{20}$, and the pole is simple in both cases by Propositions 3.2 and 3.6 in \cite{MV21}, respectively. Therefore by the Weierstrass Division Theorem (or, more directly, by \cite[Lemma 2.8]{MV21}) it follows that $T_{201}^2(\np)$ and $T_{200}^2(\np)$ are analytic in a neighbourhood of~$W\cap\{F=\frac{5}{4}\}$.
Furthermore, by \cite[Lemma 3.13]{MMV2} once again, $T_{10}(D,\frac{5}{4})=0$ if and only if $D=-1$. Finally, since $\lambda(-1,\frac{5}{4})=2,$ $T_{01}(\np)=\rho_1(\np)\mathrm B(-\lambda,\frac{1}{2})$ and $\mathrm{B}(-\lambda,\frac{1}{2})=\frac{\Gamma(-\lambda)\Gamma(\frac{1}{2})}{\Gamma(\frac{1}{2}-\lambda)}$, we have that
\[
 T_{201}^2(-1,5/4)=\lim_{\np\to (-1,\frac{5}{4})}\big(2-\lambda(\np)\big) T_{01}(\np)
 =\rho_1(-1,5/4)\textstyle\frac{\Gamma(\frac{1}{2})}{\Gamma(-\frac{3}{2})}\lim\limits_{\lambda\to 2}(2-\lambda)\Gamma(-\lambda)\neq 0,
\]
because the gamma function has simple poles with non-zero residues at~$\Z_{\leq 0}.$

\item Let us consider finally any $\np_0\in W\cap\{F=2\},$ so that $\lambda(\np_0)=\frac{1}{2}.$ Due to $T_{01}(\np)=\rho_1(\np)\mathrm B\big(-\lambda,\frac{1}{2}\big)$ with $\rho_1\neq 0$, $\lambda(\np)=\frac{1}{2(F-1)}$ and $\mathrm B\big(-\lambda,\frac{1}{2}\big)=\frac{\Gamma(-\lambda)\Gamma(\frac{1}{2})}{\Gamma(\frac{1}{2}-\lambda)}$, there exists an analytic non-vanishing function $\ell_1$ in a neighbourhood of $W\cap\{F=2\}$ such that $T_{01}(\np)=(F-2)\ell_1(\np).$ On the other hand, from~$(d')$ and arguing as in the previous case, the functions $T_{101}^{\frac{1}{2}}(\np)\!:=\bar T_{101}^{\frac{1}{2}}(\bar\np)|_{\eta=0}$ and $T_{100}^{\frac{1}{2}}(\np)\!:=\bar T_{100}^{\frac{1}{2}}(\bar\np)|_{\eta=0}$ are analytic in a neighbourhood of $\{\np\in W:\lambda(\np)=\frac{1}{2}\}$ and
\[
  T_{101}^{\frac{1}{2}}(\np)=(1-2\lambda) T_{02}(\np)
  \text{ and }
  T_{100}^{\frac{1}{2}}(\np)= T_{10}(\np)+ T_{02}(\np)
 \text{ for $\lambda(\np)\neq 1/2.$}
\]
In particular we have that the sum of residues of $T_{10}$ and $T_{02}$ along  $\{\np\in W:\lambda(\np)=\frac{1}{2}\}=W\cap\{F=2\}$
is equal to zero,
which saves us from computing the explicit value of $T_{02}.$ Indeed, since $\lambda(\np_0)=\frac{1}{2}$ and $\mathrm{B}(1-\frac{1}{\lambda},-\frac{1}{2})=\frac{\Gamma(1-\frac{1}{\lambda})\Gamma(-\frac{1}{2})}{\Gamma(\frac{1}{2}-\frac{1}{\lambda})},$ we obtain that
\begin{align*}\textstyle
 T_{101}^{\frac{1}{2}}(\np_0)&=\lim_{\np\to\np_0}(1-2\lambda(\np))T_{02}(\np)=-\lim_{\np\to\np_0}(1-2\lambda(\np))T_{10}(\np)\\
 &
 =\textstyle -\rho_2(\np_0){}_2F_1\big(-3,-\frac{1}{2};-\frac{3}{2};\frac{1-p_2}{1-p_1}\big|_{\np=\np_0}\big)\frac{\Gamma(-\frac{1}{2})}{\Gamma(-\frac{3}{2})}\lim_{\lambda\to\frac{1}{2}}(1-2\lambda)\Gamma(1-\frac{1}{\lambda})\\
 &
 =\textstyle \frac{3}{4}\rho_2(\np_0){}_2F_1\big(-3,-\frac{1}{2};-\frac{3}{2};\frac{1-p_2}{1-p_1}\big|_{\np=\np_0}\big),
\end{align*}
where in the second equality we use the expression of $T_{10}$ already proved and in the last one that 
$\frac{\Gamma(-\frac{1}{2})}{\Gamma(-\frac{3}{2})}=-\frac{3}{2}$ and $\lim_{x\to -1}(x+1)\Gamma(x)=-1,$ see for instance \cite[\S 15]{stegun}. From the same reference we get that ${}_2F_1\big(-3,-\frac{1}{2};-\frac{3}{2};x)=(x+1)(x-1)^2$. Moreover one can verify that $D\mapsto\frac{1-p_2}{1-p_1}\big|_{F=2}$ maps diffeomorphically $(-2,0)$ to $(-\infty,1)$ and that it is equal to $-1$ at $D=-\frac{1}{2}.$ Accordingly we can assert that $T_{101}^{\frac{1}{2}}(D,2)=(D+\frac{1}{2})\ell_2(D)$ where $\ell_2$ is a non-vanishing analytic function on $(-2,0)$ and, consequently, $\partial_D T_{101}^{\frac{1}{2}}(-\frac{1}{2},2)\neq 0.$ Since $\partial_DT_{01}(-\frac{1}{2},2)=0$ and $\partial_FT_{01}(-\frac{1}{2},2)=\ell_1(-\frac{1}{2},2)\neq 0,$ this proves that the gradients of $T_{01}$ and $T_{101}^{\frac{1}{2}}$ are linearly independent at $\np=(-\frac{1}{2},2)$ as desired.
\end{enumerate}
This finishes the proof of the result.
\end{prooftext}


\section{Beta and hypergeometric functions}\label{ApBeta}

In this appendix we are concerned with the integral representation of the Beta and hypergeometric functions (see \cite{Askey} for details). The Beta integral is defined for $\mathrm{Re}(z)>0$ and $\mathrm{Re}(w)>0$ by
\begin{equation}\label{eqB1}
 \mathrm B(z,w)\!:=\int_0^1t^{z-1}(1-t)^{w-1}dt=\int_0^{+\infty}t^{z-1}(1+t)^{-z-w}dt.
\end{equation}
This function can be analytically extended for $z,w\in\C\setminus\Z_{\leq 0}$ thanks to the identity
 \begin{equation}\label{eqB0}
  \mathrm B(z,w)=\frac{\Gamma(z)\Gamma(w)}{\Gamma(z+w)},
 \end{equation} 
where $\Gamma$ is the gamma function. Recall in this regard that $1/\Gamma(z)$ is an entire function with simple zeros at $z\in \Z_{\leq 0}.$ On the other hand, if we consider $a,b,c\in\C$ with $c\notin\Z_{\leq 0}$ and $z$ inside the complex open unit disc $\Dc\!:=\{z\in\C:|z|<1\}$, then Gauss hypergeometric function is defined by the power series
 \begin{equation}\label{hyper}
 _{2}F_1(a,b;c;z)\!:=\sum_{n=0}^\infty\frac{(a)_n(b)_n}{(c)_n}\frac{z^n}{n!},
 \end{equation}
where for a given $x\in\C$ we use the Pochhammer symbol $(x)_n\!:=x(x+1)\cdots(x+n-1)=\frac{\Gamma(x+n)}{\Gamma(x)}$. 

In this section by a meromorphic function of several complex variables we mean a function that locally writes as a quotient of two holomorphic functions. Recall that a function \map{f}{\Omega}{\C}, where $\Omega$ is a connected open set of $\C^n,$ is holomorphic if for each $z_0\in\Omega$ there exists an open polydisc $D_r(z_0)$ such that $f$ can be written as an absolutely and uniformly convergent power series at $z_0,$ i.e., $f(z)=\sum_{\alpha}a_{\alpha}(z-z_0)^{\alpha}$ for all $z\in D_r(z_0)$. On account of this we have the following result about uniqueness of meromorphic continuation. 

\begin{lem}\label{continuacio}
Consider two functions $\phi$ and $\varphi$ that are meromorphic on a connected open set $\Omega\subset\C^n$. If 
there exists an open subset $V$ of $\Omega$ such that $\left. \phi\right|_{V}=\left. \varphi\right|_{V}$ then $\phi=\varphi.$ 
\end{lem}

\begin{prova}
We assume without loss of generality that $\varphi=0$. In doing so the equality $\phi=0$ has to be thought only at regular points of $\phi.$ That being said, consider any two regular points of $\phi$, $z_0\in V$ and $z_1\in\Omega\setminus V$ and take a continuous path $\gamma$ joining them. Suppose that it is parameterized by $\sigma:[0,1]\longrightarrow\gamma$ with $\sigma(0)=z_0$ and $\sigma(1)=z_1.$ By compactness there exist $0\leqslant s_1< s_2<\ldots< s_k\leqslant 1$ and positive numbers $r_1,r_2,\ldots,r_k$ verifying $\gamma\subset\cup_{i=1}^kD_{r_i}(\sigma(s_i))$ and such that, for each $i\in\{1,2,\ldots,k\},$ we can write $\left.\phi\right|_{D_{r_i}(\sigma(s_i))}=f_i/g_i$ with $f_i$ and $g_i$ holomorphic on $D_{r_i}(\sigma(s_i))$.
Define $\ell=\max\{i:D_{r_i}(\sigma(s_i))\cap V\neq\emptyset\}$. Then on the regular set of $\left.\phi\right|_{D_{r_{\ell}}(\sigma(s_{\ell}))}=f_{\ell}/g_{\ell}$ (i.e., where $g_\ell\neq 0$), the equality $\left.f_\ell\right|_V=0$ implies $f_\ell=0$ by the uniqueness of analytic continuation. Accordingly $\left.\phi\right|_{D_{r_{\ell}}(\sigma(s_{\ell}))}=0.$ Next we compute $\ell$ again but replacing $V$ by $V\cup D_{r_{\ell}}$ and we iterate the process to conclude that $\phi(z_1)=0.$ Since $z_1$ is arbitrary this proves the result. 
\end{prova}

Let us remark that the previous result is also true (with the same proof) in the real setting, i.e., for functions in $\R^n$ that locally write as a quotient of real analytic functions. The following result is well-known but since we did not find its statement in its fullness we give it here for the sake of completeness.

\begin{lem}\label{2F1}
The function $(a,b,c,z)\mapsto \frac{_{2}F_1(a,b;c;z)}{\Gamma(c)}$ extends holomorphically to $\C^3\!\times\!(\C\setminus[1,+\infty))$.
\end{lem}
 
\begin{prova}
Following \cite[p. 65]{Askey}, let us show first that the function extends holomorphically to $\C^3\!\times\Dc$. To prove  this claim we write 
 \[
  \frac{_{2}F_1(a,b;c;z)}{\Gamma(c)}=\sum_{n=0}^{\infty}f_n(a,b,c,z)\text{ with $f_n(a,b,c,z)\!:=\frac{\Gamma(a+n)\Gamma(b+n)z^n}{\Gamma(a)\Gamma(b)\Gamma(c+n)\Gamma(1+n)}$.}
 \] 
Stirling's asymptotic formula $\Gamma(z)\sim \sqrt{2\pi}z^{z-\frac{1}{2}}e^{-z}$ as $\mr{Re}(z)\to +\infty$ (see \cite[Theorem~1.4.1]{Askey}) shows that
\[
 \left|\frac{\Gamma(a+n)\Gamma(b+n)}{\Gamma(c+n)\Gamma(1+n)}\right|\sim n^{\mathrm{Re}(a+b+c-1)}\text{ as $n\to +\infty.$}
\]
Fix any compact set $K\subset \C^3\!\times\Dc$ and suppose that $\mr{Re}(a+b+c-1)\leqslant m\in\N$ and $|z|\leqslant r<1$ for all $(a,b,c,z)\in K.$ Then, on account of the above asymptotic estimate and the fact that $1/\Gamma(z)$ is an entire function, there exists $C>0$ such that $|f_n(a,b,c,z)|\leqslant C n^m r^n$ for all $(a,b,c,z)\in K.$ By applying the Weierstrass M-test this proves that the series $\sum_{n=0}^{\infty}f_n$ converges uniformly 
on compact sets of $\C^3\!\times\Dc$. So the claim follows because the uniform limit of holomorphic functions is holomorphic (see \cite[Proposition~2]{Malgrange}). 

Finally the result follows by Pfaff and Kummer's formulas (see  \cite[\S 15]{stegun} or \cite[Theorem 2.3.2]{Askey}) relating the values of ${}_2F_1(a,b;c;\,\cdot\,)$ at $z$, $\frac{z}{z-1}$ and $\frac{1}{z}$, which enable to extend holomorphically to $\C^3\times(\C\setminus[1,+\infty))$ the map $(a,b,c,z)\mapsto \frac{_{2}F_1(a,b;c;z)}{\Gamma(c)}$. (These formulas are usually proved under some restrictions on the parameters $a$ $b$ and $c$ but they are always satisfied thanks to the claim and \lemc{continuacio}.) This concludes the proof of the result.
\end{prova}

It is worth to mention that by Hartogs's theorem (see~\cite[\S 1.2]{Krantz}), a function of several complex variables is holomorphic if, and only if, it is holomorphic (in the classical one-variable sense) in each variable separately. 
The main concern in this section is Euler's integral representation of ${}_2F_1$, see for instance \cite[Theorem~2.2.1]{Askey}, that is given by
\begin{equation}\label{Euler}
{}_2F_1(a,b;c;z)=\frac{1}{ \mathrm B(b,c-b)}\int_0^1t^{b-1}(1-t)^{c-b-1}(1-zt)^{-a}dt
\end{equation}
provided that $\mathrm{Re}(c)>\mathrm{Re}(b)>0$ and $z\in\C\setminus[1,+\infty)$. Our goal is to use this formula to compute~$\hat f(\alpha,x)$  (see \teoc{L8}) for some specific functions $f(x)$. Next result is addressed to this problem.
 
 \begin{prop}\label{B2F1}
The following holds:
\begin{enumerate}[$(a)$]

\item Consider $h(y;\delta;\kappa)=(1+\kappa y^2)^{\delta}$ with $\kappa>0.$ Then, for any 
         $\delta\in\R$ and $\alpha\in\R\setminus\Z_{\geq 0}$ such that $2\delta<\alpha,$
         \[
          \lim_{y\to+\infty}y^{-\alpha}\hat h(\alpha,y;\delta;\kappa)
          =\frac{\kappa^{\frac{\alpha}{2}}}{2} \mathrm B\left(-\frac{\alpha}{2},-\delta+\frac{\alpha}{2}\right).
          \]

\item Consider $g(y;\delta,\ga;x)=(1-y)^{\delta-1}(1-xy)^{-\ga}$ with $y\in (0,1)$ and $x<1$. Then, for any 
         $\delta>0$, $\ga\in\R$ and $\alpha\in\R\setminus\Z_{\geq 0}$,
         \[
          \lim_{y\to 1^-}\hat g(\alpha,y;\delta,\ga;x)= 
          \mathrm B(-\alpha,\delta){}_2F_1(\ga,-\alpha;\delta-\alpha;x).
         \]
         
\end{enumerate}
\end{prop}

\begin{prova}
In order to prove $(a)$ we define $\Omega\!:=\{(\alpha,\delta,\kappa)\in \R^3:2\delta<\alpha\text{ and }  \kappa>0\}$, which is connected. Note then that we must show the validity of the identity on $\Omega\cap\{\alpha\notin\Z_{\geq 0}\}$. We will show first the identity on an open set of $V\subset\Omega$ and then extend it by using the real version of \lemc{continuacio}. With this aim observe that if we work on $V\!:=\Omega
\cap\{\alpha<0\}$ then the application of assertion $(b)$ of \teoc{L8} with $k=0$ yields
\[
 \lim_{y\to+\infty}y^{-\alpha}\hat h(\alpha,y;\delta,\kappa)=\int_0^{+\infty}(1+\kappa u^2)^{\delta}u^{-\alpha-1}du
 =\frac{\kappa^{\frac{\alpha}{2}}}{2}\int_0^{+\infty}(1+v)^{\delta}v^{-\frac{\alpha}{2}-1}dv=\frac{\kappa^{\frac{\alpha}{2}}}{2} \mathrm B\left(-\frac{\alpha}{2},-\delta+\frac{\alpha}{2}\right),
\]
where in the second equality we perform the change of variable $v=\kappa u^2$ and in the third one we use~\refc{eqB1}. Note at this point that the function of the right hand side is meromorphic on $\Omega$ because $1/\Gamma$ is entire and $ \mathrm B(z,w)=\frac{\Gamma(z)\Gamma(w)}{\Gamma(z+w)}$. We claim that the function on the left hand side is also meromorphic on $\Omega$. To show this we work first on $\Omega\cap\{\alpha\notin\Z_{\geq 0}\}$ because in doing so we can apply assertion $(b)$ of
\teoc{L8} with any $k>\alpha$ to obtain 
 \begin{align*}
  y^{-\alpha}\hat h(\alpha, y)&=\sum_{i=0}^{k-1}\frac{h^{(i)}(0)}{i!(i-\alpha)}y^{i-\alpha}+\int_0^y\big(h(u)-T_0^{k-1}h(u)\big)\frac{du}{u^{\alpha+1}}\\
  &=\sum_{i=0}^{k-1}\frac{h^{(i)}(0)}{i!(i-\alpha)}y^{i-\alpha}+\int_0^1\big(h(u)-T_0^{k-1}h(u)\big)\frac{du}{u^{\alpha+1}}
  +\int_1^yh(u)\frac{du}{u^{\alpha+1}}-\int_1^yT_0^{k-1}h(u)\frac{du}{u^{\alpha+1}}\\
  &=\sum_{i=0}^{k-1}\frac{h^{(i)}(0)}{i!(i-\alpha)}+\int_0^1\big(h(u)-T_0^{k-1}h(u)\big)\frac{du}{u^{\alpha+1}}
  +\int_1^yh(u)\frac{du}{u^{\alpha+1}}\\
  &=\hat h(\alpha,1;\delta,\kappa)+\int_1^y(1+\kappa u^2)^{\delta}u^{-\alpha-1}du. 
 \end{align*}
Here we denote $\partial^i_yh(y;\delta,\kappa)=h^{(i)}(y)$ for shortness. Consequently, since $2\delta<\alpha$,
 \begin{align*}
  \lim_{y\to+\infty}y^{-\alpha}\hat h(\alpha,y;\delta,\kappa)&=\hat h(\alpha,1;\delta,\kappa)+\int_1^{+\infty}(1+\kappa u^2)^{\delta}u^{-\alpha-1}du
  \\
  &=\hat h(\alpha,1;\delta,\kappa)+\kappa^{\delta}\int_0^{1}(1+\kappa^{-1} v^2)^{\delta}v^{\alpha-2\delta-1}dv\\
  &=\hat h(\alpha,1;\delta,\kappa)+\kappa^{\delta}\hat h(2\delta-\alpha,1;\delta,\kappa^{-1}),
 \end{align*}
where in the second equality we make the change of variable $v=1/u$ and in the last one we apply $(b)$ in \teoc{L8} with $k=0$. By $(c)$ in \teoc{L8} the second summand is analytic on $\Omega$, whereas the first one is meromorphic on~$\Omega$. This shows the validity of the claim and so the result follows by applying the real version of \lemc{continuacio}.

In order to prove $(b)$ we fix $\delta,$ $\ga$ and $x$ and apply
$(b)$ in \teoc{L8} to the function $y\mapsto g(y;\delta,\ga;x).$ Then, taking any $k\in\N$ with $k>\alpha$ and setting $\partial^i_yg(y;\delta,\ga;x)=g^{(i)}(y)$ for shortness, we get
\begin{align*}
\lim_{y\to 1^-}\hat g(\alpha,y)=&\sum_{r=0}^{k-1}\frac{g^{(r)}(0)}{r! (r-\alpha)}+\int_0^1
\left(g(u)-\sum_{r=0}^{k-1}\frac{g^{(r)}(0)}{r!}u^r\right)\frac{du}{u^{\alpha+1}}
\\
=&\, 2^\alpha\left(\sum_{r=0}^{k-1}\frac{g^{(r)}(0)}{r! (r-\alpha)2^r}+2^{-\alpha}\int_0^{\frac{1}{2}}\left(g(u)-\sum_{r=0}^{k-1}\frac{g^{(r)}(0)}{r!}u^r\right)\frac{du}{u^{\alpha+1}}\right)
\\
&+\int_{\frac{1}{2}}^1(1-u)^{\delta-1}(1-xu)^{-\ga}u^{-\alpha-1}du-\sum_{r=0}^{k-1}\frac{g^{(r)}(0)}{r!}
\int_{\frac{1}{2}}^1u^{r-\alpha}\frac{du}{u}+\sum_{r=0}^{k-1}\frac{g^{(r)}(0)}{r! (r-\alpha)}(1-2^{\alpha-r})
\\
=&\, 2^\alpha\hat g\left(\alpha,1/2;\delta,\ga;x\right)+\int_0^{\frac{1}{2}}s^{\delta-1}(1-x(1-s))^{-\ga}(1-s)^{-\alpha-1}ds
\\
=&\, 2^\alpha \hat g\left(\alpha,1/2;\delta,\ga;x\right)+2^{-\delta}(1-x)^{-\ga}\hat g\left(-\delta,\frac{1}{2};-\alpha,\ga;\frac{x}{x-1}\right)
\end{align*}
where in the last equality we use $(b)$ in \teoc{L8}
with $k=0$ and also take $\delta>0$ into account. Thus, by applying $(c)$ in \teoc{L8}
to each summand in the last expression, this shows that the function 
\[ 
 (\ga,\alpha,\delta,x)\mapsto \lim_{y\to 1^-}\hat g(\alpha,y;\delta,\ga;x)
\]  
is meromorphic on the open connected set $\hat\Omega\!:=\R^2\times (0,+\infty)\times (-\infty,1)$. Note also that if we consider parameter values in $\hat V\!:=\hat\Omega\cap\{\alpha<0\}\cap\{\delta>0\}$ then
 \[
  \lim_{y\to 1^-}\hat g(\alpha,y;\delta,\ga;x)=\int_0^1(1-u)^{\delta-1}(1-xu)^{-\ga}u^{-\alpha-1}du
  = \mathrm B(-\alpha,\delta){}_2F_1(\ga,-\alpha;\delta-\alpha;x),
 \]
where in the first equality we apply $(b)$ in \teoc{L8}
with $k=0$ and in the second one we use Euler's integral representation~\refc{Euler}. We have just proved that the left hand side expression is a meromorphic function on $\hat\Omega$. Furthermore, by applying \lemc{2F1} and taking $ \mathrm B(-\alpha,\delta)=\frac{\Gamma(-\alpha)\Gamma(\delta)}{\Gamma(\delta-\alpha)}$ into account, we can assert that the right hand side is also a a meromorphic function on $\hat\Omega$. In view of this the identity in $(b)$ for the parameters under consideration follows by applying the real version of \lemc{continuacio}.
\end{prova}

It is worth to point out that the application of $(b)$ in \propc{B2F1} provides integral representations of the hypergeometric function in a range of parameters not covered by Euler's formula~\refc{Euler}. Indeed, by applying also $(b)$ in \teoc{L8} with any $k>\alpha$ we get
 \[
  \mathrm B(-\alpha,\delta){}_2F_1(\ga,-\alpha;\delta-\alpha;x)=\sum_{r=0}^{k-1}\frac{g^{(r)}(0)}{r! (r-\alpha)}+\int_0^1
t^{-\alpha-1}\left((1-t)^{\delta-1}(1-xt)^{-\ga}-\sum_{r=0}^{k-1}\frac{g^{(r)}(0)}{r!}t^r\right)dt,
 \]
which holds for any $\delta>0$, $\ga\in\R$ and  $\alpha\in\R\setminus\Z_{\geq 0}$, where $g(t)=(1-t)^{\delta-1}(1-xt)^{-\ga}$. We stress that~\refc{Euler} gives an integral representation for ${}_2F_1(\ga,-\alpha;\delta-\alpha;x)$ only in case that $\delta>0$ and $\alpha<0.$


\section{A technical result for the proof of \propc{prop52}}\label{apB}

\begin{prop}\label{propC1} The function
\[\textstyle
 \Phi(\al,\be)=\mathrm B(1-\al,-\frac{3}{2}){}_2F_1(-3-\al,-\frac{3}{2};-\frac{1}{2}-\al;-\be)-4\mathrm B(1-\frac{\al}{2},-\frac{1}{2}){}_2F_1(-1-\frac{\al}{2},-\frac{1}{2};\frac{1}{2}-\frac{\al}{2};-\be)
 \]
is strictly positive for all $\al,\be\in (0,1)$.
\end{prop}

\begin{prova}
In what follows given a smooth function of several variables $\psi(u)$ with $u=(u_1,u_2,\ldots,u_n)\in\R^n$, for each fixed $i\in\{1,2,\ldots,n\}$ and $m\in\N$
we denote 
by $T^m_{u_i=0}\psi(u)$ the $m$-th order Taylor polynomial of $u_i\mapsto \psi(u)$ at $u_i=0,$ i.e., $T^m_{u_i=0}\psi(u)=\sum_{k=0}^m\frac{1}{k!}\left.\frac{\partial^k\psi}{\partial u_i^k}\right|_{u_i=0}u_i^k$. 
Setting $P(\al,\be)\!:=T^3_{\be=0}\Phi(\al,\be)$ we claim that the following inequalities hold:
\begin{enumerate}[$(i)$]
\item $\Phi(\al,\be)\geqslant P(\al,\be)$ for all $\al,\be\in(0,1)$, and
\item $P(\al,\be)>0$ for all $\al,\be\in(0,1)$.
\end{enumerate}
It is clear that the result will follow once we prove this. For this purpose our first task will be to express the function $\Phi$ in terms of a definite integral and to this end we define
\begin{align*}\textstyle
 \Phi_1(\al,\be)\!:=&\textstyle\mathrm B(-\frac{3}{2},1-\al){}_2F_1(-3-\al,-\frac{3}{2};-\frac{1}{2}-\al;-\be)\\
 \intertext{and}
 \Phi_2(\al,\be)\!:=&\textstyle\mathrm B(-\frac{1}{2},1-\frac{\al}{2}){}_2F_1(-1-\frac{\al}{2},-\frac{1}{2};\frac{1}{2}-\frac{\al}{2};-\be),
 \end{align*}
so that, taking $\mathrm B(x,y)=\mathrm B(y,x)$ into account, $\Phi=\Phi_1-4\Phi_2$. Then by applying $(b)$ in 
\propc{B2F1} with $\{\alpha=\frac{3}{2}, \delta=1-\al, \ga=-3-\al, x=-\be\}$ we can assert that
$\Phi_1(\al,\be)=\lim_{y\to 1^-}\hat g_1(\frac{3}{2},y;\al,\be)$ where
$g_1(y;\al,\be)\!:=(1-y)^{-\al}(1+\be y)^{3+\al}.$ Next we apply $(b)$ in \teoc{L8}
taking $\{f=g_1, \alpha=\frac{3}{2},k=2,x=1\}$ to obtain that
\begin{align*}
\Phi_1(\al,\be)=\lim_{y\to 1^-}\hat g_1(3/2,y;\al,\be)&=
2\kappa_1+\int_0^1\big(g_1(s)-r_1(s)\big)s^{-\frac{5}{2}}ds\\
&=\int_0^1\big(g_1(s)-r_1(s)+\kappa_1s^2\big)s^{-\frac{5}{2}}ds,
\end{align*} 
where $\kappa_1=-\frac{1}{3}-\bar\kappa$ and $r_1(s;\al,\be)=1+\bar\kappa s$ with $\bar\kappa\!:=(3+\al)\be+\al.$
Similarly by $(b)$ in 
\propc{B2F1} with $\{ \alpha=\frac{1}{2}, \delta=1-\frac{\al}{2},\ga=-1-\frac{\al}{2}, x=-\be\}$ we obtain that
$\Phi_2(\al,\be)=\lim_{y\to 1^-}\hat g_2(\frac{1}{2},y;\al,\be)$ where
$g_2(y;\al,\be)\!:=(1-y)^{-\frac{\al}{2}}(1+\be y)^{1+\frac{\al}{2}}.$  Next we apply $(b)$ in \teoc{L8}
with $\{f=g_2, \alpha=\frac{1}{2},k=2,x=1\}$ to get that
\begin{align*}
\Phi_2(\al,\be)=\lim_{y\to 1^-}\hat g_2(1/2,y;\al,\be)&=2\kappa_2+\int_0^1\big(g_2(s)-r_2(s)\big)s^{-\frac{3}{2}}ds\\
&=\int_0^1\big(sg_2(s)-sr_2(s)+\kappa_2s^2\big)s^{-\frac{5}{2}}ds,
\end{align*}
where $\kappa_2=-1+\bar\kappa$ and $r_2(s;\al,\be)=1+\bar\kappa s$. Accordingly an easy computation yields 
\[
 \Phi(\al,\be)=\Phi_1(\al,\be)-4\Phi_2(\al,\be)=\int_0^1h(s;\al,\be)s^{-\frac{5}{2}}ds,
\]
where 
\[
h(s;\al,\be)\!:=(1+\be s)^{3+\al}(1-s)^{-\al}-4s(1+\be s)^{1+\frac{\al}{2}}(1-s)^{-\frac{\al}{2}}
-1+(4-\bar\kappa)s+(11/3-\bar\kappa)s^2
\]
and consequently
\[
P(\al,\be)=T^3_{\be=0}\Phi(\al,\be)=\int_0^1h_0(s;\al,\be)s^{-\frac{5}{2}}ds\text{ where $h_0\!:=T_{\be=0}^3h.$}
\]
On account of this definition and the fact that if $\frac{\partial\varphi}{\partial u_i}\equiv 0$ then $T_{u_i=0}^m(\varphi\psi)=\varphi T_{u_i=0}^m(\psi)$, we get 
\[
 \Phi(\al,\be)-P(\al,\be)=\int_0^1\big(h(s;\al,\be)-h_0(s;\al,\be)\big)s^{-\frac{5}{2}}ds
 =\int_0^1\frac{g(s;\al,\be)-g_0(s;\al,\be)}{(1-s)^\al s^{\frac{5}{2}}}ds,
\]
where
\[
 g(s;\al,\be)\!:=(1+\be s)^{3+\al}-4s(1-s)^{\frac{\al}{2}}(1+\be s)^{1+\frac{\al}{2}}
 \text{ and }
 g_0\!:=T_{\be=0}^3g.
\]
Therefore the assertion in $(i)$ will follow once we prove that $g(x;\al,\be)\geqslant g_0(x;\al,\be)$ for all $\al,\be,x\in (0,1).$ As a first step to this aim let us prove that, setting
\begin{align*}
 &\ell(\al,y)\!:=(1+y)^{3+\al}-4(1+y)^{1+\frac{\al}{2}}\text{ and }\ell_0\!:=T^3_{y=0}\ell,
\intertext{then}
&g(x;\al,\be)- g_0(x;\al,\be)\geqslant \ell(\al,\be x)-\ell_0(\al,\be x)
 \text{ for all $\al,\be,x\in (0,1)$.}
\end{align*}
In order to show this we note that
 \begin{align*}
 g(x;\al,\be)- g_0(x;\al,\be)
    &=(1+\be x)^{3+\al}-T_{\be=0}^3(1+\be x)^{3+\al}\\
    &\qquad-4x(1-x)^{\frac{\al}{2}}\left((1+\be x)^{1+\frac{\al}{2}}-T_{\be=0}^3(1+\be x)^{1+\frac{\al}{2}}\right)\\
    &\geqslant (1+\be x)^{3+\al}-T_{\be=0}^3(1+\be x)^{3+\al}
    -4\left((1+\be x)^{1+\frac{\al}{2}}-T_{\be=0}^3(1+\be x)^{1+\frac{\al}{2}}\right)
 \end{align*}
because $x(1-x)^{\frac{\al}{2}}\leqslant 1$ and, thanks to the remainder's formula in Taylor's Theorem, one can easily verify that $(1+y)^{\eta}-T_{y=0}^m(1+y)^{\eta}\geqslant 0$ for any $m$ odd and $\eta\in (1,2)$. It is clear then that a sufficient condition for the inequality in $(i)$ to hold is that $\ell(\al,y)-\ell_0(\al,y)\geqslant 0$ for all $\al,y\in (0,1).$
In order to show that this is indeed true we note that, by the remainder's formula in Taylor's Theorem again, 
\[
 \ell(\al,y)-\ell_0(\al,y)=\frac{\partial^4_y\ell(\al,y_0)}{4!}y^4\text{ for some $y_0\in(0,y)$,} 
\]
and consequently it suffices to verify that 
\[\textstyle
 \partial^4_y\ell(\al,y)=(\al+3)(\al+2)(\al+1)\al(1+y)^{\al-1}-\frac{1}{4}(\al+2)\al(\al-1)(\al-4)(1+y)^{\frac{\al}{2}-3}
 \geqslant 0
 \]
for all $\al,y\in (0,1).$ This inequality is equivalent to $(1+y)^{2+\frac{\al}{2}}\geqslant\frac{1}{4}\frac{(\al-2)(\al-4)}{(\al+3)(\al+1)}$, which is obviously true due to $(1+y)^{2+\frac{\al}{2}}\geqslant 1\geqslant \frac{1}{4}\frac{(\al-2)(\al-4)}{(\al+3)(\al+1)}$ for all $\al,y>0$. This proves the first assertion of the claim.  

We now turn to the proof of the inequality in $(ii)$. On account of the definition of the Gauss hypergeometric function, see \refc{hyper}, together with the definition of the function $\Phi(\al,\be)$ given in the statement it easily follows that
\begin{align}\notag
 &P(\al,\be)=T^3_{\be=0}\Phi(\al,\be)=\rho_0(\al)-\rho_1(\al)\be+\rho_2(\al)\be^2-\rho_3(\al)\be^3
 \intertext{with}\label{eqc1}
 &\rho_n(\al)\!:=\frac{(-3-\al)_n(-\frac{3}{2})_n}{(-\frac{1}{2}-\al)_nn!}
      \mathrm B\!\left(1-\al,-\frac{3}{2}\right)
    -4\frac{(-1-\frac{\al}{2})_n(-\frac{1}{2})_n}{(\frac{1}{2}-\frac{\al}{2})_nn!} 
     \mathrm B\!\left(1-\frac{\al}{2},-\frac{1}{2}\right).
\end{align}
Observe that $P(a,b)$ is a polynomial in $b$ for each fixed $a$. In order to prove the inequality in $(ii)$ we consider $P(a,\,\cdot\,)$ as a polynomial family depending on a parameter $a\in (0,1).$ In doing so it is clear that the following three conditions imply that the number of zeros of $b\mapsto P(a,b)$ on the interval $(0,1)$, counted with multiplicities, is the same for all $a\in (0,1)$:
\begin{enumerate}[$(a)$]
\item $P(\al,0)>0$ for all $\al\in (0,1),$
\item $P(\al,1)>0$ for all $\al\in (0,1)$ and
\item $\text{Discrim}_{\be}P(\al,\be)\neq 0$ for all $\al\in (0,1).$
\end{enumerate}
Since one can readily show that, for instance, $P(\frac{1}{2},b)>0$ for all $b\in (0,1),$ it is clear that $(ii)$ will follow once we prove that these three conditions are true. This constitutes our next task.
In order to prove the inequality in $(a)$ we first note that, from \refc{eqB0}, 
\[
 P(\al,0)=\rho_0(\al)=\mathrm B\!\left(1-\al,-\frac{3}{2}\right)-4\mathrm B\!\left(1-\frac{\al}{2},-\frac{1}{2}\right)
 =8\sqrt{\pi}\left(\frac{\Gamma(1-\frac{\al}{2})}{\Gamma(\frac{1}{2}-\frac{\al}{2})}-\frac{\al+\frac{1}{2}}{6}\frac{\Gamma(1-\al)}{\Gamma(\frac{1}{2}-\al)}\right),
\]
where we use that $\Gamma(\frac{1}{2})=\sqrt{\pi}$ and $\Gamma(z+1)=z\Gamma(z),$ see~\cite[\S 6.1]{stegun}.
Taking this into account, the fact that $P(\al,0)>0$ for all $\al\in [\frac{1}{2},1)$ is clear because $\Gamma(z)$ is negative for $z\in (-1,0)$ and positive for $z>0$ and $\lim_{z\to 0}\frac{1}{\Gamma(z)}=0$, see~\cite[\S 6.1]{stegun} again. To show that this is also true for 
$\al\in (0,\frac{1}{2})$ we use that then
\[
\frac{\Gamma(1-\frac{\al}{2})}{\Gamma(\frac{1}{2}-\frac{\al}{2})}>\frac{\Gamma(1-\al)}{\Gamma(\frac{1}{2}-\al)}>\frac{\al+\frac{1}{2}}{6}\frac{\Gamma(1-\al)}{\Gamma(\frac{1}{2}-\al)}. 
\]
The second inequality above is obvious whereas the first one follows noting that $z\mapsto\frac{\Gamma(1-z)}{\Gamma(\frac{1}{2}-z)}$ is positive and decreasing on $(0,\frac{1}{2})$. In its turn this is true due to 
 $\left(\log\frac{\Gamma(1-z)}{\Gamma(\frac{1}{2}-z)}\right)'=\Psi(\frac{1}{2}-z)-\Psi(1-z)<0$ for all $z\in(0,\frac{1}{2})$ since thedigamma function
\begin{equation}\label{eqc3}
 \Psi(z)\!:=\frac{\Gamma'(z)}{\Gamma(z)}=-\gamma+\int_0^1\frac{1-x^{z-1}}{1-x}dx
\end{equation}
is a well defined monotonous increasing function for $z>0$, see~\cite[\S 6.3]{stegun}. Here $\gamma\approx 0.577$ is the Euler-Mascheroni constant. This proves the validity of the inequality in $(a)$.

Let us turn next to the proof of the assertion with regard to $P(\al,1).$ To this end, for the sake of convenience, we introduce the function
\begin{equation}\label{eqc2}
 F(\al)\!:=\frac{3}{16}\frac{((\al-1)(\al-3)(\al-5))^2}{(2\al-3)(4\al^2-1)}\frac{\mathrm B\!\left(1-\al,-\frac{3}{2}\right)}{\mathrm B\!\left(1-\frac{\al}{2},-\frac{1}{2}\right)}=
 {\frac {\Gamma^2\left( \frac{7}{2}-\frac{\al}{2} \right)}{
{2}^{\al}\sqrt {\pi }\Gamma  \left( \frac{5}{2}-\al \right) }},
\end{equation}
where the identity follows using the so called duplication formula for $\Gamma$, see \cite[\S 6.1.18]{stegun}. This function will enable us to write $P(\al,1)=\rho_0(\al)-\rho_1(\al)+\rho_2(\al)-\rho_3(\al)$ in a more convenient form taking advantage of the fact that each $\rho_n(\al)$ is linear in $\mathrm B\!\left(1-\al,-\frac{3}{2}\right)$ and $\mathrm B\!\left(1-\frac{\al}{2},-\frac{1}{2}\right)$, see~\refc{eqc1}. In doing so, some easy computations using a symbolic manipulator (see \cite{Maple} for instance) show that 
\begin{align*}
 &P(\al,1)=\frac{40\al(\al+2)(3\al-5)}{((\al-1)(\al-3)(\al-5))^2}\mathrm B\!\left(1-\frac{\al}{2},-\frac{1}{2}\right)\Big(F(\al)-g(\al)\Big),\\
\intertext{where}
 &g(\al)\!:=\frac{(23\al-94)(\al-1)(\al-3)(\al-4)(\al-5)}{160(3\al-5)(\al+2)}.
\end{align*} 
Thus, since $\mathrm B\!\left(1-\frac{\al}{2},-\frac{1}{2}\right)=-\frac{2\sqrt{\pi}\Gamma(1-\frac{\al}{2})}{\Gamma(\frac{1}{2}-\frac{\al}{2})}$ is negative for all $\al\in (0,1)$, the assertion in $(b)$ will follow once we prove that $F(\al)>g(\al)$ for all $\al\in (0,1).$ As an intermediate step to this end we claim that if $\al\in (0,1)$ then $F(\al)>0$, $F'(\al)<0$ and $F''(\al)>0.$ The first inequality is clear from \refc{eqc2} because $\Gamma(z)>0$ for all $z>0.$ The second inequality is also easy because some computations show that
 \[
 \frac{F'(\al)}{F(\al)}=-\big(h(\al)+\ln 2\big)\text{ with $h(\al)\!:=\Psi(7/2-\al/2)-\Psi(5/2-\al)$}
 \]
and, on the other hand, $h(\al)>0$ for all $\al\in (0,1)$ due to $\Psi'(z)>0$ for $z>0.$ Finally, in order to show the third inequality we first note that
\[
  \frac{F''(\al)}{F(\al)}=\big(h(\al)+\log 2\big)^2-h'(\al).
 \]
Furthermore, due to $x^{\frac{3}{2}-\al}>x^{\frac{5}{2}-\frac{\al}{2}}>2^{-n}x^{\frac{5}{2}-\frac{\al}{2}}$ for all $x\in (0,1)$ and $\al>0,$ from \refc{eqc3} it turns out that
\begin{align*}
 h^{(n)}(\al)&=(-1/2)^n\Psi^{(n)}(7/2-\al/2)-(-1)^n\Psi^{(n)}(5/2-\al)\\[5pt]
 &=\int_0^1\left(x^{\frac{3}{2}-\al}-2^{-n}x^{\frac{5}{2}-\frac{\al}{2}}\right)\frac{(-\log x)^n}{1-x}dx>0
 \text{ for all $n\in\Z_{\geq 0}.$}
\end{align*}
Hence $h^{(n)}$ is increasing on $(0,+\infty)$ for all $n\in\Z_{\ge 0}$ and, consequently, $h^{(n)}(1)> h^{(n)}(\al)> h^{(n)}(0)$ for all $\al\in (0,1)$. Thus if $\al\in (0,1)$ then
\[
 \frac{F''(\al)}{F(\al)}>\big(h(0)+\log 2\big)^2-h'(1)= \left(\frac{2}{5}+\log 2\right)^2+\frac{27}{8}-\frac{5\pi^2}{12}\approx 0.46.
 \]
Accordingly, $F''(\al)>0$ for all $\al\in (0,1)$ and this concludes the proof of the claim. We proceed now with the proof of $(b)$, which let us recall that it will follow once we prove that $F(\al)>g(\al)$ for all $\al\in (0,1).$ With this aim in view we take the tangent lines to the graph of $F$ at $\al=0$ and $\al=1,$ that are given by 
\[\textstyle
 \ell_0(\al)=\frac{75}{16}-(\frac{15}{8}+\frac{75}{16}\log 2)\al
 \text{ and }
 \ell_1(\al)=\frac{4}{\pi}+\frac{2}{\pi}(1-6\log 2)(\al-1),
\]
respectively, see \figc{tangentes}.
\begin{figure}[t]
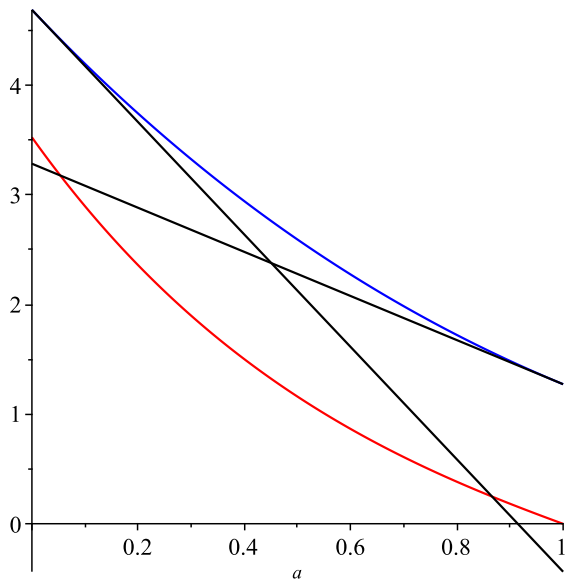

 \centering
 \begin{lpic}[l(0mm),r(0mm),t(0mm),b(0mm)]{dib_ap_c_2(0.4,0.4)}
 \end{lpic}
 \caption{The graphs of the transcendental function $F(\al)$ in blue, its tangent lines $\ell_0(\al)$ and $ \ell_1(\al)$ at $\al=0$ and $\al=1$, respectively, in black and the rational function $g(\al)$ in red.}\label{tangentes}
 \end{figure}
Since $F$ is convex, in order to show that $F(\al)>g(\al)$ for all $\al\in (0,1),$ it suffices to verify that $\max\{\ell_0(\al),\ell_1(\al)\}>g(\al)$ for all $\al\in (0,1).$ To see this we consider the unique solution of $\ell_0(\al)=\ell_1(\al)$, which one can check that it is given by
\[ 
 \al=\hat\al\!:={\frac { 75 \pi-32-192 \log2}{ \left( 75 \pi-192
 \right) \log2 +30 \pi+32}}\approx 0.45.
\]
One can also verify that, for $i=0,1,$ $\ell_i(\al)-g(\al)=\frac{p_i(\al)}{160\pi^i(2+\al)(5-3\al)}$ with
\begin{align*}
 p_0(\al)\!=&\,23{\al}^{5}-393{\al}^{4}+ \left( 3479+2250\log 2 
 \right) {\al}^{3}+ \left( -9957+750\log 2  \right) {\al}^{2}\\
  &+\left( 7688-7500\log 2  \right) \al+1860\\
  \intertext{and}
 p_1(\al)\!=&\,23\pi {\al}^{5}-393\pi {\al}^{4}+ \left( -960+2579\pi +5760
\log 2  \right) {\al}^{3}- \left(  8007\pi +1280+3840\log 2 \right) {\al}^{2}\\
 &+\left( 11438\pi +2880-
21120\log 2  \right) \al+3200-5640\pi +19200\log 2. 
\end{align*}
By applying Sturm's Theorem we can assert that $p_0$ is positive on $(0,0.46)$ and that $p_1$ is positive on $(0.44,1)$, which imply that $\max\{\ell_0(\al),\ell_1(\al)\}>g(\al)$ for all $\al\in (0,1)$ as desired. This proves $(b)$. 

Our last task is to prove the assertion in $(c)$. To this end we use a symbolic manipulator in order to show that
 \[
 \text{Discrim}_{\be}P(\al,\be)=\frac{-2(\al+2)\mathrm B(1-\frac{\al}{2},-\frac{1}{2})^4}{3((\al-1)(\al-3)(\al-5))^8}\mathscr R\big(\al,F(\al)\big),
 \]
where
\begin{align*}
\mathscr R(\al,t)=&-16384  ( 2 \al-1 )  ( 8 {\al}^{6}+36 {
\al}^{5}-126 {\al}^{4}-413 {\al}^{3}+429 {\al}^{2}+576
 \al-512 )  ( \al+3 ) ^{2} ( 2 \al-3
 ) ^{2}{t}^{4}\\
 &+3072  ( \al-1 )  ( \al-3
 )  ( \al-5 )  ( 2 \al-3 )  ( 
\al+3 ) \big( 48 {\al}^{8}+252 {\al}^{7}-1904 {
\al}^{6}-2305 {\al}^{5}+11568 {\al}^{4}\\
&-2566 {\al}^{3}-
14160 {\al}^{2}-2784 \al+11520\big) {t}^{3}-24 \al 
 ( \al+2 )  \big( 768 {\al}^{9}-7808 {\al}^{8}+
3616 {\al}^{7}\\
&+135520 {\al}^{6}-221032 {\al}^{5}-557976 {
\al}^{4}+823685 {\al}^{3}+1082256 {\al}^{2}-894960 \al
\\
&-915840 )  ( \al-5 ) ^{2} ( \al-1 ) ^
{2} ( \al-3 ) ^{2}{t}^{2}-4  ( \al-4 ) 
 ( \al+2 )  ( 320 {\al}^{8}-1400 {\al}^{7}+
1830 {\al}^{6}-5491 {\al}^{5}\\
&+4678 {\al}^{4}+32889 {
\al}^{3}-4482 {\al}^{2}-47520 \al-64800\big)  ( 
\al-1 ) ^{3} ( \al-3 ) ^{3} ( \al-5
 ) ^{4}t\\
 &+15 \al  ( \al-2 )  ( \al-4
 )  ( \al+2 )  ( 5 {\al}^{4}-15 {\al}
^{3}-5 {\al}^{2}+27 \al+36 )  ( \al-1 ) ^{
4} ( \al-3 ) ^{5} ( \al-5 ) ^{6}.
\end{align*}
Let us mention that in order to ease this computation and introduce $F$ we use 
that the coefficients of $P(\al,\be)=\rho_0(\al)-\rho_1(\al)\be+\rho_2(\al)\be^2-\rho_3(\al)\be^3$, see \refc{eqc1}, are linear in $\mathrm B\!\left(1-\al,-\frac{3}{2}\right)$ and $\mathrm B\!\left(1-\frac{\al}{2},-\frac{1}{2}\right)$ and that, on the other hand, the discriminant of a third degree polynomial is a homogeneous polynomial of degree~4 in its coefficients. On account of the above expression it is clear that~$(c)$ will follow once we prove that $\mathscr R(\al,F(\al))\neq 0$ for all $\al\in (0,1).$ To this end we note that $F(0)=\frac{75}{16}$ and $F(1)=\frac{4}{\pi}.$ Therefore, see \figc{tangentes}, the graph $t=F(\al)$ for $\al\in (0,1)$ verifies $\max\{\ell_0(\al),\frac{4}{\pi}\}<F(\al)<\frac{75}{16}$ because we previously proved that $F'<0$ and $F''>0$ on the interval $(0,1)$. Accordingly it suffices to show that $\mathscr R(\al,t)\neq 0$ for all $(\al,t)$ inside the trapezium given by $\max\{\ell_0(\al),\frac{4}{\pi}\}<t<\frac{75}{16}$ and $\al\in (0,1).$
We will prove this taking $t\in (\frac{4}{\pi},\frac{75}{16})$ as a fixed parameter and showing that the polynomial $\al\mapsto\mathscr R(\al,t)$ has not any root on $(\ell_0^{-1}(t),1),$ where $\ell_0^{-1}(t)=\frac{75-16t}{30-75\log 2}.$ 
To this effect we show first that $\mathscr R\left(\frac{75-16t}{30-75\log 2},t\right)$, $\mathscr R(1,t)$ and $\text{Discrim}_{\al}\mathscr R(\al,t)$ do not vanish for all $t\in (\frac{4}{\pi},\frac{75}{16})$. This implies that the number of roots of $\mathscr R(\al,t)=0$ on $(\ell_0^{-1}(t),1)$ does not change for $t\in (\frac{4}{\pi},\frac{75}{16})$. Taking this into account the desired result follows by checking, for instance, that this number is zero for $t=2\in (\frac{4}{\pi},\frac{75}{16})$. All these assertions can be checked systematically by applying Sturm's Theorem because only one variable polynomials are involved. This shows the validity of $(c)$ and so the inequality in $(ii)$ is
true. This concludes the proof of the result.
\end{prova}

\bibliographystyle{plain}

\end{document}